\documentclass{article}

\usepackage{amsmath}
\usepackage[francais]{babel}
\usepackage{amssymb}
\usepackage[applemac]{inputenc}
\usepackage{url}
\newcommand{\sm}{\raisebox{2pt}{~\rule{6pt}{1.2pt}~}}

\def\M{{\bf M}}

\def\F{{\bf F}}
\def\Z{{\bf Z}}
\def\R{{\bf R}}
\def\G{{\bf G}}
\def\r{{\frak r}}
\def\sch{{\rm sch}}

\def\dim{\mathop{\rm dim}}
\def\O{\Omega}
\def\o{\omega}
\def\s{\sigma}
\def\t{\mbox{\boldmath $t$}}

\def\det{\mathop{\rm det}}

\def\GL{\mathop{\bf GL}}
\def\SL{\mathop{\bf SL}}

\def\Sp{\mathop{\bf Sp}}
\def\SO{\mathop{\bf SO}}

\def\PSL{\mathop{\bf PSL}}
\def\End{\mathop{\rm End}}
\def\Gal{\mathop{\rm Gal}}
\def\Lie{\mathop{\rm Lie}}
\def\car{\mathop{\rm car}}
\def\m{\boldsymbol \mu_2}
\def\Aut{\mathop{\rm Aut}}
\def\Hom{\mathop{\rm Hom}}
\def\Ind{\mathop{\rm Ind}}

\def\IM{\mathop{\rm Im}}
\def\Tr{\mathop{\rm Tr}}

\def\Ker{\mathop{\rm Ker}}

\def\l{\ell}
\def\e{\varepsilon}
\def\n{\noindent}
\def\kbar{\overline{k}}

\def\Fbar{\overline{\F}}

\begin{document}

\qquad \qquad \qquad \qquad \qquad \qquad  à paraître dans {\it Transformation Groups} {\bf 19}  (2014)

\bigskip

\bigskip
     
\centerline {\bf  Bases normales autoduales et groupes unitaires en caractéristique 2}  

\bigskip
     
     \centerline {\rm Jean-Pierre Serre}
         
      \smallskip

      \centerline{\it \`A E.B. Dynkin pour son 90-ième anniversaire}
      
      \bigskip

    \centerline {\bf   Introduction}
      
     \medskip
     
     L'origine du présent travail est un théorème de Bayer-Lenstra  dont voici l'énoncé:  
     
     \smallskip
     
  \n {\bf Théorème A} ([BL 90, th.6.1]) - {\it Soit  $k$  un corps de caractéristique $2,$ et soit $L/k$ une extension galoisienne finie de groupe de Galois $G$. Supposons que $G$ soit commutatif. Pour que $L/k$ possède une base normale autoduale} (cf.
   ci-dessous), {\it il faut et il suffit que $G$ n'ait pas d'élément d'ordre $4$.}
   
 {\small  [Lorsque la caractéristique de $k$ est $\neq 2$, le même énoncé est valable, à condition d'y remplacer ``d'ordre 4'' par ``d'ordre 2" - ce qui revient à dire que l'ordre de $G$ doit être impair.] }
   
   \smallskip 
   
\n Rappelons (cf. [BL 90]) qu'une base normale autoduale (en abrégé : ``BNA") de $L/k$ est une $k$-base $B$ de $L$ qui est stable par l'action de $G$ et qui est telle que

\smallskip

$ {\rm Tr}_{L/k}(xy) =\begin{cases}
0 \quad  {\rm si} \ \ x,y \in B, \ x \neq y\\
1 \quad {\rm si} \ \  x,y \in B, \  x=y. 
\end{cases}$

  \bigskip
   
    Il est naturel d'essayer d'étendre le théorème A au cas où $G$ n'est pas commutatif. C'est ce que nous allons faire. Le résultat est le suivant (il est annoncé dans  [Se 05]):

    \smallskip
    
      \n {\bf Théorème B} - {\it Soit  $k$  un corps de caractéristique $2,$ et soit $L/k$ une extension galoisienne finie de groupe de Galois $G$. Pour que $L$ possède une BNA, il faut et il suffit que $G$ soit
     engendré par des éléments d'ordre impair et par des éléments d'ordre} 2.

   \n [Noter le cas particulier où $G$ est d'ordre impair, déjà traité dans [Ba 89].]  
   
   \smallskip
   
   Ce qui est surprenant dans cet énoncé, c'est que l'existence d'une BNA ne dépende  que de la structure du groupe de Galois $G$, et
   pas des propriétés du corps de base $k$, ni de celles de l'extension $L$; ce n'est pas ce qui se passe quand $\car(k) \neq 2$, comme le montrent de nombreux exemples,
   cf. [Ba 89], [BFS 94], etc.

    \smallskip
    Nous démontrerons le théorème B au \S6 (corollaire 6.1.10), comme conséquence d'un critère d'isomorphisme pour les $G$-formes trace de deux $G$-algèbres galoisiennes (théorème 6.1.5). La démonstration se fait en traduisant la question en termes de cohomologie galoisienne, comme dans [BFS 94].  Le groupe algébrique qui intervient est le groupe unitaire de l'algèbre à involution $k[G]$; cela nous amènera à étudier la structure des groupes unitaires des algèbres à involution en caractéristique 2, ce qui ne semble pas avoir été fait jusqu'à présent en dehors de cas très particuliers. De façon plus précise, soit $k$ un corps parfait de caractéristique 2, soit $A$ une $k$-algèbre à involution de dimension finie, soit  $U_A $ le groupe unitaire correspondant (vu comme groupe algébrique sur  $k$), et soit $U_A^0$ la composante neutre de $U_A$.  Le théorème de cohomologie galoisienne dont nous aurons besoin est:
         
        \medskip
     
 \n {\bf Théorème C} - {\it  On a $H^1(k,U_A^0) = 0$, autrement dit tout $U_A^0$-torseur a un point rationnel.}     
    
         \medskip
         
         Pour prouver ce résultat, on peut supposer que $U_A^0$ engendre l'algèbre $A$. Faisons cette hypothèse.
         Le cas où $A$ est semi-simple est facile. On décompose $A$ en produits d'algèbres simples, et l'on constate que
       $U_A$ est connexe et que, après extension des scalaires, c'est un produit de groupes linéaires et de groupes symplectiques. La trivialité de $H^1(k,U_A^0) = H^1(k,U_A)$ en résulte, grâce à l'hypothèse que $k$ est parfait de caractéristique 2, cf. \S4.7.
       
       \smallskip
            
       Pour passer de là au cas général, il est naturel d'introduire le radical $\r$ de $A$ et de comparer
     les groupes $U_A^0$ et $U_{A/\r}$ grâce à la projection 
      $ \pi: U_A^0 \to \ U_{A/\r}\ .$
     
    \n Le noyau de $\pi$ est un groupe unipotent connexe qui est négligeable du point de vue de la cohomologie galoisienne. La vraie difficulté est que $\pi$  n'est pas surjectif en général, contrairement à ce qui se passerait si la caractéristique de $k$ était $\neq 2$. Toutefois on peut démontrer que $\pi$ est ``presque surjectif'': son image est un sous-groupe réductif de rang maximum de $U_{A/\r}$ qui se décompose en produit comme $U_{A/\r}$; la seule différence est que certains facteurs  symplectiques $\Sp_{2n}$ de $U_{A/\r}$ sont
     remplacés par des groupes orthogonaux $\SO_{2n}$. La démonstration de ce fait occupe le \S2 et le \S3; elle repose sur une propriété
     très particulière des poids des tores maximaux des groupes unitaires: la ``propriété PL'', décrite au \S1, qui vaut, plus généralement, pour les groupes définis comme fixateurs de tenseurs de degré $\leqslant 2$, cf. \S1.4.

     Une fois ce résultat obtenu, la nullité de $H^1(k,U^0_A)$ se démontre par des arguments cohomologiques standard, cf. \S4. On peut alors passer au cas dont nous avons besoin: $A = k[G]$, cf. \S5. Ici il n'est plus nécessaire de supposer que $k$ est parfait; l'un des  principaux résultats est une caractérisation de $G \cap U^0_A$ comme le sous-groupe de $G$ engendré par les éléments d'ordre 2 et par les carrés (théorème 5.3.1). On
déduit de là le critère d'isomorphisme de $G$-formes trace mentionné plus haut; c'est l'objet du \S6, qui en donne diverses applications, par exemple au principe de Hasse, lorsque  $k$  est un corps global (théorème 6.1.18). Le \S7 contient des compléments variés, et mentionne plusieurs questions ouvertes, notamment la suivante: si $G$ est un $2$-groupe, est-il vrai que $G \to U_G/U_G^0$ est surjectif ?

     \bigskip

         \n {\it Notations.} 
         
         \smallskip
         
         Si $A$ est un anneau, on note $A^\times$ son groupe multiplicatif.
         
          On note $X \sqcup Y$  la réunion de deux ensembles disjoints $X$ et $Y$.

          Si $X$ est un schéma sur un corps $k$, et si $K$ est une extension de $k$, on note    $X_{/K}$            
               le schéma que l'on déduit de $X$ par extension des scalaires à $K$, et l'on note $X(K)$ l'ensemble des $K$-points de $X$.    
               
               Si $X$ est un schéma, on note $X^{\rm red}$ le schéma réduit associé.

          Le terme de ``groupe algébrique'' (sur un corps $k$) est utilisé au sens de ``schéma en groupes affine de type fini'', autrement dit ``sous-schéma en groupes fermé de $\GL_n$ pour  $n$  assez grand'', cf. par exemple [KMRT 98, chap.VI] et [Wa 79]; le faisceau des anneaux locaux peut avoir des éléments nilpotents $\neq 0$.  En fait, le cas le plus important est celui où le groupe est lisse, mais les éléments nilpotents interviennent parfois (par exemple pour le groupe $  \m$ des racines carrées de l'unité, ou aussi pour la définition du noyau d'un homomorphisme).
          
            Si $G$  est un groupe algébrique, sa composante neutre est notée  $G^0$, et son algèbre de Lie est notée $\Lie(G)$.
          
                      \smallskip  
            \n         Les autres notations sont standard.
         
               \medskip

      \medskip

 \centerline    {\bf \S1 - Les tores de type PL}
     
     \smallskip
     \n 1.1. {\bf La propriété PL}.
     
     \smallskip
     
     Soit  $L$ un $\Z$-module libre de type fini. Soit $\O$
     une partie finie de $L$. Nous dirons que $\O$ est ``presque libre''  (en abrégé `` PL '') s'il existe une partie $\o$ de $\O$ telle que :
       
       (1.1.1)  $\omega$ {\it est libre au sens usuel de ce terme $:$ les éléments de $\o$ sont $\Z$-linéairement indépendants}.
       
       (1.1.2) {\it On a} \ $\O \ \subset \ \{0\} \cup \o \cup -\o.$
       
         \smallskip
         
         \n {\small  [Dans cette formule, $-\o $ désigne l'ensemble des $-\alpha$ pour $\alpha \in \o$. Nous utiliserons souvent ce genre de notation dans la suite.]}
         
         \smallskip
         
         [
       \n    Autre caractérisation:
           
             \smallskip
       
       (1.1.3) {\it Toute partie $\varpi$ de $\O$ telle que $\varpi \cap  -\varpi \ = \ \varnothing$ est libre}.
       
       \medskip
       
         {\it Exemple.} Si $\{x_1,x_2,x_3\}$ est une partie libre de $L$, l'ensemble $ \{0,x_1,x_2,x_3,-x_1\}$ est presque libre.
       
       \medskip

       Cette notion a les propriétés de stabilité suivantes, qui sont immédiates:
       
       \smallskip
       (1.1.4)  {\it Si  $L'$ est un sous-groupe de  $L$ et si $\O \subset L' $, alors
     $\O$ est presque libre dans $L'$ si et seulement si il l'est dans $L$.}
     
     \smallskip
         (1.1.5) {\it Si $\O' \subset \O$ et si $\O$ est presque libre, alors $\O'$ est presque libre.}

             \medskip
             
             \n 1.2. {\bf Tores de type PL}.
             
             \smallskip
  
             Soit  $k$  un corps algébriquement clos (de caractéristique
             quelconque) et soit $V$ un espace vectoriel de dimension finie sur  $k$. On note $\GL_V$ le groupe des automorphismes de $V$, vu comme $k$-groupe algébrique.         
             
                 Soit  $T$ un $k$-tore opérant sur $V$, autrement dit muni d'un homomorphisme $\varphi : T \to \GL_V$. Soit $X(T) = \Hom(T,\G_m)$ le groupe des caractères
             de $T$; comme d'habitude, on écrit $X$ additivement, ce qui amène à noter  $t^\chi$ l'image d'un point $t$ de $T$ par le caractère $\chi$. L'espace $V$ se décompose en $V = \oplus_{\chi \in X} V_\chi$, où $V_\chi$ est le sous-espace propre relatif à $\chi$, i.e. l'ensemble des $v \in V$ tels que $\varphi(t)v = t^\chi v$ pour tout $t$. On dit que $\chi$ est un {\it poids de} $V$ si $ V_\chi \neq 0$.
             
               \medskip
  \n {\it Définition}. On dit que l'action de $T$ sur $V$ est de type PL si l'ensemble $\O$ de ses poids est de type PL au sens du \S1.1.
              
              \medskip
              
              \n {\it Remarque.}
   Si l'action de $T$ sur $V$ est de type PL, et si $W$ est un sous-espace vectoriel de $V$ stable par  $T$, alors l'action de $T$ sur $W$ est de type PL, et il en est de même de son action sur $V/W$. Cela résulte de (1.1.5).
              Dans la suite, $T$ sera le plus souvent un sous-tore de $\GL_V$, et $\varphi$
              sera l'injection canonique $T \to \GL_V$. Dans ce cas, on dira 
              simplement que $T$  est un sous-tore de $\GL_V$ de type PL.
              En fait, le cas général se ramène à ce cas particulier: en effet,
              si l'on a $\varphi : T \to \GL_V$ et si  $T'$ désigne l'image de $T$
              dans $\GL_V$, alors $(T,\varphi)$ est de type PL si et seulement si
              $T'$ est de type PL: cela résulte de (1.1.4), puisque $X(T')$ est un sous-groupe de $X(T)$.

              \bigskip
              
              \n 1.3. {\bf Premiers exemples de tores de type PL.}
              
              \smallskip
              
              1.3.1. Un tore maximal de $\GL_n$ : l'ensemble de ses poids est libre à $n$ éléments.
              
                   \smallskip

              1.3.2. Un tore maximal de $\Sp_{2n}$, ou de $\SO_{2n}$: l'ensemble de ses poids est de la forme $\o \cup -\o$, où  $\o$ est libre 
              à $n$ éléments.

                   \smallskip

              1.3.3. Un tore maximal de $\SO_{2n+1}$: l'ensemble de ses poids
              est de la forme $\{0\} \cup \o \cup -\o$, où  $\o$ est libre 
              à $n$ éléments.

                            \medskip
                            
                            Rappelons que, lorsque la caractéristique est égale à 2, le groupe orthogonal ${\bf O}_{2n}$ est lisse; le groupe $\SO_{2n}$ est sa composante neutre; c'est le noyau de l'invariant de Dickson ${\bf O}_{2n} \to \Z/2\Z$, cf. [KMRT 98, \S12.12]. Par contre, le groupe orthogonal  ${\bf O}_{2n+1}$ n'est pas lisse, comme 
                            le montre déjà le cas $n=0$, où c'est $\boldsymbol \mu_2$; le groupe $\SO_{2n+1}$ est défini comme le noyau de $\det : {\bf O}_{2n+1} \to   \G_m $ ; c'est un groupe lisse, et l'on a ${\bf O}_{2n+1} = \boldsymbol \mu_2 \times \SO_{2n+1}$.
                       
               \smallskip
              1.3.4. Si  $A$  est une algèbre à involution,
              et  $U$ le groupe unitaire correspondant (vu comme sous-groupe
              de  $\GL_A$ par l'action par les translations à gauche), alors tout tore maximal de $U$ est de type PL, cf. proposition 3.2.5.
             \bigskip
             
             \medskip
             
             \n 1.4. {\bf Exemple $:$ tores maximaux des fixateurs de tenseurs de degré $\leqslant 2$.}
             
             \smallskip
               Les différents exemples du \S1.3 sont des cas particuliers du suivant:
               
               Soit $V$  un espace vectoriel de dimension finie sur $k$ (supposé algébriquement clos, pour simplifier). Soit $(\theta_i)_{i \in I}$ une famille d'éléments de l'un des huit types suivants:
               
               \smallskip
               
               (1.4.1) \ \ un élément de $V$;
               
               (1.4.2) \ \ un élément du dual  $V'$ de $V$;
               
               (1.4.3) \ \  un élément de $\otimes^2 V$;
               
               (1.4.4) \ \  un élément de $\otimes^2 V'$;
               
               (1.4.5) \ \ un élément de $V \otimes V' = \End (V)$;
               
               (1.4.6) \ \ une forme quadratique sur $V$;
               
               (1.4.7) \ \ une forme quadratique sur $V'$;
               
               (1.4.8) \ \ un sous-espace vectoriel $W$ de $V$.
               
               \smallskip
               
               Soit $G \subset \GL_V$  le fixateur des $\theta_i$. 
               
               \n  {\small   [Précisons ce que cela signifie: si $k'$ est une $k$-algèbre commutative, un point  $g \in G(k')$ est un automorphisme de $k' \otimes V$ qui fixe les tenseurs déduits des $\theta_i$ par extension des scalaires de $k$ à $k'$ (``fixer"  a un sens clair pour les sept premiers cas, et pour (1.4.8), cela signifie que $g$ laisse stable  $k' \otimes_k W$).]}
               
               \medskip
               
               \n {\bf Proposition 1.4.9} - {\it Les tores maximaux de $G$ sont de type PL.}
               
               \smallskip

               \n {\it Démonstration}. Nous allons prouver un résultat plus précis:
               
             \smallskip  
               
               \n {\bf Proposition 1.4.10} - {\it Soit $M$ un sous-groupe de $G$ de type multiplicatif } (cf. [SGA 3, II, exposés VIII et IX]), {\it dont le groupe des caractères ne contient aucun élément d'ordre $2$. \footnote{Lorsque $k$ est de caractéristique 2 - ce qui est le cas intéressant pour la suite - cette hypothèse équivaut à dire que $M$ est {\it lisse}.} Il existe un sous-tore $\widetilde{M}$ de $G$ qui contient $M$ et qui est de type PL.}
               
               [Cet énoncé entraîne la proposition précédente, car, si on l'applique à un tore maximal  $T$  de  $G$, on obtient un tore  $\widetilde{T}$ de  $G$  qui est de type PL et qui contient $T$, donc qui est égal à $T$.]

               \smallskip
               
               \n {\it Démonstration.} Soit $X = \Hom(M,\G_m)$ le groupe des caractères de $M$, que nous écrirons additivement. L'action de $M$ sur $V$ décompose $V$ en somme directe  $V = \sum_{\chi \in X} V_\chi$ de sous-espaces propres.
           Soit $\widetilde{M}$  le groupe des automorphismes
            de $V$ qui, sur chaque $V_\chi$, sont égaux à une homothétie  $t_\chi$, satisfaisant aux deux conditions:  
            
            \smallskip
            
            (1.4.11)  $t_\chi = 1$ \ si \ $\chi = 0$;

            \smallskip
            (1.4.12) $t_\chi t_{-\chi} = 1$  \ pour tout  $\chi.$
            
            \smallskip
            
       \n     Ce groupe contient $M$. Nous allons voir qu'il convient.  Autrement dit:

            \smallskip
            
            \n {\bf Lemme 1.4.13} - {\it Le groupe $\widetilde{M}$ est un sous-tore de $G$ de type PL.}
            
            \smallskip

            \n {\it Démonstration du fait que $\widetilde{M}$ est un tore de type PL.} Soit $\O$ l'ensemble des $
            \chi \neq 0$ tels que $V_\chi \neq 0$; c'est un ensemble fini; choisissons une partie $\o$ de $\O$
            qui ne rencontre pas $-\o$, et qui est maximale pour cette propriété. Un point de $\widetilde{M}$ est déterminé par ses coordonnées $t_\chi$ avec $\chi \in \o$, et celles-ci peuvent être choisies arbitrairement. Cela montre que $\widetilde{M}$
            est un tore dont le groupe des caractères a pour base les éléments de $\o$ (vus comme caractères de $\widetilde{M}$, et non plus de $M$); de plus, l'ensemble des poids de ce tore est contenu dans $ \{0\} \cup \o \cup -\o$, donc est presque libre.
            
            \smallskip

            \n {\it Démonstration du fait que  $\widetilde{M}$ est contenu dans $G$.} On doit montrer que les $\theta_i$ sont invariants par tout point $(t_\chi)$ de $\widetilde{M}$. Il y a huit cas à considérer, suivant que $\theta_i$ est de type (1.4.1), (1.4.2), ..., (1.4.8):
            
            \smallskip
            
            Type (1.4.1). On a  $\theta_i \in V$; comme $\theta_i$ est invariant par $M$, il appartient à  $V_0$, et il est invariant par $\widetilde{M}$ puisque les éléments de $\widetilde{M}$ fixent $V_0$ d'après (1.4.11).
            
            Type (1.4.2). On a  $\theta_i \in V'$: même démonstration, avec $V$ remplacé par $V'$.
            
            Type (1.4.3). On a  $\theta_i \in V \otimes V = \oplus_{\chi_1,\chi_2} V_{\chi_1} \otimes V_{\chi_2}  $. Puisque  $\theta_i $ est invariant par $M$, il est contenu dans la somme directe des $V_\chi \otimes V_{-\chi}$; or ceux-ci sont fixés
            par  $\widetilde{M}$, d'après (1.4.12).
            
            Type (1.4.4). Même démonstration que celle du type (1.4.3), avec $V$ remplacé par $V'$.

            Type (1.4.5). On a $\theta_i \in \End(V)$. Dire que $\theta_i$ est invariant par $M$ signifie qu'il commute à l'action des éléments de $M$, ou encore que les $V_\chi$ sont stables par $\theta_i$; cela entraîne que $\theta_i$ commute à l'action de $\widetilde{M}$.
            
            Type (1.4.6). Dans ce cas, $\theta_i$ est une forme quadratique. Soit $b$ la forme bilinéaire associée. Soit 
            $t = (t_\chi)$ un $k$-point de $\widetilde{M}$. D'après le cas (1.4.4), $t$ fixe $b$ (en effet $b$  est fixée par $G$). Il en résulte que la forme bilinéaire associée à la forme quadratique  $t(\theta_i) - \theta_i$ est $0$; autrement dit, c'est le carré d'une forme linéaire $\l$ sur $V$. On doit montrer que $\l = 0$, et il suffit de le faire sur chaque $V_\chi$; c'est clair
            pour $\chi = 0$ puisque $t$ fixe $V_0$; c'est non moins clair pour $V_\chi, \ \chi \neq 0$, car la restriction de $\theta_i$ à $V_\chi$ est 0  (puisque $M$ agit sur $V_\chi$ par des homothéties de rapport arbitraire). 
            
            Type (1.4.7). Même démonstration que celle de (1.4.6), avec $V$ remplacé par $V'$.
            
            Type (1.4.8). Dans ce cas $\theta_i$ est un sous-espace vectoriel $W$ de $V$. Dire qu'il est stable par $M$ équivaut à dire qu'il est somme directe des $W \cap V_\chi$, et il est donc stable par $\widetilde{M}$ puisque ce groupe opère par des homothéties sur chaque   $V_\chi$.

            \bigskip
            
              \n {\it Remarque.} Supposons que ${\rm car}(k)=2$. La proposition 1.4.10 entraîne que, si $p$ est un nombre premier $\neq 2$, tout $p$-sous-groupe fini commutatif de $G$ est contenu dans un tore maximal (et en particulier est contenu dans $G^0$). Cela donne des renseignements non triviaux sur la structure de $G$; par exemple, d'après [St 75, th.2.27 et th.2.28], 
cela montre qu'aucun quotient de $G^{0 \ {\rm red}}$ n'est un groupe simple de type $\sf F_4,
\sf E_6, \sf E_7$ ou $\sf E_8$ (ce qui  résulte aussi des résultats du \S2.6, qui éliminent également le type $ \sf G_2$).   

\bigskip

 {\small 
\n {\it Compléments.} Mentionnons brièvement quelques autres propriétés de $G$:

 \smallskip
(1.4.14) {\it Si la caractéristique de $k$ est $\neq 2$, le groupe $G$ est lisse}. [On vérifie que, si $X \in \Lie(G)$, alors $(1+tX)(1-tX)^{-1}$ appartient à $G(k)$ pour tout  $t$ tel que $1+tX$ et $1-tX$ soient inversibles; on en déduit un morphisme d'un ouvert de la droite affine dans $G$; comme la droite affine est un schéma réduit,  l'image de ce morphisme est contenue dans $G^{\rm red}$; d'où, en dérivant en $t=0$, le fait que $2X$ appartient à $ \Lie(G^{\rm red})$; on a donc $\Lie(G^{\rm red})=\Lie(G)$, ce qui entraîne la lissité de $G$.]

\smallskip

(1.4.15) {\it Si $x \in G(k)$, alors} $x^2 \in G^0(k)$. [Utiliser le fait que $(t+x)(t+x^{-1})^{-1}$ appartient à $G(k)$ pour tout  $t$ d'un ouvert de la droite affine contenant $0$, et en déduire qu'il appartient à $G^0(k)$.]

\smallskip

(1.4.16) {\it Si $k$ est de caractéristique $2$, et si aucun des tenseurs $\theta_i$ n'est de type $(1.4.6)$ ou $(1.4.7)$, alors tout  $x \in G(k)$ d'ordre $2$ appartient à} $G^0(k)$. [Utiliser le fait que $t \mapsto 1+t(x+1)$ est un homomorphisme du groupe additif $\G_a$ dans $G$, donc dans $G^0$.]  

\smallskip

Nous reviendrons au \S3.3 sur ces deux dernières  propriétés dans le cas particulier où  $G$  est un groupe unitaire.}

             \newpage
 \centerline    {\bf \S2 - Structure des groupes linéaires à tore maximal de type PL}
                            
              \medskip
                        
         \n 2.1. {\bf Enoncé du théorème, dans le cas irréductible.}

              \medskip
           Soit $k$ un corps algébriquement clos de caractéristique 2, soit $V$ un $k$-espace vectoriel
   de dimension finie, et soit $G$ un sous-groupe algébrique {\it connexe lisse} de $\GL_V$ ayant la propriété suivante:
   
            \medskip
              \qquad  {\it  Les tores maximaux de $G$ sont de type {\rm PL} au sens du \S$1.2.$}
               
               \medskip
               
               Nous dirons alors que $(G,V)$ {\it est de type} PL.
       
       \medskip        
               
               On se propose de décrire tous les $(G,V)$ de type PL,
               sous l'hypothèse que  $V$ est un $G$-module semi-simple, ce qui entraîne que $G$ est un groupe réductif. Commençons par le cas où $V$ est irréductible (le cas général sera traité au \S2.6):      
                        
               \medskip
               \n {\bf Théorème  2.1.1} - {\it Soient $V$ et $G$ comme ci-dessus et supposons que $V$ soit un $G$-module
               irréductible. Il n'y a alors que quatre possibilités}:
               
                	(a) \ $G = 1 \ et \ \dim V = 1$.
	
	        (b) \ $G = \GL_V$.
	        
	        (c) {\it  Il existe une forme bilinéaire alternée non dégénérée sur $V$
	        telle que $G$ soit le groupe symplectique correspondant $;$ on a $G \simeq \Sp_N$ avec $N$ pair $> 0$}.

	        (d)  {\it Il existe une forme quadratique non dégénérée de rang pair $N >2$  sur $V$ telle que $G$ soit le groupe spécial  orthogonal correspondant $;$ on a $G \simeq \SO_N$.}                     
              
              \medskip
              
              La démonstration sera donnée au \S2.3 dans le cas non autodual et au \S2.5 dans le cas autodual.

              \medskip
              \n {\it Remarques.}
              
              1) Dans (d), il est nécessaire d'exclure le cas $N=2$ car
              l'action du groupe $\SO_2$ n'est pas irréductible: il y a deux droites stables. De même, on doit exclure les groupes 
              $\SO_N$ avec $N$ impair $> 1$, car leur action n'est pas semi-simple. 
              
              2) Si la caractéristique de $k$ était  $\neq 2$, on devrait modifier (d)
              en supprimant l'hypothèse que  $N$  est pair (mais en gardant celle que $N \neq 2$); on pourrait supprimer (a), car c'est le cas particulier de (d) où $N=1$.

              \medskip
              
              \n 2.2. {\bf Rappels sur les représentations irréductibles des 
              groupes réductifs.}
                          
              \medskip
              
              On se place dans une situation plus générale qu'au \S2.1: on considère un groupe réductif\footnote{Dans toute la suite, on convient que  `` réductif "
              entraîne `` connexe ''.} $G$ sur  un corps algébriquement clos $k$, et une représentation irréductible   $V$  de  $G$ (non nécessairement fidèle). On ne fait pas d'hypothèse sur la caractéristique de $k$.
              
                           Soit  $T$ un tore maximal de $G$, soit $X $ le groupe des caractères de $T$, et  soit $\O$ l'ensemble des poids de $T$ dans $V$. Soit $N$ le normalisateur de $T$ dans $G$ et soit $W = N/T$ le groupe de Weyl; le groupe $W$ agit sur $T$ et sur $X$; on a $W.\O = \O$. 
                           
                           \smallskip
                 
                 Un élément $\alpha
                 $ de $\O$ est dit {\it extrémal} si                  
                  c'est un élément extrémal de l'enveloppe convexe 
                  de $\O$ dans $X_{{\bf R}} = X\otimes_{\Z} \R$, cf. [EVT II, p.10 et p.47]. Comme $\O$ est fini, cela équivaut à:
                  
                   $\alpha$ {\it n'appartient pas à l'enveloppe convexe de}  $\O \sm \{\alpha\}$,
                   
              \n     ou encore:
              
              {\it Il existe un homomorphisme $f: X \to \R$ tel que $f(\alpha) > f(\beta)$ pour tout $\beta \in \O$ distinct de} $\alpha$.

                                \smallskip
                  Notons $\O_e$ l'ensemble des éléments extrémaux de $\O$. C'est un ensemble non vide. De plus:
                  
                  (2.2.1) {\it Le groupe $W$ opère transitivement sur} $\O_e$.
                  
                  (2.2.2) {\it Si $\chi \in  \O_e$, la multiplicité de $\chi$ est égale à} 1
                  (i.e. $\dim V_\chi = 1$).
                  
                  (2.2.3) {\it La représentation $V$ est déterminée à isomorphisme
                  près par} $\O_e$.
                  
                  (2.2.4) {\it Remplacer $V$ par sa duale
                  remplace $\O$ par $-\O$ et $\O_e$ par} $-\O_e$.
                  
                  (2.2.5) {\it Pour toute orbite  $Y$ de $W$ dans  $X$, il existe
                  une représentation irréductible de $G$ dont l'ensemble des
                  poids extrémaux est} $Y$.
                  
                  \smallskip
                  
   \n {\small            [D'après (2.2.3), cette représentation est unique, à isomorphisme près. On obtient ainsi une bijection entre les classes de représentations irréductibles de  $G$  et les orbites de $W$ dans $X$.]  }
                               
                  \medskip
            \n     {\it  Note}. Les énoncés ci-dessus sont traditionnellement énoncés d'une manière  différente, cf. par exemple  [Ch 58, exposé 16], [Ja 03, II.2],  [MT 12, \S15] et
 [St 67, \S\S12 - 14]):  on se ramène au cas où $G$  est semi-simple et l'on choisit une chambre de Weyl  $C$ dans $X_\R.$ Comme toute orbite de $W$ dans $X$ rencontre  $C$  en un point et un seul, on remplace $\O_e$ par son intersection avec $C$, que l'on  appelle {\it le plus grand poids} de $\O$ (ou de la représentation); cela conduit à classer les représentations irréductibles de $G$ par les éléments de $C \cap  X_\R.$  Dans les démonstrations qui suivent, il ne serait pas  commode d'avoir à choisir  $C$.      
          
          \medskip
          
         \n {\it Exemples}. 
         
         (2.2.6) On prend $G = \GL_n$, $T =$ tore diagonal, $X = \Z^n$ et
          $V=$ représentation standard de dimension $n$. L'ensemble $\O$ est alors la base
          canonique $\{e_1,...,e_n\}$ de $X$; son enveloppe convexe est le simplexe de sommets $e_1,...,e_n$; ses points extrémaux sont les $e_i$, de sorte que $\O_e = \O$.
                    
                    \smallskip
          (2.2.7)  On prend $G = \Sp_{2n}$ ou $G = \SO_{2n}$, et $V=$ représentation standard de dimension $2n$. L'ensemble $\O$ est de la forme $\{e_1,...,e_n,-e_1,...,-e_n\}$, où $\{e_1,...,e_n\}$ est une base de $X$; son enveloppe convexe est formé des $\sum x_ie_i$ avec $\sum |x_i| \leqslant 1$; c'est un ``hyperoctaèdre" de dimension $n$ (le polaire d'un $n$-cube). Ici encore, les points extrémaux sont les sommets et l'on a $\O_e = \O$. 
          
          \smallskip
          
 \n {\it Rappels.}
 
          Lorsque la caractéristique de $k$ est 2, le groupe $\SO_{2n}$ est conjugué à un sous-groupe de $\Sp_{2n}$. De façon plus précise, soit $C(x,y)$ une forme alternée non dégénérée sur $V$ fixée par $\Sp_{2n}$, et soit $T$ un tore
          maximal de $\Sp_{2n}$; on vérifie facilement qu'il existe une unique forme quadratique $q$ sur $V$ ayant les deux propriétés suivantes:
          
          \smallskip
          
            (2.2.7.1) La forme bilinéaire associée à $q$ est $C$.
            
            \smallskip
            
            (2.2.7.2) On a $q(v)=0$ pour tout vecteur propre  $v$ de $T$.
            
            \smallskip
            
        \n    Cette forme quadratique est invariante par $T$. Le groupe orthogonal ${\bf O}_{2n}$ qu'elle définit est contenu dans $\Sp_{2n}$ et contient $T$; c'est 
            l'unique groupe orthogonal ayant ces deux propriétés;  son algèbre de Lie est indépendante du choix de $q$:
             c'est l'ensemble des matrices de la forme $z+z^*$, où $z^*$ désigne l'adjoint de  $z$ par rapport à
             la forme alternée  $C$.

                \bigskip

                \n2.3. {\bf Démonstration du théorème $2.1.1:$ premiers cas.}
                
                \medskip
                Nous revenons aux hypothèses et aux notations du \S2.1. En particulier, $k$ est algébriquement clos de caractéristique 2, $G$ est un sous-groupe réductif de $\GL_V$, et $V$ est un $G$-module semi-simple de type PL.
                
                  \medskip

                Soit $T$ un tore maximal de $G$ et soit $\O$ l'ensemble des poids de $T$ dans $V$. Comme la représentation de $T$ dans $V$ est fidèle, $\O$ engendre le groupe $X$ des caractères de $T$. Par hypothèse, $T$ est de type PL; 
                il existe donc une partie libre $\o$ de $\O$ telle que $\O$ soit contenu dans $\{0\} \cup \o \cup -\o$, cf. \S1.1.b); il en résulte que  $\o$  est une $\Z$-base de $X$. Il y a différentes possibilités, que nous allons considérer séparément:
                
                \smallskip
                
                (2.3.1) {\it On a} \ $\o = \varnothing$.  
                
                \n Dans ce cas, $\O = \{0\}$, et 
                $T= 1$, d'où $G=1$; c'est le cas (a) du théorème  2.1.1. 
                
                \smallskip
           
       A partir de maintenant, on suppose  $\o \neq \varnothing$. Considérons d'abord le cas où $\Omega \cap -\omega =\varnothing$. Ce cas se divise en deux:
   
   \smallskip
   
           (2.3.2) {\it Le cas} \  $\O = \{0\} \cup \o$. 
           
           \n L'enveloppe convexe de $\O$ est alors le simplexe dont l'ensemble des sommets est $\O$. Tous les sommets sont extrémaux. Or, $W$ fixe $0$; comme $\o \neq \varnothing$,
           $W$ n'agit pas transitivement sur les points extrémaux. D'après (2.2.1), ce cas est impossible.
           
           \smallskip
           
           (2.3.3) {\it Le cas} \  $\O = \o$. 
           
           Tous les éléments de $\o$ sont extrémaux; leur multiplicité est
           1, d'après (2.2.2). Si l'on pose  $n = |\o|$, on a $\dim V = n = \dim T$,
           d'où le fait que $T$ est un tore maximal de $\GL_V$. De plus, le commutant de $G$ dans $\GL_V$ est réduit aux homothéties. Nous allons
           voir que cela entraîne $G = \GL_V$, autrement dit le cas (b) du théorème 1. Pour cela, remarquons que, d'après (2.2.1), $W$ permute transitivement les éléments $\{e_1,...,e_n\}$ de $\o$,
           lesquels forment une base de $X$. On peut donc identifier $W$ à un sous-groupe transitif du groupe symétrique ${ S}_n$. De plus, $W$ est engendré par des réflexions; or les seuls éléments de ${ S}_n$ qui soient des réflexions sont les transpositions, et le seul sous-groupe transitif de ${ S}_n$
           qui soit engendré par des transpositions est ${ S}_n$, cf. e.g. [Hu 67, p.171]. On a donc $W = { S}_n$. Cela entraîne que le groupe dérivé 
           $G'$ de $G$ est de type $A_{n-1}$, donc de dimension $n^2-1$.
           Comme $G$ contient le groupe des homothéties, on a
            
             \smallskip

     \qquad      $\dim G = 1 + \dim G' =  n^2= \dim \GL_{n}$, 
     
             \smallskip

           \n d'où  $G = \GL_n$; c'est le cas (b) du théorème  2.1.1.

                \bigskip

                \n 2.4. {\bf Le cas autodual $: $ structure du groupe de Weyl.}
                
                \smallskip
                
               Nous venons de démontrer le théorème  2.1.1 dans le cas particulier
               où $\O$ ne rencontre pas $-\o$. Supposons maintenant que
               $\O \cap -\o \neq \varnothing$. Dans ce cas, $0$ n'est pas extrémal, et 
               les autres éléments de $\O$ le sont. On a donc

               \smallskip
               
    \qquad           $\O_e = \O$  \ si  $0 \notin \O$    \quad  et \quad      $\O_e = \O \sm \{0\}$   \ si  $0 \in \O$.

                \smallskip
                
                Le groupe  $W$  opère transitivement sur $\O_e$; or $\O_e \cap -\O_e$ est stable par $W$, et non vide; il est donc égal à $\O_e$, ce qui signifie que  $\O$ contient $-\o$. On a, soit  $\O = \{0\} \cup \o \cup -\o$, soit $\O =  \o \cup -\o$, autrement dit il existe une base 
                $ \{e_1,...,e_n\}$ de  $X$  telle que $\O_e = \{e_1,...,e_n,-e_1,...,-e_n\}$. D'après   (2.2.4), le fait que $\O_e = - \O_e$ signifie que la représentation $V$ est {\it autoduale} (isomorphe à sa duale), ce qui n'était pas le cas au \S2.3.
                
                Le groupe $W$ permute les $e_i$, aux signes près. Il est donc contenu dans le groupe noté usuellement $\{±1\}^n.{ S}_n$, autrement dit le groupe de Weyl d'un système de racines de type $\sf C_{\it n}$  : produit semi-direct de  ${ S}_n$ par le groupe $I = \{±1\}^n$.  Les éléments de $I$ sont les $(\lambda_i )_{1 \leqslant i \leqslant n}$ avec  $ \lambda_i = ± 1$ pour tout $i$; on les identifie aux automorphismes de $X$ de la forme
           $ \{e_i \mapsto  \lambda_i e_i\}$.
   Les $(\lambda_i)$ avec $\prod  \lambda_i  = 1$ forment un sous-groupe $J$
   de $I$ d'indice 2.  D'où un sous-groupe $J.{ S}_n$ de $I.{ S}_n$ d'indice 2.
   
   \smallskip
   
   \medskip
                
                \n {\bf Proposition 2.4.1} - {\it Le groupe $W$ est égal, soit à $I.{ S}_n = \{±1\}^n.{ S}_n$,
                soit à $J.{ S}_n$.}
                
                \n [Autrement dit, $W$ est isomorphe, soit au groupe de Weyl d'un système de racines de type $\sf C_{\it n}$, soit à celui d'un système de type $\sf D_{\it n}$.]
           
           \medskip     
                Le cas $n=1$ est clair. On va donc supposer $n > 1$ dans ce qui suit.
               
               \medskip 
                La démonstration repose sur les deux propriétés suivantes de $W$:
     
     \smallskip           
   (2.4.2) {\it Il est engendré par des réflexions} (lorsqu'on le considère comme un groupe d'automorphismes de $ X_\R   \simeq \R^n$). 
   
   \smallskip
   
   (2.4.3) {\it Il opère transitivement sur}  $\O_e = \{e_1,...,e_n,-e_1,...,-e_n\}$.   
 
   \smallskip
   
   [La première propriété est commune à tous les groupes de Weyl.
   La seconde résulte de (2.2.1).]  
   
   \smallskip
   Noter, à propos de (2.4.2), que les réflexions appartenant à $I.{ S}_n$ sont de trois types (correspondant aux racines courtes et aux racines longues
   de  $\sf C_{\it n}):$
   
   \medskip
   
   (2.4.4) \ \  $t_i \  \ (1 \leqslant i \leqslant n)  \ \ \ \ \ : e_i \leftrightarrow -e_i$  \ et \ $e_\l \leftrightarrow e_\l$ \ si \ $
   \l \neq i$.
   
   (2.4.5) \  $s_{ij} \ (1 \leqslant i < j \leqslant n) : e_i \leftrightarrow e_j$ \ \ \ et \ $e_\l \leftrightarrow e_\l$ \ si \ $
   \l \neq i,j$.
   
   (2.4.6) \  $s_{ij}'  \ (1 \leqslant i < j \leqslant n) : e_i \leftrightarrow -e_j$ \  et \ $e_\l \leftrightarrow e_\l$ \ si \ $
   \l \neq i,j$.

\medskip

 L'image de $t_i$ dans ${ S}_n$ est l'identité; celle  de $s_{ij}$ (resp. de
$s_{ij}'$) est  la transposition $(ij)$.

\medskip
\n {\it Démonstration de la proposition 2.4.1}.

    \medskip
                
                \n {\bf Lemme 2.4.7} - (a) {\it La projection $W \to { S}_n$ est surjective.}
                
                (b) {\it Pour tout couple $(i,j)$ avec $1 \leqslant i < j \leqslant n$, $W$ contient $s_{ij}$ ou $s_{ij}'$.}
                
                \smallskip
                
   \n {\it Démonstration.}     La démonstration de (a) est  la même que celle utilisée pour
                 (2.3.3): l'image de $W$
                dans ${ S}_n$ est un groupe transitif qui est engendré par des transpositions,
                donc c'est ${ S}_n$. Pour (b), on remarque qu'il existe au moins un couple $(i,j)$ tels que $W$ contienne  $s_{ij}$ ou $s_{ij}'$: sinon, l'image de $W$ dans ${ S}_n$ serait triviale, ce qui est impossible puisque $n >1$. Par conjugaison, on en déduit que tous les couples $(i,j)$ ont cette propriété: en effet,  
                 ${ S}_n$ opère transitivement sur les couples
                $(i,j)$ avec $i \neq j$.
                 
                 \smallskip
                
                  \n {\bf Lemme 2.4.8} - {\it Il existe des $\lambda_i \in \{±1\}$ tels que
                  $W$ contienne le groupe des permutations de} $\{\lambda_1e_1,...,\lambda_ne_n\}.$
                  
                  \smallskip
                  
                    \n {\it Démonstration. }   Soit $i$ tel que $1 < i \leqslant n$. D'après le lemme 2.4.7 (b), appliqué au couple $(i-1,i)$, il existe dans $W$ une réflexion $\sigma_i$
                    qui permute $e_{i-1}$ et $\mu_ie_i$, avec $\mu_i = ±1$, tout en laissant fixes les  $e_\l$ avec $\l < i-1$ ou $\l > i$. Posons:
                    
                    $\e_1 = e_1$,
                    
                    $\e_2 = \sigma_1(\e_1) = \mu_1 e_2$,
                    
                    $\e_3 = \sigma_2(\e_2) = \mu_2 \sigma_2(e_2) = \mu_1 \mu_2 e_3$,
                  
                    ....
                    
                    La réflexion $\sigma_i$ permute $\e_{i-1}$ et $\e_i$, et fixe
                    les autres $\e_l$. Or les transpositions $(12), (23),(34), ...$ engendrent le groupe symétrique ${ S}_n$. D'où le lemme, en prenant :
                    
                    $\lambda_1 = 1, \  \lambda_2 = \mu_1,  \  \lambda_3 = \mu_1\mu_2,$ \ \ ...
               
               \medskip
               
               \n {\it Fin de la démonstration de la proposition 2.4.1.}  D'après le lemme 2.4.8, on peut supposer que les $e_i$ ont été choisis tels que leur groupe de 
               permutations soit contenu dans $W$. On a donc ${ S}_n \subset W$.
               Cette inclusion est stricte, à cause de (2.4.3). On en déduit, d'après (2.4.2), que $W$ contient au moins une réflexion n'appartenant pas à ${ S}_n$, i.e. du type $t_i$ ou du type
               $s_{ij}'$, cf. (2.4.4) et (2.4.6). 
               
               \smallskip
               Supposons d'abord que $W$ contienne l'une des $t_i$. Vu la transitivité de ${ S}_n$ opérant sur $[1,n]$, $W$ contient tous les $t_i$. Or ceux-ci engendrent le groupe $I = \{±1\}^n$. On a donc $W= I.{ S}_n$.
               
               \smallskip
               
               Supposons maintenant que $W$ contienne l'une des  $s_{ij}'$.
   Vu la transitivité de ${ S}_n$ sur les couples $(i,j)$ avec $i\neq j$, le groupe
              $W$ contient tous les $s_{ij}'$, et il contient aussi les                 
         produits $s_{ij}s_{ij}'=t_it_j$. Or les $t_it_j$ engendrent le sous-groupe $J$ de  $I$. Donc $W$ contient $J.{ S}_n$,
         et, comme ce groupe est d'indice 2 dans $I.{ S}_n$, on a, soit $W =
         J.{ S}_n$, soit $W = I.{ S}_n$.

                \bigskip

                \n 2.5. {\bf  Le cas autodual $: $ fin de la démonstration du théorème  2.1.1.}
                
                \smallskip
                
                On conserve les notations et hypothèses du \S2.4. 
                
                    \smallskip
                
             Puisque la représentation $V$ est autoduale, le groupe $G$ ne contient pas les homothéties. C'est donc un groupe semi-simple.
             
                 \smallskip
        Considérons d'abord le cas $n=1$.   Le groupe $G$ est de rang 1, et la représentation irréductible $V$ a pour poids, soit $\{e_1,-e_1\}$, soit  $\{e_1,0,-e_1\}$; de plus  $e_1$ est une base du groupe $X$ des caractères
        de $T$.  Le premier cas conduit à $G = \SL_2 = \Sp_2$, la représentation $V$ de $G$ étant la représentation naturelle de dimension 2.
        Le second cas est impossible, car toute représentation irréductible 
 non triviale de $\SL_2$ en caractéristique 2 est un produit tensoriel de transformées de Frobenius de la représentation standard, et  0
        n'en est pas un poids. 
        
            \smallskip
        
        A  partir de maintenant, nous supposons  $n > 1$. D'après la proposition 2.4.1, il y a deux possibilités:
        
        \smallskip
                      
             (2.5.1) {\it On a} $W = I.{ S}_n$; {\it le système de racines de $G$ est de type $\sf B_{\it n}$ ou} $\sf C_{\it n}$.
             
             (2.5.2) {\it On a} $W = J.{ S}_n$; {\it le système de racines de $G$ est de type} $\sf D_{\it n}$.
             
             [Pour $n=2$, le type $\sf D_{2}$ doit être interprété comme $\sf A_{1} \times \sf A_{1}$; cela traduit le fait que $\SO_4$ est isogène à $\SL_2 \times \SL_2$. C'est le seul cas où le système de racines n'est pas irréductible.]
             
             \smallskip
             
             Il reste à préciser la structure du groupe $G$ (pas seulement à isogénie près), et celle de sa représentation $V$. Nous allons pour cela déterminer le système de racines $R$ de $G$. Soit $\alpha$ une racine; 
             c'est un poids de $T$ dans $\Lie(G) \subset \End(V)$, et c'est donc la différence de deux éléments de $\O$; on en déduit que  $\alpha$ est de l'un des types suivants: 
             
             \smallskip            
             
             (2.5.3)  $\alpha = \pm e_i \pm e_j$      avec $i\neq j$;
             
              (2.5.4) $ \alpha = \pm 2 e_i$;
              
              (2.5.5)        $ \alpha = \pm e_i$.
              
              \smallskip
              
           \n   [Le troisième cas ne peut se présenter que si $0$ est 
              un poids de $V$.]
              
          \medskip
          
             \n {\bf Lemme 2.5.6} - {\it Supposons $n \geqslant 3$. Il y a au plus trois
             possibilités pour le système de racines} $R$:
             
             (i) {\it c'est l'ensemble de tous les éléments de $X$ de type} (2.5.3).
             
             (ii) {\it c'est l'ensemble de tous les éléments de $X$  de type $(2.5.3)$, et de tous ceux de type} (2.5.4).
             
             (iii) {\it c'est l'ensemble de tous les éléments de $X$ de type $(2.5.3)$, et de tous ceux de type} (2.5.5).
             
             \smallskip
             
             \n {\it Démonstration.} La proposition 2.4.1 montre que $W$ opère transitivement sur chacun des trois types. De plus, les types (2.5.4) et (2.5.5) ne peuvent pas coexister, car aucune racine n'est le double d'une autre. Le type
             (2.5.3) doit être présent, car sinon le système de racines serait 
             de type $\sf A_{1} \times   ... \times \sf A_{1}$, ce qui n'est pas le cas,
             on l'a vu [c'est ici que l'hypothèse $n \geqslant 3$ est utilisée: lorsque $n=2$, le type $\sf D_{2}$ est isomorphe à $\sf A_{1} \times \sf A_{1}$].  D'où le lemme.
             
             \smallskip

  \n {\bf Lemme 2.5.7} - {\it Conservons les hypothèses et notations du lemme $2.5.6$.
  Alors~$:$
  
  Dans le cas} (i), {\it on a $G \simeq \Sp_{2n}$ et $V$ est la représentation
  naturelle de $G$ de dimension $2n$.
  
  Dans le cas} (ii), {\it on a $G \simeq \SO_{2n}$ et $V$ est la représentation
  naturelle de $G$ de dimension $2n$.
  
  Le cas }(iii) {\it est impossible.}
  
  \smallskip

 \n {\it Démonstration.} Dans le cas (i), le lemme 2.5.6 montre que le système de racines $R$ est celui du type $\sf C_{\it n}$, les poids étant ceux de la représentation naturelle. D'où le résultat. Même argument pour le cas (ii),
 avec $\sf C_{\it n}$ remplacé par $\sf D_{\it n}$. Quant au cas (iii), il donne un système
 de racines de type $\sf B_{\it n}$, avec pour poids ceux de la représentation
 naturelle de dimension $2n+1$; en caractéristique 2, cette représentation
 n'est pas irréductible: la représentation irréductible correspondant aux
 poids extrémaux $\pm e_i$ est le quotient de la précédente par l'unique sous-espace stable de dimension 1; l'élément $0$ de $X$ n'en est pas un poids, et
 par conséquent les  $\pm e_i$ ne peuvent pas en être des poids [d'ailleurs
 l'image de $\SO_{2n+1}$ dans cette représentation n'est pas de type $\sf B_{\it n}$: c'est le groupe $\Sp_{2n}$, qui est de type $\sf C_{\it n}$].
 
 \smallskip
 
  Le lemme 2.5.7 entraîne le théorème  2.1.1, lorsque  $n\geqslant 3$. Le cas
  $n=2$ se traite de façon analogue. La seule différence est que, dans l'énoncé du lemme 2.5.6, il faut ajouter les deux cas suivants: 
   
  (iv) {\it le système $R$ est formé des éléments de $X$ de la forme} $\pm 2e_i$;
  
  (v) {\it le système $R$ est formé des éléments de $X$ de la forme} $\pm e_i$.
  
 Les deux cas conduisent à un groupe $G$ de type $\sf A_{1} \times \sf A_{1}$; le premier donne pour $V$ le produit tensoriel des deux représentations naturelles de dimension 2; cela donne le groupe $\SO_4$. Le second cas
 ne serait possible que si  $0$ était un poids de $V$, ce qui ne peut pas se produire. On a donc bien obtenu tous les cas énumérés dans le théorème 2.1.1.
 
 \medskip
 
 \n {\it Remarque.} Une curieuse conséquence du théorème 2.1.1 est que, 
 si $V$ est irréductible non triviale, alors $0$ {\it n'est pas un poids de} $T$. 
 Ou, de façon équivalente:
 
(2.5.8)    {\it Si $V$ est une représentation semi-simple de type PL d'un groupe linéaire connexe lisse $G$, tout élément de $V$ qui  est fixé par un tore
   maximal de $G$ est fixé par $G$.}      

  C'est là une propriété spéciale à la caractéristique 2. Il serait intéressant d'en avoir une démonstration {\it a priori}.

                \bigskip

                \n 2.6. {\bf Le cas semi-simple.}
                                
                \smallskip
                
             Passons  maintenant au cas où $V$ est une représentation semi-simple, mais non nécessairement irréductible, du groupe $G$. Comme 
             au \S2.1, nous supposons que la représentation considérée
             est {\it  fidèle}, i.e. qu'elle donne un plongement de $G$ dans $\GL_V$; d'après [SGA  3, ${\rm VI_B}$, cor.1.4.2], cela revient à
             dire
             que son noyau, au sens de la théorie des schémas, est trivial.
             
             \smallskip
             
                          De façon plus précise, considérons le cas où la représentation $V$ se décompose en:
                                     
                \smallskip

             (2.6.1) $ V = \oplus_{i\in I} V_i$,
                          
                \smallskip
                             
         \n    où chaque $V_i$ est irréductible, et:
                       
                \smallskip
                
             (2.6.2) {\it Si $i \neq j, \ V_i$ n'est isomorphe, ni à $V_j$, ni à sa duale.}
                                 
             %\medskip
                
          %          \n {\small [L'hypothèse (2.6.2) n'est là que pour éliminer
              %   des cas triviaux: si un facteur irréductible d'une représentation semi-simple $V$ a une multiplicité $> 1$ (ou est isomorphe au dual d'un autre facteur), on peut l'éliminer sans changer le groupe $G$.] }
                
                \medskip

                Notons $G_i$ l'image de $G$ dans $\GL_{V_i}$. Le groupe $G$ est un sous-groupe du produit direct  $\prod_{i \in I} G_i$.
                                            
                \smallskip
         
                         \n {\bf Théorème 2.6.3 - } {\it Si la représentation $V$ est de type PL, on a $G = \prod_{i \in I}  G_i;$ en particulier, $G$ est isomorphe à un 
                         produit de groupes de type $\GL, \Sp$ et $\SO$.}
                                     
\smallskip

\smallskip

\n {\it Démonstration}. Soit  $T$ un tore maximal de $G$, et soit  $T_i$
l'image de $T$ dans  $G_i$. On a $T \subset \prod_{i \in I}  T_i$.
                
\medskip

\n {\bf Lemme 2.6.4} - {\it On a} $T = \prod_{i \in I}  T_i$.

         \smallskip
         
         \n {\it Démonstration du lemme 2.6.4}. Soit $X_i$ le groupe des caractères
         de $T_i$, et soit $X$ celui de $T$. La projection $T \to T_i$ donne une
         injection $X_i \to X$, ce qui nous permet d'identifier les $X_i$ à des sous-groupes de  $X$; le fait que $T \to \prod  T_i$ soit injectif montre que $X$
         est engendré par les $X_i$. Le lemme 2.6.4 revient à dire que $\oplus \ X_i \to X$ est injectif, i.e. que $X$ est somme directe des $X_i$; c'est ce que nous allons démontrer.

         Puisque les $V_i$ sont de type PL, le théorème 2.1.1 montre qu'il existe
         une base $\o_i$ de $X_i$ qui est formée de poids de $V_i$, donc de
         poids de $V$. On a :
         
         \smallskip
         
            (2.6.5) {\it Si $i\neq j$, on a \  $\o_i \cap \o_j = \varnothing$ \ et} \ $\o_i \cap -\o_j = \varnothing$. 
                
                  \smallskip
En effet, supposons que $\o_i$  et $ \o_j$ (resp. $-\o_j$) ait un élément commun. Comme cet élément est un élément extrémal de $V_i$ et de $V_j$
 (resp. de la duale $V'_j$ de $V_j$), cela entraîne, d'après (2.2.3) et (2.2.4)
 que $V_i$ est isomorphe à $V_j$ (resp. à sa duale), ce qui contredit l'hypothèse (2.6.2).
 
 \smallskip
   En particulier, les $\o_i$ sont disjoints. Soit $\o = \cup \ \o_i$. D'après
   (2.6.5), on a $\o \cap - \o = \varnothing$. D'après (1.1.3), cela entraîne
   que $\o$ est une partie libre de $X$, d'où le fait que $X$ est somme directe des $X_i$.
   
   \smallskip
   
   \n {\it Fin de la démonstration du théorème 2.6.3.} Le lemme 2.6.4 montre que le rang  (dimension d'un tore maximal) de $G$ est égal à celui de
 $\prod   G_i$. L'égalité $G = \prod G_i$ en résulte, en vertu du résultat suivant:
 
 \smallskip
 
 \n {\bf Proposition 2.6.6} - {\it Soit $H_i$ une famille finie de groupes réductifs, et soit $H$ un sous-groupe réductif de $\prod H_i$ ayant les deux propriétés suivantes}:
   
   (2.6.7) {\it Les projections $H \to H_i$ sont surjectives.}
   
   (2.6.8) {\it Le rang de $H$ est égal à la somme de ceux des $H_i$.
   
  \n On a alors} $H = \prod  H_i$.
   
   \smallskip
   
   \n {\it Démonstration}. Il suffit de traiter le cas où il y a deux groupes $H_1$ et $H_2$. Si l'on pose $N = \Ker (H \to H_1) = H \ \cap \ \{1\} \! \times \! H_2$ , le rang de $N$ est égal à la différence de ceux de $H$ et de $H_1$, autrement dit à celui de $H_2$. Mais $N$ est un sous-groupe
   normal de $H_2$. Le rang  de $H_2/N$ est donc $0$, ce qui n'est possible que si $N= H_2$, auquel cas $H = H_1 \times H_2$.
              
              \smallskip      
              
              \n {\bf Corollaire 2.6.9} - {\it Tout sous-groupe réductif  $V$ de $\prod  H_i$ de rang maximum est produit direct des $V \cap H_i$}.
              
         \n {\small     [Rappelons qu'un sous-groupe d'un groupe réductif $G$ est dit {\it de rang maximum} si son rang est égal à celui de $G$.]}
               
               \smallskip
               
     \n {\it Démonstration}. Soit $V_i$ la projection de $V$ sur le $i$-ième facteur. En appliquant la proposition précédente 
     au produit des $V_i$, on voit que $V = \prod  V_i$  et $V_i = V \cap H_i$.
         
                      \bigskip

 \centerline          {\bf \S3 - Groupes unitaires sur un corps algébriquement clos de caractéristique 2}
                            
              \medskip
              
            \n  A partir de maintenant, le corps de base $k$ est supposé de caractéristique 2.
              
              \medskip
              
        \n 3.1. {\bf Algèbres à involution et groupes unitaires - préliminaires.}

     \smallskip

     Soit $A$ une {\it algèbre à involution} sur  $k.$ Rappelons que cela signifie que $A$ est une $k$-algèbre associative, à élément unité, munie d'une application $k$-linéaire $a \mapsto a^*$ telle que $a^{**}=a$  et
          $(ab)^* = b^*a^*$ quels que soient $a,b \in A$. On suppose en outre
que $\dim_k A < \infty$. 

\smallskip
  Nous associerons à $A$ deux groupes algébriques:
  
  \smallskip
  
  (3.1.1) Son {\it groupe multiplicatif} $\GL_{1,A}$. Il représente le foncteur
  
  \smallskip
  
    $k' \mapsto (A\otimes_k k')^\times$, où  $k'$ parcourt la catégorie des $k$-algèbres commutatives;
    
    \smallskip
    
    \n en particulier le groupe de ses $k$-points est $A^\times$. 
    
    C'est un sous-groupe algébrique du groupe
    $\GL_V$, où  $V$ est le $k$-espace vectoriel sous-jacent à  $A$. C'est un groupe connexe et 
    lisse, qui est réductif si  $A$  est une algèbre absolument semi-simple; il est muni d'une anti-involution $g \mapsto g^*$, déduite de celle de $A$.
    
     \smallskip
     
    (3.1.2) Son {\it groupe unitaire schématique} $U_A^{\sch}$; c'est le sous-groupe du précédent formé des points  $u$  tels que $uu^*=1$. Ce n'est pas en général un groupe lisse: lorsque  $A = k$, on a  $\GL_{1,A} = \G_m$ et $U_A^{\sch}= \m$.
    
    \smallskip
    
    Soit $U_A$ le schéma réduit associé à $U_A^{\sch}$. Lorsque c'est un sous-schéma en groupes de $U_A^{\sch}$, nous l'appellerons 
    le {\it groupe unitaire} de $A$; cela se produit lorsque $k$ est parfait (cf. [SGA 3, exposé VI$_{\rm A}$, \S0.2]), ou, plus généralement, lorsque  $A$  est déduit, par extension des scalaires, d'une algèbre à involution sur un corps parfait (exemple: $A = k[G]$,
    où $G$ est un groupe fini, cf. \S5).
    
    \smallskip
\n 3.1.3. {\it Algèbres de Lie}.     L'algèbre de Lie  de $U_A^{\sch}$ est l'espace vectoriel $H_A$ des
     {\it éléments hermitiens} de $A$, autrement dit des éléments $a\in A$ tels que $a^*=a$, le crochet étant défini par la formule $[a,a'] = aa' + a'a$: cela
     résulte de la définition de  $U_A^{\sch}$, appliquée à l'algèbre $k'= k[x]/(x^2)$ des nombres duaux.   L'algèbre de Lie de $U_A$ est une sous-algèbre de $H_A$. Nous en donnerons quelques propriétés plus loin (cf. proposition 3.3.1),  mais je n'en connais pas de description générale, mis à part le cas, dû à Merkurjev, où $A$ est l'algèbre d'un groupe fini; dans ce cas, on verra au \S5.1 que
     $\Lie(U_A)$ est un hyperplan de $\Lie(U_A^{\sch})$.
     
     \medskip
                    \medskip                 
        \n 3.2. {\bf Sous-algèbres étales et tores maximaux.}
        
        \smallskip
        
        Jusqu'à la fin du \S3, on suppose  que $k$ est {\it algébriquement clos} (et, bien sûr,
        de caractéristique 2).
        \smallskip
        
        Soit $A$ une $k$-algèbre à involution de dimension finie, et soit $C$ une $k$-sous-algèbre commutative
        de $A$ satisfaisant aux deux conditions suivantes:
         
          \smallskip
         (3.2.1) {\it Elle est stable par l'involution de} $A$.
         
          \smallskip
         (3.2.2) {\it C'est une algèbre étale}, autrement dit (puisque $k$ est algébriquement clos), c'est une algèbre diagonalisable, i.e. un produit de copies de $k$, cf. [A V, \S6, n$^{\rm o}$ 3].
          \smallskip
         
         D'après (3.2.2), on peut écrire $C$ sous la forme $C = k^I$, où $I$
         est un ensemble fini que l'on peut interpréter, soit comme le spectre
         de $C$, soit comme $\Hom_{\rm alg}(C,k)$, soit comme l'ensemble des idempotents indécomposables de $C$. L'involution de $C$ opère sur $I$;
         on peut donc décomposer $I$ en deux parties:
         
         \smallskip
          $I_0 = $  éléments de $I$ fixés par $i \mapsto i^*$;
          
            \smallskip

          $I_1 = I \sm I_0 = $  éléments $i$ de $I$ tels que $i^* \neq i$.
          
            \smallskip

          Le groupe unitaire $U_C$ de $C$ est formé des familles $(u_i)_{i\in I}$ telles que  $u_i= 1$ si $i \in I_0$, et $u_i u_{i^*} = 1$ si $i\in I_1$. Si l'on décompose $I_1$ sous la forme $I_1 = J \cup J^*$, avec $J \cap J^* = \varnothing$, on voit que les valeurs de $u_i$ pour $i\in J$ déterminent toutes les autres, et peuvent être choisies arbitrairement. En particulier,
          $U_C$ {\it est un sous-tore de $U_A$}, isomorphe à  $(\G_m)^J$. De façon plus précise, pour tout $i\in I$, l'application
          $e_i : u   \mapsto u_i$ est un caractère de $U_C$; on a  $e_i + e_{i^*} = 0$ \ et $e_i = 0 $ si $i \in I_0$; les $(e_i)_{ i\in J}$ forment une base $\o_J$ du groupe des caractères $X(U_C)$ de $U_C$.
          
          \smallskip

          Soit maintenant $V$ un $C$-module de dimension finie
          sur  $k$. Le groupe  $U_C$ opère sur $V$.

 \smallskip

 \n {\bf Lemme 3.2.3 } - {\it L'action du tore $U_C$ sur $V$ est de type PL,
 au sens du \S1.2.}
 
 \smallskip

 \n {\it Démonstration.} La décomposition de $C$ en produit donne une décomposition correspondante de $V$ en produit de $V_i$. Les poids
 de $U_C$ dans $V$ sont les $e_i$, avec multiplicité $\dim V_i$; ils sont contenus dans l'ensemble   $     \{0\} \cup \o_J \cup -\o_J$, où $\o_J$ est la base de   $X(U_C)$ définie ci-dessus. La condition (1.1.2) est donc satisfaite.
 
 \medskip
 
 \n{\bf Lemme 3.2.4} - {\it Si $T$ est un tore maximal de $U_A$, il existe une sous-algèbre $C$ de $A$, satisfaisant aux conditions $(3.2.1)$ et $(3.2.2)$, telle
 que $T = U_C$.}
 
 \smallskip
                
                 \n {\it Démonstration.} Soit $C$ la sous-algèbre de $A$ engendrée par le groupe $T(k)$ des $k$-points de $T$. Comme $T$ est un tore, cette algèbre est étale: c'est là une propriété générale des sous-tores de $\GL_n$.
                 Comme $T$ est stable par l'involution de $A$, il en est de 
                 même de $C$. Les conditions (3.2.1) et (3.2.2) sont donc satisfaites, et il est clair que $T$ est contenu dans $U_C$; comme $T$ est un tore maximal, on a donc $T = U_C$.

                  \smallskip
           Les lemmes 3.2.3 et 3.2.4 entraînent:
           
           \smallskip
           
           \n {\bf Proposition 3.2.5} - {\it Soit $V$ un $A$-module à gauche de dimension finie
           sur $k$. Si $T$ est un tore maximal de  $U_A$, l'action de   $T$
           sur $V$ est de type PL.} 
                                
                                \n [Précisons que l'on fait opérer $A^\times$ et $U_A$ sur $V$ par multiplication à gauche.]
                                
                                \smallskip
                                
                             Noter que cela s'applique en particulier au cas où $V= A$, avec sa structure naturelle de $A$-module à gauche: on a $U_A \subset \GL_{1,A} \subset \GL_V$, et l'on voit que {\it les tores maximaux de $U_A$ sont de type} PL, comme annoncé au \S1.3.4.

                             \smallskip 
                             
                             \n {\it Remarque} 3.2.6. Ce dernier résultat peut aussi se déduire de la proposition 1.4.9, appliquée au groupe $U_A^{\sch}$. En effet, ce groupe est le fixateur d'une famille de tenseurs quadratiques de $V=A$. De façon plus précise, soit $(e_i)$ une base de $V$ et soit $(\l_j)$ une base du dual  $V'$  de  $V$; notons $\theta_i \in \End(V)$  la multiplication à droite par $e_i$ et notons $\theta_j \in \otimes^2 V'$ la forme bilinéaire $(a,b) \mapsto \l_j(a^*b)$ . On vérifie facilement que :
                             
                             le fixateur des $\theta_i$ est $\GL_{1,A}$;
                             
                             le fixateur des $\theta_i$ et des $\theta_j$ est $U_A^{\sch}$.
                                                                   
                \medskip
                
                \n3.3. {\bf La composante neutre du groupe unitaire.}
                
                \smallskip
                
                Soit $U_A^0$ la composante neutre du groupe $U_A$. C'est la structure de ce groupe qui va nous intéresser. On a tout d'abord:
                
                \smallskip
                
                \n {\bf Proposition 3.3.1} - (a) {\it Si $u \in U_A$ est tel que $u^2 = 1$, on a \ $u \in U_A^0$ et} $1+u \in \Lie(U_A).$

\smallskip
                      
                      (b) {\it Si $x \in A$ commute à $x^*$\!, on a \ $x+x^* \in \Lie(U_A);$ si de plus $x$ est inversible, on a $x^{-1}x^* \in U_A^0$}.
                                  
                  \smallskip
                  
      \n           [Dans cet énoncé, ainsi que dans les suivants, on identifie les groupes lisses $U_A$
      et $U^0_A$ à l'ensemble de leurs $k$-points, de sorte que ``$u \in U_A$'' signifie $u \in U_A(k)$, autrement dit ``$u \in A$ et $uu^*=1$''.]
      
      \smallskip
                  
                  \n {\it Démonstration de} (a). Ecrivons $u$ sous la forme
                  $u=1+ \e$; comme  $u^2=1$ , on a $\e^2 = 0$.  Puisque $u$ est unitaire, on a $u^* = u^{-1} = u$,
                  d'où $\e^* = \e$. Si $t$ est un élément de $k$,
                  posons  $u_t = 1 + t\e$; on a $u_tu_t^* = u_t^2 = 1$.                   L'application $f: t \mapsto u_t$ se prolonge en un homomorphisme du groupe additif $\G_a$ dans le groupe $U_A^{\rm sch}$; comme $\G_a$ est connexe et lisse, l'image de $f$ est contenue dans $U_A^0$; on a donc $u_t \in U_A^0$ pour tout $t$, d'où $u = u_1 \in U_A^0$; comme la dérivée de $f$ en $0$  est $\e$, on a $\e  \in \Lie(U_A).$

                   \smallskip
                  
                  \n {\it Démonstration de} (b). Soit $C$ la sous-algèbre de $A$ engendrée par $x$ et $x^*$. C'est une algèbre commutative, qui est stable par l'involution de $A$. Pour tout $y \in C^\times$, l'élément $ y^{-1}y^*$ est unitaire. L'application $y \mapsto y^{-1}y^*$ définit un homomorphisme $\varphi_C : \GL_{1,C} \to U_A$ dont l'application tangente en l'élément neutre est $y \mapsto y + y^*$. Cela montre que $\Lie(U_A)$ contient tous les
                  $y+y^*$, et en particulier contient $x+x^*.$ De plus $\GL_{1,C} $
                  est connexe, puisque c'est un ouvert dense d'un espace affine; l'image de $\varphi_C$ est donc contenue dans $U_A^0$. 
                  
                  \smallskip
                  
                  \n {\bf Corollaire 3.3.2} - {\it Pour tout $u \in U_A$, on a $u^2 \in U_A^0$ \ et  \ $u+u^{-1} \in  \Lie(U_A)$}.                 
                  
                   \smallskip
                   
                   On applique (b) à $x = u^{-1}$, d'où $x^*=u$, ce qui donne $x^{-1}x^* = u^2$.
                  
 \smallskip
                  
                  \n {\bf Corollaire 3.3.3} - {\it Le groupe $U_A/U_A^0$ est un $2$-groupe
                  abélien élémentaire. } 
                  
                  \smallskip
         En effet, c'est un groupe fini dont tous les éléments sont de carré 1 d'après le corollaire précédent; on sait qu'un tel groupe est abélien (utiliser l'identité $xyx^{-1}y^{-1} = x^2(x^{-1}y)^2(y^{-1})^2 $.)     
                  
                  \medskip
                  
        \n {\it La condition d'engendrement.}
                                   
                  \smallskip 
                  
                   Soit $A'$ la sous-algèbre de $A$ engendrée par les $k$-points de $U_A^0$.
                   Cette algèbre est stable par l'involution de $A$, et l'on a $U_{A'}^0 = U_A^0$; autrement dit, on ne change pas le groupe $U_A^0$
                   lorsqu'on remplace $A$ par $A'$. Nous pourrons donc par la suite nous borner aux algèbres $A$ satisfaisant à la condition:
                   
                   \smallskip
                   
                   (3.3.4) {\it On a $A=A'$, autrement dit $A$ est engendrée comme $k$-algèbre par les $k$-points du groupe} $U_A^0$.                 
                  
                  \medskip

 \n3.4. {\bf Exemple $:$ le cas où $A$ est une algèbre semi-simple.}
 
  \smallskip
    Ce cas ne présente pas de difficulté: on décompose $A$ en produit d'algèbres de matrices; les différents facteurs sont, soit stables par l'involution, soit permutés deux-à-deux; on est ainsi ramené à déterminer les
  involutions sur une algèbre de matrices $\M_n(k), \ n \geqslant 1$. Lorsque $n=1$  on a $\M_1(k) = k$, avec l'involution triviale. Pour $n >1$,
le fait que la caractéristique soit 2 entraîne que  l'involution est associée à une forme bilinéaire  $B$ qui est symétrique et non dégénérée. Si $B$ est alternée, $n$ est pair, on a $U_A^0 \simeq \Sp_n$, et la condition (3.3.4) est satisfaite. Si $B$ n'est pas alternée, l'ensemble des $x \in k^n$ tel que $B(x,x)=0$ est un hyperplan qui est stable par $U_A$;  la condition (3.3.4) n'est donc pas satisfaite. D'où:

  \smallskip
  \n {\bf Proposition 3.4.1} - {\it Si l'algèbre à involution $A$ est semi-simple, et satisfait à la condition d'engendrement $(3.3.4)$, elle est produit direct d'algèbres à involution $A_i$ de l'un des trois types suivants} :
  
   \smallskip
  (3.4.2) {\it Le corps $k$. On a alors \ $U_{A_i}^{\sch} =   \boldsymbol \mu_2$ \ et} \ $U_{A_i} = 1$. 

 \smallskip
(3.4.3) {\it Le produit $\M_{n_i} \times \M_{n_i}^{\rm opp}$ d'une algèbre de matrices par l'algèbre opposée, l'involution étant $(a,b) \mapsto (b,a)$.
On a alors} \ $U_{A_i}^{\sch} = U_{A_i} = \GL_{n_i}$.

 \smallskip
(3.4.4) {\it Une algèbre de matrices $\M_{2n_i}$ munie d'une involution symplectique. On a alors} \ $U_{A_i}^{\sch} = U_{A_i} \simeq \Sp_{2n_i}.$
  
  \smallskip
  
  Noter que cela entraîne que $U_A$ est un groupe réductif; en particulier,
  c'est un groupe connexe.

                  \medskip
                  
                  \n 3.5. {\bf Structure de certains sous-groupes de $U_A$, pour $A$ semi-simple.}
                  
                  \smallskip
                  
                    On conserve les notations du \S \ précédent; en particulier, $A$ est semi-simple, et on la décompose en produit d'algèbres à involution indécomposables:
                    
                    \smallskip
                    
                      $ A = \prod A_i$,
                      
                      \smallskip
                      
                      \n où les $A_i$ sont de l'un des trois types décrits dans la proposition 3.4.1. On a $U_A = \prod U_{A_i}$.
                                           
                      \smallskip
                  
                  \n {\bf Théorème 3.5.1} - {\it  Soit $H$ un sous-groupe algébrique lisse et connexe de $U_A$ satisfaisant à la condition } :
                  
                  \smallskip
                  (3.5.2) {\it Les éléments de $H(k)$ engendrent l'algèbre $A$.}
                  
                  \smallskip
                  
             \n     {\it Les propriétés suivantes sont alors équivalentes}:
                  
                  (3.5.3) {\it L'action de $H$ sur $A$ par multiplication à gauche est de type} PL.
                  
                  (3.5.4) {\it Le groupe $H$ est réductif de rang égal à celui de $U_A$}.
                  
                  (3.5.5) {\it Le groupe $H$ est produit direct de ses projections $H_i$ sur les $U_{A_i}$, et l'on a $:$
                  
                   $H_i = U_{A_i}$ dans les cas (3.4.2) et (3.4.3), 
                   
                     $H_i = U_{A_i}$ ou \ $H_i \simeq \SO_{2n_i}$ dans le cas (3.4.4) où \  $U_{A_i} \ \simeq \Sp_{2n_i}$.}
                  
         \smallskip
         
         \n {\it Démonstration.} S'il y a des facteurs de $A$ de type (3.4.2), on peut les supprimer sans changer ni $U_A$ ni $H$. Supposons donc qu'il n'y ait aucun tel facteur.
           
           On a évidemment (3.5.5) $ \Rightarrow$ (3.5.4).          D'autre part, si (3.5.4) est satisfaite,
         les tores maximaux de $H$ sont des tores maximaux de $U_A$, et l'on a vu que ceux-ci sont de type PL, cf. proposition 3.2.5; cela montre que (3.5.4) entraîne (3.5.3). Reste à prouver que (3.5.3) implique
         (3.5.5). Supposons donc que (3.5.3) soit satisfaite, et, 
pour chaque indice  $i$, soit $V_i$ l'espace vectoriel défini de la manière suivante:

          Si $A_i$ est de type (3.4.3), on prend $V_i = k^{n_i}.$
          
          Si $A_i$ est de type (3.4.4), on prend $V_i = k^{2n_i}.$
          
       \n    Dans chaque cas, il y a une action naturelle de $U_{A_i}$ sur $V_i$. On en déduit une action de $U_A^0$ sur $V = \oplus_{i \in I} V_i$. De plus, les $V_i$ sont des $H$-modules
          irréductibles d'après (3.5.2). On peut donc appliquer au 
          $H_i$-module $V_i$ le théorème  2.1.1 du \S2.1; cela donne les assertions sur la structure de $H_i$.
       D'autre part, un argument analogue à celui fait ci-dessus montre que, si $i \neq j$,   la représentation $V_i$ de $G$ n'est isomorphe, ni à $V_j$, ni à sa duale. D'après le théorème 2.6.3, cela entraîne que   $H$ est le produit des $H_i$, ce qui achève la démonstration.

                  \medskip

                     \n 3.6. {\bf Structure de $U^0_A$ dans le cas général. }
                  
                  \smallskip
                  
                  Soit $A$ une $k$-algèbre à involution de dimension finie.  
                  
                  \smallskip
                                   
                  \n {\bf Théorème 3.6.1} - {\it Il existe une suite exacte}
                  
                  \smallskip
                  
                  (3.6.2) \ \ $ 1 \ \to N \ \to \  U^0_A \ \to \ H \ \to 1$,

                  \smallskip
                  
        \n          {\it où  $N$  est unipotent connexe, et où  $H$  est un produit de groupes $H_i$
        isomorphes, soit à $\GL_{n_i}$, soit à   $\Sp_{2n_i}$,  soit à} $\SO_{2n_i}$.

        [Rappelons que  $k$  est supposé algébriquement clos ; le cas plus général où  $k$  est parfait sera examiné au \S4.]  
                  \smallskip

        \n {\it Démonstration.}  On peut supposer que $A$ satisfait à la condition d'engendrement (3.3.4). Faisons cette hypothèse, et soit $\r$ le radical de $A$. L'algèbre $A/\r$ est une algèbre à involution semi-simple; son groupe unitaire $U_{A/\r}$ est un groupe réductif, cf. \S3.4. Soit $\pi : U^0_A \to     U_{A/\r}$ l'homomorphisme défini par $A \to A/\r$,
        soit  $N$  le noyau de  $\pi$ (au sens schématique) et soit  $H$ son image (``image schématique'', i.e. plus petit sous-groupe de   $U_{A/\r}$ contenant l'image par $\pi$ des $k$-points de $U^0_A$). On a alors la suite exacte (3.6.2); en effet, la proposition 2.5.2 de [SGA 3, exposé VI A, édition révisée, p.319] montre que l'homomorphisme $U^0_A/N \to U_{A/\r}$ est une immersion fermée, et son image est égale à $H$ car $U^0_A$ est réduit (parce que quotient d'un schéma réduit), et donc lisse (car c'est un schéma en groupes sur un corps parfait) \footnote{Je dois ces explications à Michel Raynaud. Par la suite, on se servira seulement de l'exactitude de la suite $1 \to N(k) \to  U^0_A(k) \to H(k) \to 1$, qui est évidente.}.  Il reste à prouver:    
        \smallskip
        
        \n {\bf Lemme 3.6.3} - {\it Le groupe $N$ est unipotent connexe. Le groupe  $H$  satisfait aux propriétés du théorème 3.5.1 relativement à $A/\r;$ en particulier, c'est  un produit de groupes  $H_i$
        isomorphes, soit à $\GL_{n_i}$, soit à   $\Sp_{2n_i}$,  soit à} $\SO_{2n_i}$. 
        
        \smallskip
        
        \n {\it Démonstration du lemme 3.6.3}. Le groupe  $N$  est un sous-groupe du groupe dont les $k$-points sont de la forme $1 + x$ avec $x \in \r$; or ce groupe est unipotent, comme on le voit par dévissage à l'aide des puissances de $\r$. D'autre part,  l'action de  $U_A^0$ sur $A$ est de type PL, cf. proposition 3.2.5; il en est donc de même de celle de $U_A^0$
        sur $A/\r$, cf. \S1.2; même chose pour celle de $H$ sur $A/\r$. On peut appliquer le théorème 3.5.1 à $H$ et à $A/\r$: en effet, les conditions (3.5.2) et (3.5.3) sont satisfaites; d'après (3.5.5), le groupe $H$ a les propriétés voulues. Reste à montrer que $N$ est connexe. Cela résulte du lemme suivant :        
        \smallskip
        
        \n {\bf Lemme 3.6.4} - {\it Soit $1 \to G_1 \to G_2 \to G_3 \to 1$ une suite exacte de groupes algébriques sur $k$. 
        Supposons que  $G_1$ soit unipotent, que $G_2$ soit lisse et connexe, et que $G_3$ soit réductif. Alors $G_1$
        est connexe.}
        
        \smallskip
        
         \n {\it Démonstration du lemme 3.6.4}.  Après passage au quotient par $G^0_1$ on est ramené au cas où  $G_1$ est fini étale
         et d'ordre une puissance de 2. Comme $G_2$ est connexe, son action sur $G_1$ est triviale, autrement dit
         $G_1$ est contenu dans le centre de $G_2$  D'autre part, il est clair que tout sous-groupe unipotent lisse normal connexe de $G_2$ est trivial; comme $G_2$  est lisse, cela signifie que $G_2$ est réductif. Or le centre d'un groupe réductif est de type multiplicatif; ainsi, $G_1$ est à la fois unipotent et de type multiplicatif; il est donc trivial.        
        
        \smallskip
        
        \n {\it Remarque}. Le groupe $N$ construit dans la démonstration du théorème 3.6.1 n'est pas nécessairement lisse (on en verra un exemple dans 5.5.19). Le groupe réduit correspondant  $N^{\rm red}$ est le {\it radical unipotent $R_u(U^0_A)$ de}  $U^0_A$ . Le groupe $\widetilde {H}  =
        U^0_A/ N^{\rm red}$  est le {\it plus grand quotient réductif de} $U^0_A$: ce groupe est lié à  $H$  par une isogénie    
        $\widetilde{H} \to H$ dont le noyau est infinitésimal et unipotent; en fait, d'après [Va 05],   $\widetilde {H}$ se déduit de $H$ en remplaçant certains des facteurs $\Sp_{2n_i}$ par des facteurs $\SO_{2n_i+1}$. Ce genre de phénomène est spécial à la caractéristique 2.
                                
                    \medskip
                    
                    \n 3.7. {\bf Exemple  de remplacement  de \ $\Sp_{2n} $ par \ $\SO_{2n}$ (et même par ${\bf O}_{2n}$)}    
                                    
                    \smallskip

        Soit $V$ un $k$-espace vectoriel de dimension $2n$, avec $n > 0$,
        muni d'une forme alternée non dégénérée, notée $C(v,v')$. On munit
        l'algèbre $E = \End(V)$ de l'involution correspondante : si  $e \in E$, on a
        \smallskip
        
 (3.7.1) \quad        $C(ev,v')\ = \ C(v,e^*v')$ \ pour tous \   $v,v'\in V$. 
         
           \smallskip
         
         Le groupe unitaire $U_E$ correspondant est le groupe symplectique $\Sp(V)$. On note $H$ l'ensemble des éléments hermitiens de $E$, et l'on note $H_o$ le sous-espace de $H$ formé des éléments de la forme $e+e^*$ avec 
 $e\in E$.
 
 \smallskip
 
 Considérons une algèbre à involution $R$ de radical $\r$ telle que  $\r^2 = 0$ et $R/\r = E$. Il y a une structure naturelle de $E$-module à gauche sur $\r$. Faisons l'hypothèse:
 
 \smallskip
 (3.7.2) {\it Il existe une base $\o$ du $E$-module $\r$, avec $\o$ central dans $R$, et hermitien.}
 
 \smallskip   
 
 \n {\it Exemple} 3.7.3. On prend pour $R$ l'algèbre à involution déduite de $E$ par extension des scalaires à $k[t]/(t^2 )$,
 et l'on prend $\o = 1 \otimes t$. 
 
 \smallskip
 
 Nous allons classer les couples $(R,\o)$:
 
 \smallskip
 \n  {\bf Proposition 3.7.4} - (a) {\it Les classes d'isomorphisme des couples $(R,\o)$ ci-dessus correspondent bijectivement
 } (par une bijection définie plus loin) {\it aux formes quadratiques sur $V$ dont la forme bilinéaire associée est un multiple de $C$. }
 
 (b) {\it Si $(R,\o)$ correspond à la forme quadratique $q$, l'image de $U_R \to \Sp(V)$ est le groupe orthogonal ${\bf O}(q)$.}
 
 \medskip
 
 \n  {\bf Corollaire 3.7.5} - {\it Il y a trois possibilités pour l'image de $U_R \to \Sp(V)$} :
 
  (i) {\it C'est $\Sp(V);$ ce cas ne se produit que pour l'exemple 3.7.3.}
  
  (ii) {\it C'est le sous-groupe de  $\Sp(V)$ fixant une droite de $V$.}
  
  (iii) {\it C'est le groupe orthogonal ${\bf O}(q)$ d'une forme quadratique $q$ dont la forme bilinéaire associée est égale à $C$.} 
  
  [Noter que, dans le cas (iii), l'image de $U_R \to \Sp(V)$ est strictement plus grande que $\SO(q)$.]
  
  \smallskip
  \n  Ces trois cas correspondent aux trois possibilités suivantes pour  $q$:
    
    (i) $q = 0$;
    
    (ii) $q \neq 0$ et la forme bilinéaire associée est  0; il existe alors $v_0 \in V, v_0 \neq 0$, tel que $q(v) = C(v,v_0)^2$ pour tout  $v\in V$.
  
  (iii) la forme bilinéaire associée à $q$ est un multiple non nul de $C$.
  
  \medskip
  
  \n {\it Démonstration de la proposition 3.7.4} (a). Soit $f : E \to R$ un {\it relèvement} de $E$ dans $R$, autrement dit un homomorphisme tel que le composé  $E \to R \to E$ soit l'identité. On sait ([A VIII, \S13, n$^{\rm o}$ 6]) qu'un tel $f$ existe, et que, si $f'$ est un autre relèvement, il existe  $z \in E$ such that $f'(x) = (1+z\o )f(x)(1+z\o )^{-1}$ pour tout $x\in E$, i.e.
        $f'(x) = f(x) + (zx + xz)\o$. Si l'on applique ceci à $f'(x) = f(x^*)^*$, on obtient l'existence de $z \in E$ tel que:
        
        \smallskip
        
 (3.7.6)        $f(x^*)^* = f(x) + (zx + xz)\o$.
 
 \smallskip
 
\n En remplaçant  $x$ par $x^*$ cela donne:
        
        \smallskip
                 (3.7.7) \ $ f(x)^* = f(x^*)+ (zx^*+x^*z) \o$                 \smallskip
                 
            \n     En conjuguant (3.7.7) et en ajoutant (3.7.6) on obtient:
                 
                 \smallskip
                 
                $zx+xz+xz^*+z^*x = 0$, \ i.e. \ $(z+z^*)x = x(z+z^*)$\  pour tout $x \in E$.

           \smallskip
           
      \n     Cela montre que  $z+z^*$ appartient au centre de  $E$, i.e.:
         
         \smallskip
       (3.7.8)  \ $z+z^* = \lambda$ \ avec  $\lambda \in k$.
       
       \smallskip
       
      \n  Soit $q_z$ la forme quadratique $v \mapsto C(zv,v)$. La forme bilinéaire associée est $\lambda C$: cela résulte de
       (3.7.8). De plus, si l'on change les choix de $f$ et de $z$, on constate que $z$ est remplacé par $z + h$ avec $h \in H_o$, et cela ne change pas $q_z$. Ainsi, la forme  $q_z$ est {\it canoniquement associée à}  $(R,\o)$. Cette forme
       {\it détermine $(R,\o)$ et peut être choisie arbitrairement}. En effet, elle détermine le choix de $z$, à l'addition près d'un élément de $H_o$; et, lorsqu'on connaît  $z$, le couple $(R,\o)$ est isomorphe à l'algèbre $E \otimes k[t]/(t^2)$, munie de l'élément $\o = 1 \otimes t$ et de l'involution $x + yt \mapsto x^* + (y^*+zx^* + x^*z)t$.
 
 \smallskip
 
  \medskip
  
  \n {\it Démonstration de la proposition 3.7.4} (b).
  
   Soit $u\in U_E$. Pour que $u$ appartienne à l'image de $U_R \to U_E$, il faut et il suffit qu'il existe $b\in E$ tel que $(f(u)+b\o)(f(u)^*+b^*\o) = 1$. En utilisant (3.7.7), cela s'écrit:
  
  \smallskip
  
  (3.7.9) \ $uzu^* = z + bu^*+ub^*$.
  
    \smallskip

 \n  On a $bu^*+ub^* \in H_o$ , et inversement tout élément de $H_o$ peut s'écrire sous la forme $bu^*+ub^*$. La formule 
  (3.7.9) est donc équivalente à :
  
    \smallskip

  (3.7.10)  \ $uzu^* \equiv z $ \ (mod $H_o$), \ i.e.  $q_{uzu^*} = q_z$,
  
    \smallskip

\n  ce qui équivaut à dire que $u$ fixe la forme quadratique $q_z$.

  \smallskip
  
  \n {\it Remarque}. On pourrait supprimer l'hypothèse (3.7.2), et supposer seulement que  $\r^2=0$ et $R/\r = E$.
  On trouve alors que {\it l'image de $U_R \to \Sp(V)$ est une intersection de groupes orthogonaux} ${\bf O}(q_i)$ 
  relatifs à des formes $q_i$ dont les formes bilinéaires associées sont des multiples de $C$. En prenant $R$ convenable, on peut s'arranger pour que cette intersection soit triviale,  autrement dit que {\it tous les éléments unitaires de $R$ sont} $\equiv 1$ (mod $\r$).

                                   \bigskip
                      
 \centerline    {\bf \S4 - Nullité de $H^1(k,U^0_A)$}
     
     \bigskip
     
     \smallskip

Le but de ce \S \ est de démontrer le théorème C de l'introduction. Rappelons son énoncé:

\smallskip

                  \n {\bf Théorème C} - {\it Soit $k$ un corps parfait de caractéristique $2$, et soit $A$ une $k$-algèbre à 
                  involution de dimension finie. Si $U_A^0$ est la composante neutre du groupe unitaire $U_A$ de $A$, on a
                  $H^1(k,U^0_A) = 0$.}
                  
                  \bigskip

  La démonstration sera donnée au \S4.9. Elle procède par une série de réductions, basées sur le \S3, ainsi que sur des résultats de cohomologie galoisienne qui sont rappelés aux  \S\S \ 4.2 à 4.6.
  
  \bigskip
  
  \bigskip
  
     \n 4.1. {\bf Notations.}
     
     \smallskip
  Dans la suite de ce \S, $k$ est un corps parfait de caractéristique 2; on note $\kbar$ une clôture algébrique de $k$; le
  groupe de Galois Gal($\kbar/k$) est noté $\Gamma_k$.  Si $G$ est un groupe algébrique sur  $k$, on note 
  $H^1(k,G)$ l'ensemble de cohomologie  $H^1(\Gamma_k, G(\kbar))$, cf. e.g. [Se 64, I, \S5 et III, \S1]; la classe    
  dans $H^1(k,G)$ du cocycle unité est notée 0.

 \medskip
  
\n 4.2. {\bf Une première réduction.}  

\smallskip
\n {\bf Lemme 4.2.1} - {\it Si $G' \to G$  est un homomorphisme de groupes algébriques tel que $G'(\kbar) \to G(\kbar)$ soit bijectif,
 l'application correspondante  $H^1(k,G') \to H^1(k,G)$ est bijective. }
 
 \smallskip
   C'est clair.
   
   \smallskip
   
   \n En particulier:
   
   \smallskip
   
   \n {\bf Lemme 4.2.2} - {\it L'application $H^1(k,G^{\rm red})  \to H^1(k,G)$ est bijective.}
   
   \smallskip
   
                      \medskip
                 
                 \n 4.3. {\bf Quotient par un sous-groupe unipotent.}
                 
                 \smallskip
                                
                                \n {\bf Lemme 4.3.1} - {\it  Si $N$ est un groupe  unipotent connexe, on a $H^1(k,N)=0$}.
                                
                                \smallskip
  \n {\it Démonstration.} Lorsque  $N$ est lisse, c'est un résultat bien connu (que l'on démontre en prouvant que  $N$
  a une suite de composition dont les facteurs sont isomorphes au groupe additif $\G_a$, cf. e.g. [Se 64, III, prop.6]). Le cas général s'en déduit grâce  au lemme 4.2.2.
  
  \smallskip          
                  
                    \n {\bf Proposition 4.3.2} - {\it Soit $N$ un sous-groupe unipotent normal connexe d'un groupe algébrique $G$. L'application naturelle $\pi : H^1(k,G) \to H^1(k,G/N)$ est injective.}
                    
                    \smallskip
                    
                     \n {\it Démonstration.} Soit  $x \in H^1(k,G/N)$. On doit montrer que $\pi^{-1}(x)$ a au plus un élément.
                     Lorsque $x=0$, cela résulte de la suite exacte de cohomologie non abélienne ([Se 64, I, prop.38]) et du lemme 4.3.1. Le cas général se ramène à celui-là par ``torsion'', cf. [Se 64, I, cor.2 à la prop.39].
                     
                     \smallskip
                     
                     {\small L'application $\pi$ est en fait {\it bijective}. La surjectivité est vraie même si $N$ n'est pas connexe. On la démontre en se ramenant au cas où $N$ est commutatif, et en utilisant la proposition 41 de [Se 64, I] combinée avec le fait que $H^2(k,N)=0$. Comme nous n'utiliserons pas ce résultat, nous laissons les détails de la démonstration au lecteur. Voir aussi [GM 13, lemme 7.3],  qui ne fait pas d'hypothèse sur $k$, mais suppose que $N$ est {\it scindé} (``split'', cf. [Sp 98, chap.14]).}              
                     
                          \medskip
                     
                     \n 4.4. {\bf  Extensions quadratiques.}

                     \smallskip
                     
                        \n {\bf Proposition 4.4.1} ([Se 64, III, \S2.3, exerc.2]) - {\it Soit $k'$ une extension quadratique de $k$, et soit $G$ un groupe algébrique linéaire connexe. L'application naturelle $\iota: H^1(k,G) \to H^1(k',G)$ est injective.}
                        
                        \smallskip  
                        
                        \n {\it Démonstration.}  Si $k$ est fini, on sait, d'après Lang ([La 56]), que $H^1(k,G) = 0$. On peut donc supposer que $k$ est infini.  Vu le lemme 4.2.2, on peut supposer que $G$ est lisse; d'après la proposition 4.3.2, appliquée au radical unipotent de $G$, on peut aussi supposer que  $G$ est réductif. D'autre part, l'argument de ``torsion" utilisé plus haut montre qu'il suffit de prouver
                        que $x \in H^1(k,G)$ et $\iota(x) = 0$ entraîne $x=0$. 
                        
                        L'hypothèse $\iota(x) = 0$ entraîne que  $x$  provient d'un cocycle du groupe Gal$(k'/k)$, agissant sur $G(k')$. Si  $z \ \mapsto \ ^sz$ désigne la 
                       $k'/k$-conjugaison dans $G(k')$, un tel cocycle équivaut à la donnée d'un élément $a$ de $G(k')$ tel que
                       $a.^sa=1$, et il nous faut prouver qu'il existe $b\in G(k')$ tel que $a =b.^sb^{-1}.$ Comme $k$ est infini, 
 la proposition 3.2.1 de [Se 62] montre qu'il existe $z\in G(K')$ tel que l'élément $a' =z^{-1}a \ ^sz$ soit un élément semi-simple {\it régulier} de  $G(k')$, autrement dit appartienne à un tore maximal  
                       de  $G$ et à un seul, cf. [Bo 91, \S12.2]). Quitte à remplacer  $a$ par $a'$, on peut donc supposer que $a$ est semi-simple régulier. Or, on a le lemme suivant:
                       
                       \smallskip
                       
                       \n {\bf Lemme 4.4.2} - {\it Soit  $F$  un corps parfait de caractéristique $ p > 0$, et soit $H$ un $F$-groupe réductif. Si $x \in H(F)$ est semi-simple régulier, il existe un unique élément semi-simple $y$ de  $H(F)$ tel que
                       $x=y^p;$ de plus $y$ est régulier}.
                       
                       \smallskip
                       
                       \n {\it Démonstration du lemme 4.4.2.} Soit $T$ le tore maximal de $H$  contenant  $x$. Le groupe  $T(F)$ est un groupe abélien dans lequel l'application $t \mapsto t^p$ est bijective (cela résulte par descente galoisienne du cas où $F$ est algébriquement clos). Il existe donc un unique
                       $y \in T(F)$ tel que $y^p=x$.  Soit $y'$ un élément semi-simple de $H(F)$ tel que $y'^p=x$; tout tore maximal contenant $y'$ contient $x$, donc est égal à $T$, d'où $y'=y$, ce qui montre à la fois que  $y$  est unique et qu'il est régulier.
                       
                       \smallskip
                       
                       \n {\it Fin de la démonstration de la proposition 4.4.1.} D'après le lemme ci-dessus, appliqué à  $F=k',  p=2, H=G_{/k'}, x=a$, il existe un unique élément semi-simple $b$ de $G(k')$ tel que $b^2=a$; on a $(^sb)^2 =  \ ^sa = a^{-1}
                       = (b^{-1})^2$, d'où $^sb = b^{-1}$ et 
                       $a = b.b = b.^sb^{-1}$, comme on le désirait.   
                       
                       \medskip                    
                       
                       \n 4.5. {\bf Le cas des groupes linéaires, symplectiques et orthogonaux.}
                       
                       \smallskip
                       
                       \n {\bf Proposition 4.5.1 -} {\it On a $H^1(k,G)=0$ lorsque $G$ est l'un des groupes suivants}:
                       
                       (a) {\it le groupe multiplicatif $\GL_{1,S}$ d'une $k$-algèbre $S$  de dimension finie}.
                                              
                       (b) {\it le groupe symplectique $\Sp_{2n}$ défini par une forme bilinéaire alternée non dégénérée de rang} $2n$, $n \geqslant 1$;
                       
                       (c) {\it le groupe spécial orthogonal $\SO(q)$ associé à une forme quadratique non dégénérée  $q$ de rang pair.}
                       
 \smallskip
 
 \n {\it Démonstration.} L'assertion (a) équivaut à dire que, si un $S$-module devient isomorphe à $S$ après extension du corps de base à $\kbar$, alors il est isomorphe à $S$, ce qui est un résultat standard sur les modules (non nécessairement libres), cf. [A VIII, \S2, th.3]. Pour (b) voir par exemple [KMRT 98, (29.25)]. 
 
 \smallskip
 
 D'après 
  [KMRT 98, (29.29)]   l'assertion (c) signifie que deux formes quadratiques non dégénérées de même rang pair et de même invariant d'Arf (appelé ``discriminant'' dans [KMRT 98, xix-xxi]) sont isomorphes, ce qui résulte facilement de l'hypothèse que  $k$ est parfait, cf. [Arf 41].

                                          \n {\small [Variante: utiliser la proposition 4.4.1 pour se ramener au cas où le corps $k$  n'a aucune extension quadratique, et prouver  qu'alors toutes les formes quadratiques non dégénérées sont hyperboliques.]} 
                                            \medskip

                       \n 4.6. {\bf Restriction des scalaires.}
                       
                       \smallskip
                       
                       Soit $K$ une extension finie de $k$. Si $X$ est une $K$-variété quasi-projective, on lui associe (cf. e.g. [BoS 64, \S2.8 à \S2.10], ou [KMRT 98, \S20.5 à \S20.9]) une $k$-variété quasi-projective $Y$, munie d'un $K$-morphisme $p: Y_{/K} \to X$. Le couple $(Y,p)$ est caractérisé par les propriétés équivalentes suivantes:
              
              \smallskip         
                       (a) (à la Weil, cf. [We 61, \S 1.3]) Soit $\Sigma$ l'ensemble des $k$-plongements de $K$ dans $\kbar$; si $\sigma \in \Sigma$, soit $X_\sigma$ la $\kbar$ variété déduite de $X$ par l'extension des scalaires $\sigma: K \to \kbar$. Les conjugués $p^\sigma: Y_{/\kbar} \ \to X_\sigma $
                       de $p$ définissent un {\it isomorphisme}
                  
                  \smallskip     
                       
                         (4.6.1)  $Y_{/\kbar} \ \to \ \prod_{\sigma \in \Sigma} X_\sigma.$
                         
                         \smallskip
                         
                         (b) (à la Grothendieck, cf. [CGP 10, A5]) Le couple $(Y,p)$ représente le foncteur qui, à une $k$-variété $Z$,
                         associe l'ensemble Mor$_K(Z_{/K}, X)$ des $K$-morphismes de $Z_{/K}$ dans $X$. On a donc une bijection naturelle:
                         
                         \smallskip
                         
                           (4.6.2)  \  Mor$_k(Z,Y) \ = \ $ Mor$_K(Z_{/K}, X)$.
                           
                           \smallskip
                           
                           \n Pour $Z = \ $Spec$(k)$, cela donne:
                           
                           \smallskip
                           (4.6.3)  $Y(k) \ = \ X(K)$.
                           
                           \smallskip
                           
                           On dit que  $Y$ se déduit de $X$ par {\it restriction des scalaires de $K$ à $k$}, et l'on écrit 
                           $Y =R_{K/k} (X)$.
                           
                           \smallskip
                           
                           Soit $X$ un $K$-groupe algébrique; alors  $R_{K/k} (X)$  a une structure naturelle de $k$-groupe algébrique (car le foncteur  $R_{K/k} $ commute aux produits), et l'on a (cf. [BoS 64, cor.2.10] et [KMRT 98, (29.6)]):
                           
                           \smallskip
                           
                           (4.6.4)    $H^1(k,  R_{K/k} (X)) \ = \ H^1(K, X)$.   
                           
                           \smallskip
                           
                           De plus:
                           
                           \smallskip
                           
                           \n {\bf Proposition 4.6.5} - {\it Soit  $Y'$    un $k$-sous-groupe algébrique de $Y =   R_{K/k} (X)$. Pour que
                           $Y' $ soit de la forme  $R_{K/k} (X')$, où  $X'$ est un $K$-sous-groupe algébrique de $X$, il faut et il suffit que $Y'_{/k}$ soit compatible avec la décomposition en produit (4.6.1), i.e. soit un produit de sous-groupes algébriques
 des}    $X_\sigma$.
 
 \n {\small [Il y a un énoncé analogue pour les variétés qui ne sont pas munies d'une structure de groupe.]}
                             
                             \smallskip

                             \n {\it  Démonstration}. Voir [BT 65, 6.18] , qui fait des hypothèses de connexion (et - implicitement - de lissité) qui sont inutiles.  [Vu (4.6.1), le point essentiel est la remarque que, si  $Z = \prod_{i\in I} Z_i$, et si $Z'$ est un sous-schéma fermé non vide de $Z$ qui est décomposable en $Z' = \prod_{i \in I} Z'_i$, avec $Z'_i \subset Z_i$, alors une telle décomposition est {\it unique}, i.e. les $Z'_i$ sont déterminés de manière unique par  $Z$. Noter que ce ne serait pas vrai si $Z$ était vide, car l'ensemble vide se décompose en produit de plusieurs façons.]
                             
                             \smallskip
                             
                             \n {\bf Corollaire 4.6.6} - {\it Soit $X$ un $K$-groupe réductif, et soit $V$ un $k$-sous-groupe réductif de $Y =   R_{K/k} (X)$. Si le rang de $V$ est égal au rang de $Y$, il existe un $K$-sous-groupe $W$ de $X$ tel que
                             $V = R_{K/k}(W)$.}
                             
                             \smallskip
                             
                             \n {\it Démonstration}. Cela résulte de la proposition ci-dessus, combinée avec le corollaire 2.6.9.

                                                  \medskip
                                                  
                                                  {\it Exemple de restriction des scalaires.} Soit $A$ une $K$-algèbre à involution sur $K$,
                                                  et soit  $U_{A,K}$ (resp.  $U_{A,K}^{\sch}$) son groupe unitaire (resp. son groupe unitaire schématique) sur  $K$. Notons $U_A$ (resp.  $U_A^{\sch}$) les groupes analogues sur  $k$. On a
                                                  
                                                  \smallskip
                                                  
                   (4.6.7)                               $U_A = R_{K/k}(U_{A,K})$ \ et $U_A^{\sch} = (U_{A,K}^{\sch})$.
                   
                   \smallskip
                   
                   \n Cela résulte facilement de la définition (b) ci-dessus (en comparant les points à valeurs dans une $k$-algèbre commutative). 
                   
                     D'après (4.6.4), on en déduit :
                     
                     \smallskip
                     
                     (4.6.8) $ H^1(k, U_A) = H^1(K, U_{A,K})$.

                                                 \medskip

                       \n 4.7. {\bf Le groupe unitaire d'une algèbre à involution semi-simple.}
                       
                       \smallskip

                   Nous allons démontrer un cas particulier du théorème C:
                   
                   \smallskip
                   
                   \n {\bf Proposition 4.7.1} - {\it Soit $A$ une $k$-algèbre à involution de dimension finie satisfaisant aux deux conditions suivantes}:
                   
                    (a) {\it Elle est semi-simple}.

                     (b) {\it Elle est engendrée par $U^0_A$}.
                                         
           \n          {\it On a alors $H^1(k,U_A) = 0$.}                   

                     \n [La condition (b) signifie que (3.3.4) est satisfaite après extension des scalaires à $\kbar$. Noter aussi que $U_A$ est réductif, donc connexe,  cf. \S3.4.]
                     
                     \smallskip
                     
                     \n {\it Démonstration}. On peut supposer que $A$ est indécomposable comme algèbre à involution, autrement dit, est, soit le produit de deux algèbres simples $S$ et $S'$  échangées par l'involution, soit une algèbre simple. Dans le premier cas, on a  $U_A \simeq \GL_{1,S}$, d'où $H^1(U_A) = 0$ d'après la proposition 4.5.1 (a). Dans le second cas, notons $L$ le centre de $A$; c'est un corps, sur lequel agit l'involution $x \mapsto x^*$. Il y a deux 
                    possibilités:
                     
                     \smallskip
                     
                     (i) {\it L'involution agit trivialement sur } $L$ (involution ``de première espèce''). Si $[A]$ désigne la classe de $A$ dans le groupe de Brauer de $L$, on a $2 [A] = 0$ car $A$ est isomorphe à son opposée; puisque $L$ est parfait de caractéristique~2, la 2-dimension cohomologique de $\Gamma_L$
                     est $\leqslant 1$, cf. [Se 64, II,  prop.3] et la 2-composante de son groupe de Brauer est 0. On a donc $[A] = 0$, autrement dit $A$ est isomorphe à une algèbre de matrices $\M_n(L)$, avec $n \geqslant 1$. Distinguons deux cas:
                     
                     \smallskip
                     
          ({\rm i}') On a $n=1$, i.e.      $A = L$. Comme l'involution est triviale sur  $L$, on a $U_A = 1$ d'où $H^1(k,U_A) = 0$.   
                           
                                          \smallskip
                        ({\rm i}")                On a $n > 1$ et dans ce cas $n$ est pair, et l'involution est définie par une forme alternée non dégénérée (cf. \S3.4). Sur le corps de base $L$, le groupe unitaire correspondant $U_{A,L}$ est $\Sp_n$; on en déduit que $H^1(k,U_A) = H^1(L,\Sp_n) = 0$ d'après (4.6.8) et la proposition 4.5.1 (b).
                     
                     \smallskip
                     
                (ii) {\it L'involution agit non trivialement sur $L$}   (involution ``de seconde espèce"). Soit $K$ le sous-corps
                de $L$ fixé par l'involution. L'extension $L/K$ est une extension quadratique. L'algèbre $B = L \otimes_K A$   est produit de deux algèbres simples permutées par l'nvolution. On a vu plus haut que cela entraîne $H^1(L, U_{A,K}) = 0$; d'après
    la proposition 4.4.1 on a donc $H^1(K, U_{A,K}) = 0$, d'où $H^1(k,U_A) = 0$  d'après (4.6.8).

                     \medskip
                     
                 \n 4.8. {\bf  Nullité  de la cohomologie de certains sous-groupes de $U_A$.}
                 
                 \smallskip
    Le résultat suivant généralise la proposition 4.7.1:
   
   \smallskip
   
   \n {\bf Proposition 4.8.1} - {\it Soit $A$ une $k$-algèbre à involution semi-simple de dimension finie, et soit $H$ un sous-groupe réductif de $U_A$ satisfaisant aux deux conditions suivantes}:
   
     (c) {\it Le rang de $H$ est égal à celui de $U_A$.}
     
 (d) {\it L'algèbre $A$ est engendrée par   $H$.}
 
 \n {\it Alors $H^1(k,H) = 0.$}
 
 [Noter que (d) entraîne la condition (b) de la proposition 4.7.1, puisque  $H$  est contenu dans $U^0_A$.]
 
 \smallskip
 
 \n {\it Démonstration}.  Décomposons $A$  en produit d'algèbres à involution indécomposables, comme au  \S  \ précédent. D'après le théorème 3.5.1, $H$ est compatible avec cette décomposition sur $\kbar$, donc aussi sur $k$,  et l'on est ramené au cas où $A$ est indécomposable. Faisons cette hypothèse, notons  $L$ le centre de $A$, et notons $K$ la sous-algèbre de $L$ fixée par l'involution. On a vu au \S4.7 que  $K$ est un corps,  et que  $L$  est étale sur $K$ de degré 1 ou 2.  Commençons par le cas particulier où $K=k$:
 
 \smallskip
 
 \n {\bf Lemme 4.8.2} - {\it  Supposons que  $K=k$. On a alors, ou bien $H = U_A$, ou bien $U_A \simeq \Sp_n$ et 
 $H \simeq \SO(q)$, où  $q$  est une forme quadratique non dégénérée de rang}  $n \ pair  \geqslant 4$.
 
 \n [Noter que cela entraîne $H^1(k,H) = 0$ d'après la proposition 4.7.1 et la proposition 4.5.1 (c).]
 
 \smallskip
 \n {\it Démonstration}. Pour prouver l'égalité   $H = U_A$, il suffit de la démontrer après extension des scalaires à $\kbar$; or cela  a été fait dans le théorème 3.5.1,  à la seule exception du cas où $U_{A  /\kbar}$ et $H_{/\kbar}$ sont respectivement isomorphes à $\Sp_n$ et $\SO_n$ avec $n$ pair $\geqslant 4$. Dans ce cas, on a vu au \S  \ précédent que $U_A$ est $k$-isomorphe à $\Sp_n$.
 Si $B$ est une forme alternée non dégénérée invariante par $U_A$, il existe une unique $\kbar$-forme quadratique $q$
 invariante par $H$ telle que  $q(x+y)=q(x)+q(y) + B(x,y)$; comme cette forme est unique, elle est définie sur $k$,
 et les groupes  $H$ et $\SO(q)$ deviennent égaux sur  $\kbar$, donc sont égaux sur  $k$.
 
 \smallskip
 
   Revenons au cas général où $K$ est une extension finie quelconque de $k$. D'après (4.6.7), on a $U_A =  R_{K/k}(U_{A,K})$. Si l'on étend les scalaires à $\kbar$,  $A$ se décompose en produit d'algèbres correspondant aux plongements de $K$ dans $\kbar$, et le théorème 3.5.1 montre que  $H$  est compatible avec cette décomposition.  D'après le lemme 4.6.5, cela entraîne que  $H$ est de la forme $ R_{K/k}(H')$, où  $H'$  est un $K$-sous-groupe réductif  de   $U_{A,K}$, satisfaisant aux propriétés (c) et (d) sur le corps  $K$. En appliquant à $H'$   le lemme ci-dessus, on voit que,
   ou bien  $H = U_A$, ou bien $H \simeq R_{K/k}(\SO(q)) $ et dans les deux cas, on a $H^1(k,H) = H^1(K,H') = 0$.

                     \medskip
                     
            \n  4.9. {\bf Démonstration du théorème C.}
            
            \smallskip
                              
                              Soit $B$ la sous-algèbre de $\kbar \otimes_k A$ engendrée par les $\kbar$-points de $U^0_A$.   Comme cette algèbre est stable par le groupe de Galois $\Gamma_k$, elle provient par extension des scalaires d'une sous-algèbre à involution $A'$ de $A$. On a $U^0_{A'} = U^0_A$. Cela nous permet de remplacer $A$ par $A'$. Autrement dit, {\it nous pouvons supposer que  la condition d'engendrement $(3.3.4)$ est satisfaite sur} $\kbar$.   
                              
                                Soit $\r$ le radical de $A$, et soit $\pi : U^0_A \to U_{A/\r}$ l'homomorphisme défini par la projection 
                                $A \to A/\r$. Soit $N$ le noyau de $\pi$: comme au \S3.6, on a                                une suite exacte:
                                
                                \smallskip
                                
                                (4.9.1) \ $ 1 \ \to N \ \to U^0_A \ \to \ H \ \to \ 1.$
                                
                              \smallskip
                              D'après le lemme 3.6.3,  $N$ est unipotent connexe. D'autre part $H$ satisfait aux conditions (c) et (d) de la proposition 4.8.1, relativement à l'algèbre à involution $A/\r$: en effet, cela a été démontré au \S3.6 sur $\kbar$.
                              D'après la proposition (4.8.1) on a donc $H^1(k,H) = 0$, d'où $H^1(k,U^0_A)=0$ d'après la proposition 4.3.2.
                              
                              \smallskip
                              
                              Cela achève la démonstration du théorème C.
                              
                              \medskip
                              
                              \n {\bf Corollaire 4.9.2} - {\it L'application $H^1(k, U_A) \ \to \ H^1(k,U_A/U^0_A)$  est injective.}
                              
                              \smallskip
                              
                \n              {\it Démonstration}. Le théorème C entraîne que l'image réciproque de  0  est $\{0\}$. L'injectivité en résulte par ``torsion", cf. [Se 64, I, prop.39, cor.2].
                
                \smallskip
                
                            \n {\bf Corollaire 4.9.3} - {\it Si $k_1$ est une extension finie de degré impair de $k,$ l'application
                             $ H^1(k,U_A) \ \to \ H^1(k_1,U_A)$ est injective.}

                                                             \smallskip
                               
                               \n {\it Démonstration}. D'après le corollaire précédent, il suffit de démontrer l'injectivité de                               
                                  $ H^1(k,U_A/U^0_A) \ \to \ H^1(k_1,U_A/U^0_A)$; or celle-ci résulte de ce que le groupe des $\kbar$-points de $U_A/U^0_A$ est un 2-groupe abélien, cf. corollaire 3.3.3.

                                      \bigskip
                      
 \centerline    {\bf \S5 - Le groupe unitaire de l'algèbre d'un groupe fini}
     
     \smallskip

Dans ce \S, ainsi que dans les deux suivants,  $G$ est un groupe fini, et $k$ est un corps de caractéristique 2. On note  $A$ l'algèbre $k[G]$ du groupe $G$; on munit $A$ de son involution canonique, caractérisée
par le fait que $g^* = g^{-1}$ pour tout $g\in G$.

                  \medskip
                  
                  \n 5.1. {\bf Le groupe unitaire $U_G$.}
                  
                  \smallskip
 Soit $U_G^{\sch}$ le groupe unitaire schématique de $A$; c'est un schéma en groupes dont le groupe des $k$-points  contient $G$. 

Si $E$ est un groupe réduit à un élément, on a $U_E^{\sch}=\m$. Les homomorphismes évidents $E \to G \to E$ donnent
des homomorphismes $ \m \ \to \ U_G^{\sch} \ \to \m$ dont le composé est l'identité. Nous noterons $U_G$ le noyau de $U_G^{\sch} \to \m$; on a une décomposition de $U_G^{\sch}$ en
produit:

                  \medskip

(5.1.1)   \ $U_G^{\sch} \ = U_G \ \times \ \m.$

\medskip

\n {\bf Théorème 5.1.2} (Merkurjev) - {\it Le schéma en groupes $U_G$ est lisse.}

%\n[C'est donc le schéma réduit de $U_G^{\sch}$, autrement dit le schéma
%noté $U_A$ au $\S4$. C'est lui que nous appellerons le {\it groupe unitaire de $G$.]

\medskip

\n {\it Démonstration} (d'après une lettre de Merkurjev du 19/5/2002). Soit $k'$ une $k$-algèbre commutative. Un élément $\sum x_g g$
de $k' \otimes_k A = k'[G]$ est un point de $U_G$ si et seulement si il satisfait aux équations 
suivantes:

\smallskip
 $(A_1) \quad \sum_{g\in G} x_g = 1$,
 
 \smallskip
 
 $(A_s) \quad \sum_{g\in G} x_g x_{sg} = 0   $ pour tout $s \in G$ tel que $s\ne 1$.
 
 \medskip
 
 Lorsque $s$ est d'ordre 2, la relation $(A_s)$ se réduit à $0 = 0$,
 car les termes  $x_gx_{sg}$ et $x_{sg}x_g$ ont pour somme $0$. Un argument
 analogue montre que, pour tout $s$, l'équation $(A_s)$ est équivalente
 à l'équation $(A_{s^{-1}})$.
 
 Soit $G_2$ l'ensemble des éléments $s$ de $G$ tels que $s^2=1$ et
 choisissons une partie $\Sigma$ de $G \sm G_2$ telle que $G \sm G_2 = \Sigma\ \sqcup \ \Sigma^{-1}$; pour tout élément $s$ de $G \sm G_2$, on a, soit $s\in \Sigma$, soit $s^{-1}\in \Sigma$, mais pas les deux à la fois. Le schéma en groupes $U_G$ est défini par l'équation linéaire $L = 1$ et par les équations quadratiques
 $P_s = 0$ pour $s\in \Sigma$, où $L = \sum_{g\in G} x_g$ et $P_s = \sum_{g\in G} x_g x_{sg}$.  Les différentielles  des polynômes $L$ et $P_s$ en l'élément neutre $1 = (1,0,0,...)$ sont:
 
 \smallskip
   
  $dL = dx_1$  \  et  \  $dP_s = dx_s + dx_{s^{-1}}$.
  
  \smallskip
  
  Elles sont linéairement indépendantes. D'après le critère jacobien, 
  cela montre que $U_G$ est lisse au point 1, donc lisse partout puisque
  c'est un schéma en groupes. 
  
  \medskip
  
  \n {\it Autre démonstration.}                 
 Plaçons-nous sur le corps de base $\F_2$; comme $\F_2$ est parfait, on peut définir  le groupe algébrique
      noté $U_A$ au $\S4$; c'est un groupe lisse.  On a $U_A \subset U_G$.   D'après la proposition 3.3.1, l'algèbre $\Lie(U_A)$ contient
      les éléments $1+g$ ($g$ d'ordre 2) et $s+s^{-1}$ ($s \in S$) ; elle
      coïncide donc avec $\Lie(U_G)$, ce qui entraîne\footnote{Plus généralement, on a le critère de lissité suivant: soit $X$ un schéma de type fini sur un corps infini $k$, soit $x\in X(k)$, et soit $T_x(X)$ l'espace tangent à $X$ au point $x$. Supposons que, pour tout $t \in T_x(X)$, il existe un $k$-schéma lisse $V$, un point $v\in V(k)$ et un morphisme $f:V \to X$ tel que $f(v)=x$ et que $t$ appartienne à l'image de $T_v(V) \to T_x(X)$. Alors $X$ est lisse en $x$; cela se voit en prouvant que le cône tangent à $X$ en $x$ est égal à $T_x(X)$.}   que $U_G = U_A$, donc que $U_G$ est lisse.
    
    \medskip

\n {\it Remarque.} Soit $r =   \frac{1}{2} (|G| - |G_2|).$ La démonstration de Merkurjev donnée plus haut montre que $U_G$ est une intersection transversale de $r$ quadriques dans un espace affine. On en déduit, grâce au théorème de Bézout (cf. e.g. [Fu 84, Example 8.4.6]), que {\it le nombre de composantes connexes de $U_G$ est au
plus égal à} $2^r$, l'égalité n'étant possible que si  $U^0_G$ est une variété linéaire ; cette borne est certainement grossière; il serait intéressant de l'améliorer.

\medskip
  
  \n {\it Compléments}.
  
    La décomposition (5.1.1) montre que $\Lie(U_G)$
  est formée des éléments $\sum a_g g$ de $A$ qui sont hermitiens et  tels que
  $\sum a_g = 0.$ On obtient une base de cette algèbre en
  prenant les éléments  $1+g$ ($g$ d'ordre 2) et $s+s^{-1}$ ($s \in \Sigma$). Le nombre de ces éléments est
   $ |G_2| - 1 + |\Sigma|  $.  Comme $|G| - |G_2| = 2 |\Sigma|$, on en tire:

  \medskip
  
     \n {\bf Proposition 5.1.3} - {\it On a}:
     
     \smallskip
  
\quad   \ $ \dim U_G \ = \ \dim \Lie(U_G) \ =  \ |G_2| - 1 + |\Sigma|  \ = \ \frac{1}{2} (|G| + |G_2|)- 1$.

                  \medskip

         Notons $U^0_G$ la composante neutre de $U_G$. Si $H$
                  est un sous-groupe de $G$, on a $U_H^0 \subset U_G^0$.
                  
                  \medskip
                  
                  \n {\bf Proposition 5.1.4} - {\it Soit $(H_i)$ une famille de sous-groupes de $G$ de réunion égale à $G$. Alors $U^0_G$ est engendré} (comme groupe algébrique)
                  {\it par les} $U^0_{H_i}$.
                  
                  \smallskip
                  
                  \n {\it Démonstration.} Cela résulte du fait que $\Lie(U_G)$
                  est engendré comme espace vectoriel par les sous-espaces $\Lie (U_{H_i})$.
                  
                  \smallskip
                    Cet énoncé s'applique en particulier  à la famille des sous-groupes cycliques maximaux de $G$.
                    
                    \smallskip
                    
                    (5.1.5) Rappelons, cf. [CR 62, \S55] et [Fe 82, \S I.7],  que les facteurs indécomposables{ de $A$ sont appelés les {\it blocs}\footnote{Les blocs peuvent être vus, soit comme des idéaux bilatères de $A$, soit comme des algèbres quotients de $A$; dans ce qui suit, nous choisirons le point de vue ``quotients''.} de $A$.  Parmi ceux-ci figure le {
                    \it bloc principal} $B$, caractérisé par le fait que la représentation unité  $ \t: A \to k$ se factorise en $A \to B \to k$. Ce bloc est auto-dual,
                   i.e.  stable par l'involution de $k[G]$. Si $B'$ est le produit des autres blocs, $B'$ est également auto-dual, et $A = B \times B'$. D'où une décomposition du groupe unitaire:
                   
                   \smallskip

               (5.1.6)  $U_G^{\sch} = U_B^{\sch} \times U_{B'}^{\sch}$ \ et  \ $U_G = U_B \times U_{B'}$.

               \smallskip
               
               Comme $\t$ se factorise par $B$, l'homomorphisme $U_G^{\sch} \to \m$ se factorise par $U_B^{\sch}$;   il en résulte que     $U_{B'}^{\sch}= U_{B'}$, autrement dit {\it que $U_{B'}$ est lisse} (c'est une intersection transversale de quadriques).      
                                              \medskip
                  
                  \n 5.2. {\bf Les caractères essentiels.}
                  
                  \smallskip
                  
                    Soit $C = \{1,c\}$ un groupe d'ordre 2. L'homomorphisme  $ t \mapsto 1 + t(1+c)$ est un isomorphisme du groupe additif $\G_a$ sur le groupe $U_C$. Dans cet isomorphisme, le sous-groupe $\{0,1\}$ de $\G_a$ correspond au sous-groupe $C$ de $U_C$.
                    
                    \smallskip
                    
                      Soit $\e: G \to C$ un homomorphisme de $G$ dans $C$, et soit $\varphi_\e: U_G \to U_C$ l'homomorphisme correspondant de $U_G$ dans $U_C$.
                      Lorsqu'il existe $g\in G$ avec $g^2=1$ et $\e(g) = c$, l'extension $1 \to \Ker(\e) \to G \to C \to 1$ est scindée, et cela entraîne que $\varphi_\e$ est surjectif. Dans le cas contraire, on a:
                      
                      \smallskip
                      
                      \n {\bf Théorème 5.2.1} - {\it Supposons que $\e(g)=1$ pour tout $g\in G$ tel que $g^2=1$.
                      L'image de  $\varphi_\e: U_G \to U_C$ est alors égale à $\e(G)$, i.e. à $\{1\}$ si $\e = 1$ et à $C$ si $\e \ne 1$. On a $\varphi_\e(U^0_G) = \{1\}.$}
                  
                  \smallskip
                  
                  \n {\it Démonstration.} Le cas où $\e = 1$ est clair. Supposons que $ \e \ne 1$ et notons $G_1$ (resp. $G_c$) l'ensemble des $g\in G$ tels que $\e(g)=1$ (resp. $ \e(g) = c$). L'hypothèse faite sur $\e$ équivaut à dire que, pour tout $g\in G_c$, on a $g \ne g^{-1}$. On peut donc décomposer $G_c$ comme réunion disjointe $G_c = S \sqcup S^{-1}$. Soit $x = \sum x_g g$ un point de $U_G$, à valeurs dans un corps $k$ de caractéristique $2$. On a:
                  
                  \smallskip
                  
                  (5.2.2) \ $\varphi_\e(x) = \lambda(x).1 + \mu(x).c$, \ 
                 avec $\lambda(x) = \sum_{a\in G_1} x_a$ et  $\mu(x) = \sum_{b\in G_c} x_b$.
                  
                  \smallskip
                  
                  \n On a:
                  
                  \smallskip
                  
                  (5.2.3) \ $\lambda(x) + \mu(x) = 1$ d'après la propriété $(A_1)$ du \S5.1.
                  
                  \smallskip
                  
                  \n D'autre part:
                  
                  \smallskip
                  
                  (5.2.4)  \ $\lambda(x)\mu(x) = \sum_{a\in G_1, b\in G_c} x_ax_b.$
                  
                  \smallskip
                  
                  \n Définissons une application $\theta: S \times G \to G_1 \times G_c$ par:
                  
                  \smallskip
                  
                  $(s,g) \mapsto (g,sg)$ \ \ si $g\in G_1$ \ et \ $(s,g) \mapsto (sg,g)$ \ si $g\in G_c$.
                  
                  \smallskip
                  
             \n      On vérifie que $\theta$ est bijective, l'application réciproque $G_1 \times G_c \to S \times G$ étant:
                  
                  \smallskip
                  
      $(a,b) \mapsto   (ba^{-1},a)$  \ si $ba^{-1} \in S$ \ \ et \ \ $(a,b) \mapsto (ab^{-1},b)$ \  si $ab^{-1} \in S$.
      
      \smallskip
      
       Si les couples $(s,g)$ et $(a,b)$ se correspondent par $\theta$ et $\theta^{-1}$, on a $x_ax_b=x_gx_{sg}$. Cela permet de récrire (5.2.4) sous la forme:
       
       \smallskip
       
       (5.2.5) $\lambda(x)\mu(x) = \sum_{s\in S} \sum_{g\in G} x_gx_{sg}.$
       
       \smallskip
       
       D'après la formule $(A_s)$ du \S5.1, on a $\sum_{g\in G} x_gx_{sg} = 0$ pour tout $s \in S$; d'après (5.2.5) on a donc $\lambda(x)\mu(x) =  0$, et comme $\lambda(x)+\mu(x) = 1$ d'après (5.2.3), on en déduit que $(\lambda(x), \mu(x)) = (1,0)$ ou $(0,1)$. L'image de $x$ dans $U_C$ est donc égale à $1$ ou à $c$.
       Comme l'image de $\varphi_\e$ est finie, on a $\Ker(\varphi_\e) \supset U^0_G$.
                  
                  \medskip
                  
                  Un homomorphisme $\e: G \to C$ tel que $\e(g)=1$ pour tout $g\in G_2$ sera appelé un {\it caractère essentiel} de $G$. D'après le théorème ci-dessus l'homomorphisme $\varphi_\e : U_G \to U_C$ est à valeurs dans $C$, autrement dit définit un homomorphisme de groupes algébriques $\tilde{\e}: U_G \to C$; ici, $C$ est vu comme un groupe algébrique étale, de dimension 0. Il est commode d'interpréter $\tilde{\e}$ comme l'homomorphisme $U_G \to \Z/2\Z$ défini par:

                 \smallskip
                    
                  (5.2.6)  \ $\tilde{\e}(x) = \mu(x)$ \ pour tout $x \in U_G(k)$,
                  
                    \smallskip
                  \n ce qui a un sens puisque $\mu(x) = 0$ ou $1$, comme on vient de le voir.
                  
                    Soit $U^0_G$  la composante neutre de $U_G$. L'homomorphisme $\tilde{\e} : U_G \to C \simeq \Z/2\Z$ est trivial sur $U^0_G$, et peut donc être vu comme un caractère du groupe quotient $U_G/U_G^0$.
                    
                    \smallskip
                    
                    {\small {\it Remarque.} L'homomorphisme $k[G] \to k[C]$ défini par $\e$ se factorise par le bloc principal  $B$, cf. (5.1.15). Il en résulte que $\tilde{\e} : U_G \to C$ se factorise par $U_B$. Autrement dit, les caractères essentiels n'apportent aucun renseignement sur les groupes unitaires des blocs distincts du bloc principal.}
                    
                    \medskip 
    \n {\bf Théorème 5.2.7} - {\it Soient $\e_1, \e_2: G \to C$ deux caractères essentiels 
    de $G$  et soit $\e_3= \e_1\e_2$ leur produit. Alors $\e_3$ est essentiel et l'on a}:
    
   (5.2.8)  $\tilde{\e}_3 = \tilde{\e}_1\tilde{\e}_2$; 
   
   (5.2.9) {\it l'homomorphisme $\varphi: U_G \to U_{C \times C}$ défini par $(\e_1,\e_2): G \to C \times C$ est trivial sur $U^0_G$, et son image est contenue dans $C \times C$.}
    
    \smallskip        
    
    Le fait que $\e_3$ soit essentiel est clair. Le reste de la proposition est évident si  l'un des $\e_i$ est égal à 1. On peut donc supposer qu'ils sont $\ne 1$, autrement dit que $G \to C \times C$ est surjectif.
    
    \medskip
     \n {\it Démonstration de} (5.2.8).
  Notons $\mu_i : G \to \Z/2\Z$ la fonction $\mu$ associée à $\e_i$ comme dans (5.2.2):
    
    \smallskip   
    
  (5.2.10)      $\mu_i(x) = \sum_{\e_i(g) = c} x_g$ \ si   $x = \sum x_gg$  est un point de $U_G$.
    
    \smallskip
    
    D'après (5.2.6), la formule (5.2.8) équivaut à   $\mu_1(x)+\mu_2(x)+\mu_3(x) = 0$ 
    pour tout point $x$ de $U_G$. Notons $H_1$ l'ensemble des $g\in G$ tels que
    $\e_1(g)=1$ et $\e_2(g)=\e_3(g)= c$; définissons de même $H_2$ et $H_3$.     
    Posons :
    
    \smallskip
    
    (5.2.11) $Y =  \sum_{g\in H_1} x_g, \ Z =  \sum_{g\in H_2} x_g, \ T =  \sum_{g\in H_3} x_g.$
    
        \smallskip
        
        L'ensemble des $g$
    tels que $\e_1(g)=c$ est $H_2 \sqcup H_3$. Avec les notations de (5.2.10), on a donc:

\smallskip
      (5.2.12)  \ $\mu_1(x) = Z + T.$      
      \smallskip
      
  \n    De même:

    \smallskip
      (5.2.13)  \ $\mu_2(x) = T + Y$.
      
      \smallskip

      (5.2.14)  \ $\mu_3(x) = Y+Z$.
      
      \smallskip
      
 \n     En ajoutant ces trois relations, on obtient l'égalité cherchée: 
 
 \smallskip
 $\mu_1(x)+\mu_2(x)+\mu_3(x) = 0.$
 
 \medskip

   \n {\it Démonstration de } (5.2.9).  Conservons les notations ci-dessus et notons $H_0$ l'ensemble des $g\in G$ tels que $\e_1(g)=\e_2(g)=1$; on a 
   $G = H_0 \sqcup H_1 \sqcup H_2 \sqcup H_3$. Posons:
   
   \smallskip
   (5.2.15) $X = \sum_{g\in H_0} x_g.$
   
   \smallskip
   
   L'image de  $x$  par $\varphi: U_G \to U_{C \times C}$ a pour coordonnées $(X,Y,Z,T)$; de façon plus précise, on a:
   
   \smallskip 
   
    $ \varphi(x) = X.(1,1) \ + \ Y.(c,1) \ + \ Z.(1,c) \ + \ T.(c,c)$,
   
   \smallskip
   
 \n   où $(1,1), (c,1), (1,c), (c,c)$ sont les quatre éléments de $C \times C$. Le fait que  $ \varphi(x) $ appartienne à $C \times C$ équivaut à dire que $X,Y,Z$ et $T$ sont tous 0, à l'exception de l'un d'eux qui est égal à 1. Pour le démontrer, nous allons d'abord calculer  $XY + ZT$. D'après (5.2.11) et (5.2.15), on a:
 
 \smallskip
 
 (5.2.16) \ $XY+ZT = \sum x_ux_v$, 
 
 \smallskip
 \n la somme étant étendue aux $(u,v) \in \ H_0\times H_1 \ \sqcup \ H_2\times H_3.$
 
 Choisissons une partie $\s$ de $H_1$ telle que $H_1 = \s \sqcup \s^{-1}$. C'est possible car $H_1$ est stable par $g\mapsto g^{-1}$, et ne contient aucun élément de carré 1. Si $(u,v) \in \ H_0\times H_1 \ \sqcup \ H_2\times H_3$, on a $vu^{-1} \in H_1$. Posons $(a,b) = (u,v)$ si $vu^{-1} \in \s$
 et $(a,b)=(v,u)$ si $vu^{-1} \in \s^{-1}$. L'application $(u,v) \mapsto (a,b)$ est une bijection de $ H_0\times H_1 \ \sqcup \ H_2\times H_3$ sur l'ensemble des couples $(a,b) \in G \times G$ tels que $ba^{-1} \in \s$. Cela permet de récrire la formule (5.2.16) sous la forme:
 
  \medskip
 \! (5.2.17) $XY + ZT = \sum_{(a,b) \in \s'} x_ax_b = \sum_{a\in G, s\in \s} x_ax_{sa} = \sum_{s\in \s} \sum_{a\in G}x_ax_{sa} $. 
   
   \smallskip
   
   D'après la formule $(A_s)$ du \S5.1, on a $\sum_{a\in G}x_ax_{sa} = 0 $ pour tout  $s\in \s$. On en tire:

   \smallskip
   (5.2.18) \ $XY+ZT = 0$, \ autrement dit \ $XY = ZT$.
  
  \smallskip 
   Un argument analogue montre que $XZ = YT$ et $XT = YZ$. D'autre part la formule $(A_1)$ du \S5.1 montre que $X+Y+Z+T=1$. On conclut la démonstration de (5.2.9) en appliquant le lemme suivant :
   
   \smallskip
   
   \n {\bf Lemme 5.2.19} - {\it Soient $X,Y,Z,T$ des éléments d'un corps de caractéristique $2$ tels que
   $X+Y+Z+T = 1, XY = ZT, XZ = YT$ et $XT = YZ$. Alors trois de ces éléments sont égaux à  $0$, et le quatrième est égal à $1$.}
   
\smallskip
   
   \n {\it Démonstration.} Si aucun des $X,Y,Z,T$ n'est nul, les équations $XY=ZT, XZ=YT$ et $XT=YZ$ entraînent $X^2=Y^2=Z^2=T^2$, i.e. $X=Y=Z=T$, ce qui est incompatible avec $X+Y+Z+T=1$. Si par exemple $X$ est 0, l'équation $XY= ZT$ montre que, soit $Z=0$ ou $T=0$; si $X=Z=0$, l'équation $XZ=YT$ montre que $Y$ ou $T$ est 0; de même, si $X=T=0$, alors $Y$ ou $Z$ est 0. La seule possibilité est donc que trois des éléments soient 0, le quatrième étant égal à 1.
\medskip

\medskip

\n 5.3. {\bf Le sous-groupe $G_0$.}

\smallskip

  Soit $X_G$ le groupe des caractères essentiels de  $G$, et soit $G_0= \bigcap_{\e\in X_G} \Ker(\e).$
  Le groupe $G/G_0$ est un $2$-groupe abélien élémentaire, dont le dual est $X_G$.

  \smallskip

  \n {\bf Théorème 5.3.1} - (i) {\it Le groupe $G_0$ est le sous-groupe de $G$ engendré par les éléments
  d'ordre $2$ et par les carrés.}
  
  (ii) {\it On a $G_0 = G \cap U^0_G.$}
    
    \smallskip
    
  {\small  [Dans (ii), $G_0$ est vu comme sous-groupe algébrique constant de  $U_G$; la formule $G_0 = G \cap U^0_G$ signifie que $G_0 = G \cap U_G^0(k)$ pour tout corps $k$ de caractéristique 2, ce qui équivaut d'ailleurs à
    $G_0 = G \cap U_G^0(
    \F_2)$.]}
      
      \smallskip
      
      \n {\it Démonstration de} (i). Soit $H$ le sous-groupe de  $G$ engendré par les carrés et par les éléments d'ordre $2$ de $G$. Tout caractère essentiel de $G$ est trivial sur $H$. On a donc $H \subset G_0$. D'autre part, $H$ est normal dans $G$ et tout élément de $G/H$ est de carré 1. Il en résulte que $G/H$ est un 2-groupe abélien élémentaire, et $H$ est donc l'intersection des noyaux des caractères $\e: G \to C$ qui sont triviaux sur $H$. Comme $H$ contient $G_2$, ces caractères sont essentiels. On a donc $G_0 \subset H$, d'où $G_0=H$.
                  \smallskip
                  
                  \n {\it Démonstration de} (ii). D'après la proposition 3.3.1 et le corollaire 3.3.2, les éléments d'ordre 2 de $G$, ainsi que les carrés, sont dans $U^0_G$. Vu (i), cela entraîne $G_0 \subset U^0_G$. D'autre part, si $g\in G$ n'appartient pas à $G_0$, il existe un caractère essentiel $\e$ de $G$ tel que
                  $g \notin \Ker(\e)$; comme $\Ker(\e)$ contient $U^0_G$, cela montre que $g$ n'appartient pas à $U^0_G$. On a donc bien $G_0 = G \cap U^0_G$.
                  
                  \smallskip

                    \smallskip
  
        \n {\bf Proposition 5.3.2} -  {\it Pour que $G_0 = G$, il faut et il suffit que $G$ soit engendré par des éléments d'ordre $2$ et par des éléments d'ordre impair.}
 
      \smallskip
      
      \n {\it Démonstration}. Soit $K$ le sous-groupe de $G$ engendré par les éléments d'ordre 2 et par ceux d'ordre impair. Comme tout élément d'ordre impair est un carré, on a $K \subset G_0$; si $K=G$ on a donc $G_0=G$. D'autre part, si $K \ne G$, le quotient $G/K$ est un 2-groupe non trivial; il existe donc un homomorphisme surjectif $ G \to C$ dont le noyau contient $K$; un tel homomorphisme est essentiel, ce qui montre que $G_0 \ne G$.

  \smallskip
           
 \n {\it Exemples de groupes  $G$  tels que} $G = G_0$ : un groupe simple (abélien ou non), un groupe de Coxeter, un groupe commutatif sans élément d'ordre 4.
  
  \smallskip
  
 \n {\it Exemples de groupes $G$ tels que} $G \ne G_0$: un groupe quaternionien d'ordre $2^n \ (n \geqslant 3)$, un groupe ayant un quotient cyclique d'ordre 4, un groupe $\widehat{{ S}}_n$, avec $n = 3,4,5.$ \footnote{Rappelons que $\widehat{{ S}}_n$ désigne une extension centrale de $ S_n$ par un groupe d'ordre 2, dans laquelle les transpositions, et les produits de deux transpositions disjointes, se relèvent en des éléments d'ordre 4.}

  \medskip
                  
 \n 5.4. {\bf Décomposition de $U_G/U^0_G$ en produit.}
 
   Le groupe $U_G/U_G^0$ est un groupe étale. Plaçons-nous d'abord sur le corps algébriquement clos $k = \Fbar_2$, de sorte que nous pouvons identifier         $U_G/U_G^0$ au groupe de ses $k$-points. D'après le corollaire 3.3.3,
   c'est un 2-groupe abélien élémentaire; il contient le groupe $G/G_0$, cf. théorème 5.3.1. 
   
   Soit $\widetilde{X}_G = \Hom(U_G/U_G^0,C)$ le dual de $U_G/U_G^0$. Comme le dual de $G/G_0$ est le groupe $X_G$ des caractères essentiels, l'injection $G/G_0 \to  U_G/U_G^0$ donne par dualité une surjection $r: \widetilde{X}_G \to X_G$. D'autre part,
   on a vu que tout $\e \in X_G$ définit un élément  $\tilde{\e}$ de $\widetilde{X}_G$, et que l'application
   $i :  X_G \to \widetilde{X}_G$ définie par $ \e \mapsto \tilde{\e}$ est un homomorphisme (théorème 5.2.7); le composé $r \circ i$ est l'identité car la restriction de $\tilde{\e}$ à $G/G_0$ est $\e$. On déduit de ceci que le groupe $\widetilde{X}_G$ se décompose en produit :
   
   \smallskip
   (5.4.1) $\widetilde{X}_G = X_G \times \IM(i)$.
   
    \smallskip
    
\n   Par dualité, cela donne une décomposition analogue pour $U_G/U_G^0$:

    \smallskip
   (5.4.2) $U_G/U_G^0 = G/G_0 \times E_G$ , 
   
    \smallskip
   
   \n où $E_G$ est le dual de $\IM(i)$, c'est-à-dire l'intersection des noyaux des  $\tilde{\e}$ (les $\tilde{\e}$
   étant identifiés à des caractères de $U_G/U_G^0$).
   
   Cette décomposition, étant canonique, descend au corps de base $\F_2$. On obtient ainsi:
   
    \smallskip
   \n {\bf Proposition 5.4.3} - {\it Soit $E_G$ l'intersection des noyaux des homomorphismes $\tilde{\e}:U_G/U_G^0 \to C$ associés aux caractères essentiels $\e$ de $G$. Le groupe    $U_G/U_G^0$ est produit direct de ses sous-groupes  $G/G_0$ et $E_G$.}
    
    \smallskip
    \n {\it Remarques.} 
    
    (5.4.4) On aurait également pu définir $E_G$ comme le quotient $U_G/G.U_G^0$; toutefois, le point essentiel de la construction ci-dessus est que $E_G$ est canoniquement un facteur direct du groupe
    $U_G/U_G^0$; ce sera utile au \S6.3.
    
    (5.4.5) Il serait intéressant de trouver un procédé de calcul explicite de $E_G$. On en verra quelques exemples ci-dessous; dans chacun d'eux, on trouve que $E_G = \{1\}$, sauf dans les cas 5.5.6  et 5.5.17 : $G = \widehat{{ S}}_3$ et $G = \widehat{{ A}}_5$.
       \medskip
    
    \n 5.5. {\bf Exemples de détermination des groupes $U_G/U^0_G$ et $E_G$.}
    
    \medskip
      Pour simplifier, on suppose que le corps de base  $k$  est algébriquement clos ; en fait, cette hypothèse est inutile dans les cas 5.5.2 et 5.5.5, et dans les cas 5.5.6, 5.5.16 et 5.5.17 elle peut être remplacée par celle que
      $k$ contient une racine primitive cubique de l'unité.
      
      \smallskip
    \n 5.5.1. {\it $G$ d'ordre impair.} L'algèbre $k[G]$ est semi-simple.  Les 2-blocs de $G$ sont les facteurs simples de $k[G]$. Le bloc principal est
    le corps $k$. Aucun autre bloc n'est auto-dual ; cela résulte du théorème analogue en caractéristique 0, dû à Burnside [Bu 11, \S222, th.II] sous la forme équivalente suivante: le seul caractère irréductible de $G$ à valeurs réelles est le caractère unité. On en déduit, cf. \S3.4, que $U_G$ est isomorphe à un produit $ \prod  \GL_{n_i}$,
    avec $|G| = 1 + 2 \sum n_i^2$; en particulier, $U_G$ est  connexe et l'on a donc $E_G = \{1\}$. 
    
    \medskip    
    
    \n 5.5.2. {\it $G$ cyclique d'ordre une puissance de 2}.  Soit  $n$ l'ordre de $G$. On a vu plus haut que
    $U_G =   \{1\}$ lorsque  $n=1$ et $U_G \simeq \G_a$ lorsque $n=2$. Supposons $n \geqslant 4$, et notons $\e$ l'unique caractère essentiel non trivial de $G$.
    On peut identifier $k[G]$ à l'algèbre $A  = k[t]/(t^n)$, de telle sorte que $s = 1+t$ soit un générateur de  $G$.
    L'involution de  $A$ est telle que  $s^* = s^{-1}$, autrement dit $t^* = t/(1+t) = t+t^2 + ... + t^{n-1}$.
    Tout élément unitaire  $x$ de  $U_A$ s'écrit comme un polynôme:
   
   \smallskip 
    \quad  $x = 1 + a_1t + ... + a_{n-1}t^{n-1}$.
    
    \smallskip 
    
  \n  Dans le produit $x^*x$, le coefficient de $t^2$ est $a_1^2 + a_1$; puisque $x$ est unitaire, on a donc  $a_1^2 + a_1=0$, i.e. $a_1 = 0$ ou $1$. Le caractère $U_G \to \Z/2\Z$ défini par  $ x \mapsto a_1$ est le caractère $\tilde{\e}$ associé à $\e$. Son noyau $U_1$ est formé des $x$ tels que $a_1 = 0$, i.e. $x \equiv 1\pmod{t^2}$; nous allons voir que $U_1$ est connexe.
    
    \smallskip
    \n {\bf Lemme 5.5.3} - {\it Tout élément $x$ de $U_1$ s'écrit sous la forme $x = y^{-1}y^*.z$, où $y$ est un élément inversible de  $A$ et} $z \equiv 1\pmod{t^{n-2}}$.
    
    \n [Noter que tout élément de $A$ de la forme $ y^{-1}y^*$ est unitaire; il en est de même des   $z\in A$
    tels que  $z \equiv 1\pmod{t^{n-2}}$, car on a $z=z^*$ et $z^2=1$.]

\smallskip

\n {\it Démonstration.} Soit $x \in U_1$, et soit $b = t^*+x^*t$. On a:

\smallskip

 \ (5.5.4) $ bx  =  t^*x + xx^*t = t^*x + t =  b^*.$
 
 \smallskip

\n Comme $x \equiv 1\pmod{t^2}$, on a $b \equiv t^*+ t \equiv t^2 \equiv tt^* \pmod{t^3}$. Il existe donc
$y\in A$ tel que $b=tt^*y$ et $y \equiv 1\pmod{t}$. L'équation (5.5.4) entraîne $ tt^*xy = tt^*y^*$, autrement dit
$tt^*(xy+y^*) = 0$, ce qui équivaut à  $xy+y^* \equiv 0\pmod{t^{n-2}}$. Si l'on pose $z= xy(y^*)^{-1}$, on a donc
$z \equiv 1\pmod{t^{n-2}}$ et $x = y^{-1}y^*.z$; d'où le lemme.

\smallskip
  Comme l'ensemble des  $ y^{-1}y^*$ est connexe, ainsi que celui des $z$ avec $z \equiv 1\pmod{t^{n-2}}$,
  le lemme montre que $U_1$ est connexe, et comme il est d'indice 2 dans $U_G$, on a $U_1 = U^0_G$.
 Cela montre que l'on a $E_G =  \{1\}$ (d'ailleurs il est clair que  $U_G$ est engendré par $G$ et $U_1$).    
 
  % Noter que l'exposant de $U^0_G$ est $n/2$, alors que celui de $U_G$ est $n$.      
                  
                  \medskip
                  
                  \n 5.5.5. {\it $G$ quaternionien d'ordre 8.} Les éléments de $G$ sont traditionnellement notés
                  $\{\pm1,\pm i, \pm j,
                  \pm k\}$; pour éviter tout risque de confusion,  dû au corps
                  de base $k$ ainsi qu'au fait que $-1=1$ dans $k$, nous les noterons $\{1,\o,u,\o u, v, \o v,w, \o w \}$; on a
                  $\o^2=1, uv = w, vw=u, wu=v, u^2=v^2=w^2 = \o$.
                     Il y a trois caractères essentiels $\ne 1$, dont les noyaux sont les trois sous-groupes cycliques d'ordre 4 de $G$; on se trouve donc dans la situation du théorème 5.2.7, que l'on va utiliser.
                     
                     Soit $U_1$ l'ensemble des éléments de $k[G]$ de la forme                       $x = 1 + (1+\o)a$, avec $a \in k[G]$.
                     C'est un sous-groupe de  $U_G$: cela résulte du fait que $1+\o$ est central, de carré nul, et que $\o a=\o a^*$ pour tout $a$. Un élément de  $U_1$ s'écrit de façon unique :
                     
                       $ x= 1 + (1+\o)(\alpha + \beta u + \gamma v + \delta w)$, \ \ avec $\alpha, \beta, \gamma, \delta \in k$.
                       
                      \n  Cela montre que $U_1$ est  isomorphe à $\G_a \times \G_a \times \G_a \times \G_a$. D'autre part, $U_1$ est le noyau de l'homomorphisme $U_G \to U_{G/\{1,\o\}}$; d'après le théorème 5.2.7, l'image de cet homomorphisme
                       est d'ordre 4. On conclut de ceci que $U_1 = U^0_G$,  $(U_G:U_G^0) = 4$ et  $U_G = G.U^0_G$, autrement dit $E_G =  \{1\}$.
                  \medskip

          \n 5.5.6. {\it $G = \widehat{{ S}}_3$, d'ordre 12.}  Le groupe $G$ est produit semi-direct d'un groupe 
          cyclique $C_4= \{1,s,s^2,s^3\}$ d'ordre 4, par un groupe cyclique $C_3 =\{1,r,r^2\}$ d'ordre 3, l'action du premier sur le second étant donnée par  $srs^{-1}= r^2$. Le quotient de $G$ par son centre     $\{1,s^2\}$ est
          isomorphe au groupe symétrique  ${ S}_3$. 
          
          \smallskip
                    
            L'algèbre $A=k[G]$ se décompose en produit de deux blocs $B$ et $B'$. On a $B = A/\pi A$ et $B'= A/\pi'A$, où
            $\pi$ et $\pi'$ sont les idempotents centraux $\pi = r+r^2; \ \pi' = 1+r+r^2$. On a $\dim B = 4$ et $\dim B' = 8.$            

\smallskip
            
             Le bloc $B$ est le bloc principal; il est donné par la projection $A \to k[C_4]$.  Le groupe $U_G$ se décompose en $U_B \times U_{B'}$; on a donc   $U_G/U^0_G \ \simeq \ U_B/U^0_B \times U_{B'}/U^0_{B'}$.  Nous allons voir que {\it chacun des quotients  $U_B/U^0_B$   et $U_{B'}/U^0_{B'}$  est  d'ordre} 2. Comme $G$ n'a qu'un seul caractère essentiel  $\ne 1$, cela montrera que {\it le groupe
 $E_G$ est d'ordre}~2; de façon plus précise, c'est le facteur $U_{B'}/U^0_{B'}$ de $U_G/U_G^0$.
 
 \smallskip
   Le cas du bloc $B \simeq k[C_4]$ a déjà été traité, cf. 5.5.2, avec $n=4$.
   
   \smallskip
                  
               Passons à $B'$. Soient $r'$ et $s'$ les images de $r$ et de $s$ dans $B'$, et notons $K$ la sous-algèbre de $B'$ engendrée par $r'$. C'est une algèbre quadratique, isomorphe à $k[r']/(r'^2+r'+1) \simeq k \otimes_{\F_2} \F_4 $; comme on a supposé $k$ algébriquement clos,
               on a $K \simeq k \times k$; de plus, si $z=(z_1,z_2)$ est un élément de $K$, on a $z^* = (z_2,z_1)= s'zs'^{-1}$. Tout $x \in B'$ s'écrit de façon unique sous la forme  
              
              \smallskip 
                 $x = a + bs' + cs'^2+ds'^3$, \  avec $a = (a_1,a_2) \in K, ..., d = (d_1,d_2) \in K$.
                  
                  \smallskip
 \n  Lorsqu'on écrit que $x$ est unitaire, on trouve les relations suivantes :   
 
 \smallskip
    
    (5.5.7) \quad $ a_1a_2 + b_1b_2 +  c_1c_2 + d_1d_2=1,$  
    
    (5.5.8) \quad $ (a_1+c_1)(b_1+d_1) = 0$,
    
    (5.5.9) \quad $(a_2+c_2)(b_2+d_2) = 0$,
    
    (5.5.10) \quad  $a_1c_2+c_1a_2 = b_1d_2+b_2d_1$. 
    
    \smallskip     
                  
         \n       Les relations  (5.5.8) et (5.5.9) donnent {\it a priori} quatre possibilités:
  \smallskip              
                
                (5.5.11) \quad  $b_1=d_1$ et $b_2=d_2$ \ (autrement dit $b=d$),  
                  
                     (5.5.12) \quad  $a_1=c_1$ et $a_2=c_2$ \ (autrement dit $a=c$),  

(5.5.13) \quad  $b_1=d_1$ et $a_2=c_2$,

(5.5.14) \quad  $a_1=c_1$ et $b_2=d_2$.

\smallskip

  En fait, (5.5.13) est impossible, car, en ajoutant (5.5.7) et (5.5.10), on trouverait $0=1$. Même   chose pour (5.5.14).
  Il ne reste donc que les seules possibilités (5.5.11) et (5.5.12). Dans le cas de (5.5.11),  les équations
  (5.5.7) et (5.5.10) entraînent:
  
  \smallskip
  
  (5.5.15) \quad $a_1a_2 + c_1c_2 = 1$ \ et \ $a_1c_2=a_2c_1$.
  
  \smallskip
  
  On obtient ainsi un sous-groupe ouvert d'indice 2 de $U_{B'}$, qui
  est produit semi-direct de $\G_m$ (le groupe des $(a,0,0,0)$ avec $a_1a_2=1$) par 
  un groupe isomorphe à  $\G_a \times \G_a \times \G_a$. Ce groupe est connexe; c'est donc la composante neutre de $U^0_{B'}$.
  
  Le cas (5.5.12) donne l'autre composante connexe de $U_{B'}$: celle contenant~$s'$.
  
  \smallskip
  
  \n {\it Remarques}.
  
  1) Lorsque $k$ ne contient pas de racine primitive cubique de l'unité, $K = k\otimes \F_4$ est une extension
  quadratique de $k$, et le groupe $\G_m$ ci-dessus doit être remplacé par son ``$K/k$-tordu'', autrement dit par
  $\Ker(R_{K/k}(\G_m) \to \G_m)$.
  
  2) Le bloc $B'$ est du type décrit au \S3.7: son radical $\r'$ est engendré par l'élément hermitien  $\o = 1+s'^2$, qui est de carré nul, et le quotient  $B'/\r'$ est isomorphe à l'algèbre de matrices $\M_2$. L'image de $U_{B'} \to U_{B'/\r'}\simeq \SL_2$ est le groupe orthogonal ${\bf O}_2$ associé à la forme quadratique $x^2+xy+y^2$: c'est la situation de 3.7.5 (iii).
  
    \medskip
  
  \n 5.5.16. {\it $G = {A}_5 = \PSL_2(\F_5) = \SL_2(\F_4)$,  d'ordre} 60. Ici encore, il y a deux blocs, $B$ et $B'$. 
  
  Le bloc principal $B$  est de dimension 44;  le quotient de $B$ par son radical $\r_B$ est isomorphe
  à $k \times \M_2(k) \times \M_2(k)$, ces trois facteurs correspondant aux représentations irréductibles de $G$ de degré 1, 2 et 2. On a $U_{B/\r_B} \simeq \SL_2 \times \SL_2$. L'homomorphisme $U_B \to U_{B/\r_B}$ est lisse (et donc surjectif): cela se voit en remarquant  que $\Lie U_G \to \Lie(\SL_2 \times \SL_2) $ est surjectif. Un calcul sur ordinateur fait par M. Barakat (cf. [Ba 13]) montre que {\it le noyau de  $U_B \to U_{B/\r_B}$ est connexe}; cela entraîne que  $U_B$ est connexe, et que son plus grand quotient réductif est isomorphe à $\SL_2 \times \SL_2$.
  
    Le bloc $B'$ est isomorphe à $\M_4(k)$ : c'est un bloc de défaut 0, cf. [CR 62, \S86] et [Se 68, \S16.4]; il correspond à la représentation irréductible de degré 4 de  $G$, qui est un $k[G]$-module projectif. On a $U_{B'} \simeq \Sp_4$. Comme $U_G = U_B \times U_{B'}$, on en déduit que $U_G$ est connexe, d'où $E_G = \{1\}$.

  \smallskip
  
  \n {\it Remarque}. Lorsque $k$ ne contient pas de racine primitive cubique de l'unité, on doit remplacer
  $\SL_2 \times \SL_2 $ par $R_{K/k}(\SL_2) $, où  $K= k\otimes \F_4$ désigne comme ci-dessus
  l'extension quadratique de $k$ engendrée par les racines cubiques de l'unité. 
  
  \medskip
  
  \n 5.5.17. {\it $G =  \widehat{{ A}}_5 = \SL_2(\F_5)$, d'ordre 120}. Il y a deux blocs  $B$ et $B'$ qui correspondent à ceux de  ${A}_5$, cf. [Fe 82, Chap.V, \S4]. De façon plus précise, si  $\o = 1+c$, où $c$  est l'unique élément d'ordre 2 de $G$, les quotients $B/\o B$ et $B'/\o B'$ sont les blocs de ${A}_5$ . D'où des homomorphismes: 
  
   \smallskip
\quad  $ \varphi: U_B \ \to \ \SL_2 \times \SL_2 $  \  \  et  \ \ $\varphi' : U_{B'}  \to \Sp_4$.
    
    \smallskip
    
    Ces homomorphismes ont  les propriétés suivantes:
    
          \smallskip

    (5.5.19) {\it L' homomorphisme $\varphi: U_B \ \to \ \SL_2 \times \SL_2$ est surjectif, mais n'est pas lisse.} 
    
          \smallskip

    (5.5.20) {\it L'homomorphisme $\varphi': U_{B'} \ \to \  \Sp_4$ est lisse, mais n'est pas surjectif$;$ son image est
    un groupe orthogonal ${\bf O}_4$}.
   
      \smallskip
      
    \n {\it Démonstration de}  (5.5.19). Si $\varphi$ n'était pas surjectif, son image aurait un facteur ${\bf O}_2$ et un tel groupe n'a pas de sous-groupe isomorphe à ${A}_5$. D'autre part, un calcul facile montre que l'image de $\Lie(U_B) \to \Lie(\SL_2 \times \SL_2)$ est égale au centre de   $ \Lie(\SL_2 \times \SL_2)$; d'où le fait que $\varphi$  n'est pas lisse. D'après [Va 05], cela entraîne:
     
     (5.5.21) {\it L'homomorphisme $U_B \to \SL_2 \times \SL_2$ se factorise par $\SO_3 \times \SO_3$.}
     
     \smallskip      
      
         \n {\it Démonstration de}  (5.5.20). Le bloc $B'$ est du type décrit au \S3.7; de plus, il n'est pas isomorphe à $\M_4 \otimes_k k[t]/(t^2)$; cela se voit,  par exemple, en remarquant que certains hermitiens de  $B'/\o B'$ ne sont pas images d'hermitiens de $B'$. On conclut en utilisant le corollaire 3.7.5. (La lissité de $\varphi'$ se vérifie par un calcul explicite.)

               \smallskip

     Noter deux conséquences de (5.5.19) et (5.5.20):
     
         \smallskip

     (5.5.22) {\it Le groupe $U^0_G/R_u(U^0_G)$  est isomorphe à} $\SO_3 \times \SO_3 \times \ {\bf SO}_4$.
     
         \smallskip

     (5.5.23) {\it On a\ $E_G \neq \{1\}$.} 
     
    \n  On aurait $|E_G| = 2$ si $U_B$ était connexe; j'ignore si c'est le cas.

     \medskip

   \n 5.6. {\bf Nombres de points sur les corps finis.}
    
    \smallskip

    Soit $q= 2^n, n \geqslant 1$, une puissance de $2$, et soit  $n_G(q)$ le nombre d'éléments du groupe
   $U_G(\F_q)$.  On trouve dans la littérature un certain nombre de cas où ce nombre a été calculé, cf. par exemple
   [BR 00]  (cette référence m'a été signalée par M. Barakat). Ce calcul peut être utilisé pour déterminer l'ordre du groupe $U_G/U^0_G$. Voici comment (pour simplifier on se borne au cas où  $G$  est un $2$-groupe, de sorte que  $U_G$ est unipotent):

   On remarque d'abord que, si  $N$  est un groupe unipotent connexe sur  $\F_q$, tout $N$-torseur a un point rationnel, et le nombre de ses points rationnels est $q^{d(N)}$, où  $d(N)$ est la dimension de  $N$; cela se voit par dévissage, à partir du cas $N = \G_a$.
   
   On applique ceci au groupe $N = U^0_G$; on en déduit que, si  $d$  désigne la dimension de $U_G$,
   on a  $U_G(\F_q) = 2^{c(q)}q^d$, où  $2^{c(q)}$ est le nombre de $\F_q$-points du groupe étale $I = U_G/U^0_G$.
   Le groupe des $\Fbar_2$-points de $I$ est un groupe de type (2,...,2), autrement dit un $\F_2$-espace vectoriel; soit  $c$ sa dimension et soit   $\sigma$
   l'automorphisme de Frobenius de ce groupe, relativement au corps $\F_2$. Si $q= 2^n$, l'entier $c(q)$ défini ci-dessus  est la dimension de l'espace des points fixes du $n$-ième itéré $\sigma^n$ de $\sigma$; on a donc $0 \leqslant c(q) \leqslant c$, et l'égalité $c(q) = c$ a lieu lorsque  $n$  est un multiple de l'ordre de  $\sigma$. D'où:
   
   \smallskip
   
   \n {\bf Proposition 5.6.1} - {\it On a $|U_G(\F_q)| = 2^{c(q)}q^d$, où  $d$  est la dimension de
   $U_G$ et  $c(q)$ est un entier compris entre  $0$ et $c$, avec $|U_G/U^0_G| = 2^c;$ de plus, $c(q) = c$   pour une infinité de valeurs de  $q$.}
   
   \smallskip
   
      \n     \n {\bf Corollaire 5.6.2} - {\it On a $U_G= U^0_G$ si et seulement si $ |U_G(\F_q)| = q^d $ pour tout~$q$.}
      
      \smallskip
   
   \n     \n {\bf Corollaire 5.6.3} - {\it On a $|U_G/U^0_G|  \ = \ \sup_q |U_G(\F_q)|/q^d$.}
   
   \medskip
   
  \n {\it Exemple}. Le corollaire 5.6.2, combiné avec les résultats de  [BR 00], montre que $U_G$ est connexe lorsque  $G$  est un 2-groupe diédral ou extra-spécial, et que $|U_G/U^0_G| = 4$ lorsque $G$ est un 2-groupe quaternionien; dans les deux cas, cela entraîne  $E_G = \{1\}$.

\bigskip
                      
 \centerline    {\bf \S6 - Structure des $G$-formes trace}
     
     \smallskip
     
       Comme au \S5,  $G$ est un groupe fini, et $k$ est un corps de caractéristique  2 (non nécessairement parfait). On note $\kbar$ une clôture algébrique de $k$ et $k_s$ la plus grande extension séparable de  $k$  contenue dans $\kbar$; le groupe
       de Galois $\Gal(k_s/k) = \Aut_k(\kbar)$ est noté $\Gamma_k$.
       On rappelle que $G_0$ désigne le sous-groupe de $G$ engendré par les éléments d'ordre 2 et par les carrés, cf. \S5.3.
       
         \smallskip

       \n 6.1. {\bf Le théorème principal  et ses corollaires.}
       
       \smallskip
       
       \n 6.1.1. {\it La $G$-forme trace associée à une $G$-algèbre galoisienne.}
       
         \smallskip

        Soit $L$ une $G$-algèbre galoisienne sur  $k$. Rappelons  que cela signifie que $L$ est une $k$-algèbre commutative, qui est étale (i.e. produit fini d'extensions finies séparables de $k$),
        et qui est munie d'une action de $G$ qui en fait un $k[G]$-module libre de rang 1. D'autres caractérisations se trouvent dans [A VIII, \S16.7] et [BFS 94, \S1.3]; l'une d'elles est: ``$L$ est l'algèbre affine d'un $G$-torseur sur $k$''. Lorsque $L$ est un corps, il revient au même de dire que $L$ est une extension galoisienne de $k$, munie d'un isomorphisme $G \simeq \Gal(L/k)$; cela explique la terminologie choisie.
        
          \smallskip

    Nous noterons $q_L: L \times L \to k$ la forme bilinéaire symétrique définie par:
    
      \smallskip

     \  (6.1.2)  $q_L(x,y) = \Tr_{L/k}(xy)$,  \ où  $\Tr_{L/k}$ désigne la trace.
     
       \smallskip

     Comme $L/k$ est étale, cette forme bilinéaire est non dégénérée. Elle est invariante par $G$: on a
     
       \smallskip

     \ (6.1.3)  $q_L(gx,gy) = q_L(x,y)$ \  quels que soient $g \in G, x\in L, y\in L$.
     
     \smallskip
     
      On dit que $q_L$ est la {\it $G$-forme trace} de $L$. Deux telles $G$-formes  $q_L$ et $q_{L'}$
      sont isomorphes si et seulement si il existe une bijection $k$-linéaire $f:L \to L'$, compatible à l'action de $G$, et telle que  $q_{L'}(f(x),f(y)) = q_L(x,y)$ quels que soient $x,y \in L$. On écrit alors $q_L \simeq_G q_{L'}$.

        \smallskip

      \n 6.1.4. {\it Le théorème principal.} 
      
        Si $L$ est une $G$-algèbre galoisienne, et si $H$ est un sous-groupe normal de  $G$, l'ensemble
        $L^H$ des éléments de $H$ fixés par $H$ est une sous-algèbre de $L$ sur laquelle $G/H$ opère, et c'est une $G/H$-algèbre galoisienne. Nous allons utiliser ceci pour $H=G_0$:
        
          \smallskip

\n {\bf Théorème 6.1.5} - {\it Soient $L$ et $L'$ deux $G$-algèbres galoisiennes sur $k$. Les deux propriétés suivantes sont équivalentes}:

 (i) {\it Les $G$-formes  $q_L$ et $q_{L'}$ sont isomorphes.}
 
 (ii) {\it Les $G/G_0$-algèbres galoisiennes $L^{G_0}$ et $L'^{G_0}$ sont isomorphes.}

\smallskip

La démonstration sera donnée au \S6.3.

\smallskip

\n 6.1.6. {\it Traduction du théorème 6.1.5 en termes d'homomorphismes.}

\smallskip
        
  Une $G$-algèbre galoisienne est déterminée à isomorphisme près  par un homomorphisme continu
  $\varphi_L : \Gamma_k \to G$; deux homomorphismes donnent des $G$-algèbres  isomorphes si et seulement si ils sont $G$-conjugués, cf. [BFS 94, \S1.3.1].  Le théorème 6.1.5 peut donc être reformulé de la façon suivante:
  
  \smallskip
  
   \n {\bf Théorème 6.1.7} - {\it  Soient $\varphi_L, \varphi_{L'} : \Gamma_k \to G$  les homomorphismes associés à $L$ et $L'$, et soient $\varphi^0_L, \varphi^0_{L'} : \Gamma_k \to G/G_0$ leurs composés avec $G \to G/G_0$. Pour que
   les $G$-formes traces $q_L$ et $q_{L'}$ soient isomorphes, il faut et il suffit que $\varphi^0_L = \varphi^0_{L'}.$} 
    
    \smallskip
    
 \n {\it Remarque}. Bien que $\varphi_L$  ne soit défini qu'à $G$-conjugaison près, son composé $\varphi^0_L$ avec $G \to G/G_0$ est défini sans ambiguïté, puisque $G/G_0$ est commutatif. Comme de plus
 $G/G_0$ est un 2-groupe élémentaire (cf. \S5.3), on peut utiliser la théorie d'Artin-Schreier ([A V, \S11.9]) pour interpréter $\varphi^0_L$ en termes du quotient $k/\wp(k)$, où $\wp: k \to k$ est $x \mapsto x^2+x$; de façon plus précise, si $X_G$ désigne le groupe des caractères essentiels de $G$, on peut identifier $\varphi_L$ à un homomorphisme $\alpha_L : X_G \to k/\wp(k)$ et le théorème 6.1.5 dit que $\alpha_L$ détermine $q_L$.  
 
 \smallskip
 \n 6.1.8. {\it Application à l'existence d'une BNA et démonstration du théorème B.}
 
 \smallskip
  Prenons pour $L'$ une $G$-algèbre galoisienne scindée (``split''), autrement dit un produit $k \times ... \times k$ de copies 
  de $k$, permutées de façon simplement transitive par $G$; du point de vue galoisien, cela revient à 
  dire que $\varphi_{L'} : \Gamma_k \to G$ est trivial. La $G$-forme trace correspondante est la {\it $G$-forme unité}; elle a une base normale autoduale canonique: l'ensemble des idempotents indécomposables $(1,0,..,0),...,(0,...,0,1)$. Le théorème 6.1.7 entraîne donc:
  
  \smallskip
  
    \n {\bf Théorème 6.1.9} - {\it Pour que $L$ possède une base normale autoduale, il faut et il suffit 
    que $\varphi^0_L = 1$, autrement dit que l'image de  $\varphi_L : \Gamma_k \to G$ soit contenue dans $G_0$,
    ou encore que $L$ soit isomorphe à une algèbre induite $\Ind^G_{G_0} M$, où $M$ est une $G_0$-algèbre 
    galoisienne.}
    
    [Pour la notion d' {\it induction}  des algèbres galoisiennes (appelée {\it coïnduction} dans [A VIII, \S16.7]), voir par exemple [BFS 94, \S1.3.2].]

\smallskip

  \n {\bf Corollaire 6.1.10} - {\it  Supposons que $L$ soit un corps. Pour que $L$ ait une BNA, il faut et il suffit que $G = G_0$, autrement dit} (cf. proposition 5.3.2), {\it que $G$ soit engendré par des éléments d'ordre $2$ et par des éléments d'ordre impair.}
  
 [C'est le théorème B de l'introduction.]
 
 \smallskip
 
 Cela résulte du théorème 6.1.9: en effet, l'hypothèse que $L$ est un corps équivaut à dire que $\varphi_L$ est surjectif; son image ne peut être contenue dans $G_0$ que si $G=G_0$.
 
 \smallskip

\smallskip

\n 6.1.11. {\it Application aux extensions de degré impair.}

\smallskip

\n   \n {\bf Corollaire 6.1.12} - {\it Soient $L$ et $L'$ deux $G$-algèbres galoisiennes, et soit $k_1$ une extension finie de $k$ de degré impair. Si les $G$-formes $q_L$ et $q_{L'}$ deviennent isomorphes après extension des scalaires à $k_1$, elles sont $k$-isomorphes.}

\smallskip

\n {\it Démonstration.} Le groupe $\Gamma_{k_1}$ est un sous-groupe d'indice impair de $
\Gamma_k$; cela entraîne que, si $\psi, \psi'$ sont deux homomorphismes de $\Gamma_k$ dans un $2$-groupe abélien qui coïncident sur $\Gamma_{k_1}$, alors $\psi = \psi'$. D'où    $\varphi^0_L = \varphi^0_{L'} $, et l'on applique le théorème 6.1.7.

\smallskip
\n {\it Variante}. Utiliser le corollaire 4.9.3.

 \smallskip

\n {\it Remarque.} J'ignore si cet énoncé s'étend à des $G$-formes bilinéaires symétriques quelconques (comme c'est le cas en caractéristique $\ne 2$, cf. [BFL 90, th.4.1]).

\medskip
 \n 6.1.13. {\it Un théorème d'existence.}

 \smallskip

 \n {\bf Théorème 6.1.14} - {\it  Soit  $H$ un sous-groupe de $G$ tel que $H.G_0=G$ et soit $L_0$ une $G/G_0$-algèbre galoisienne. Il existe une $H$-algèbre galoisienne $M$ telle que, si l'on pose $L' = \Ind^G_H M$,
 on ait $L'^{G_0} \simeq L_0$.}
 
 \smallskip
 
 \n {\it Démonstration}. Soit $S$ un 2-Sylow de $H$; puisque  $H \to G/G_0$ est surjectif, et que $G/G_0$ est un 2-groupe, l'homomorphisme $S \to G/G_0$ est surjectif.  On a d'autre part $cd_2(\Gamma_k) \leqslant 1$, cf.  [Se 64, II, \S2.2]; d'après [Se 64, I, \S3.4], cela entraîne
 que l'homomorphisme $ \varphi_{L_0}: \Gamma_k \to G/G_0$ se relève à $S$, donc {\it a fortiori} à $H$. On obtient
 ainsi une $H$-algèbre galoisienne qui a la propriété requise.
 
 \smallskip
 
 \n {\bf Corollaire 6.1.15} - {\it  Pour toute $G/G_0$-algèbre galoisienne $L_0$, il existe une $G$-algèbre galoisienne $L$ telle que $L^{G_0} \simeq L_0$.}
 
 \smallskip
 
   C'est le cas particulier $H = G$. 
   
   \smallskip
 
 \n {\bf Corollaire 6.1.16} - {\it Soit $L$ une $G$-algèbre galoisienne et soit $H$ un sous-groupe de $G$
 tel que} $H.G_0 = G$ (par exemple un 2-Sylow de $G$). {\it Il existe une $H$-algèbre galoisienne $M$
 telle que les $G$-formes trace de $L$ et de $\Ind^G_H M$ soient isomorphes.}
 
\smallskip

{\it Démonstration.} On applique le théorème 6.1.14 à $L_0 = L^{G_0}$; on obtient une $H$-algèbre $M$, d'où une algèbre induite $L'=\Ind^G_H M$, avec $L'^{G_0} \simeq L^{G_0}$. D'après le théorème 6.1.5, cela entraîne  que les
$G$-formes trace de $L$ et de $L'$ sont isomorphes.

\smallskip

\n 6.1.17. {\it Application au principe de Hasse pour les $G$-formes trace.}

\smallskip
  Supposons que $k$ soit un {\it corps global} de caractéristique 2, autrement dit une extension de type fini de $\F_2$ dont le degré de transcendance est 1. Le ``principe de Hasse'' est l'énoncé suivant (démontré en caractéristique $\ne 2$ dans [BPS 13]):
  
  \smallskip
  
   \n {\bf Théorème 6.1.18} - {\it Soient $L$ et $L'$ deux $G$-algèbres galoisiennes sur $k$. Si,
   pour toute place $v$ de $k$, les $G$-formes trace  $q_L$ et $q_{L'}$ deviennent isomorphes sur le complété
   $k_v$ de  $k$ en $v$, alors elles sont isomorphes.}
   
   \smallskip
   
   \n {\it Démonstration}. Soient $\varphi^0_L, \varphi^0_{L'} : \Gamma_k \to G/G_0$
   les homomorphismes définis dans le théorème 6.1.7. Vu le théorème en question, il nous faut prouver que $\varphi^0_L = \varphi^0_{L'}$. Or, par hypothèse, cette égalité devient vraie lorsqu'on restreint ces homomorphismes aux groupes de décomposition des places de $k$; comme ces groupes engendrent topologiquement $\Gamma_k$
   (cela résulte par exemple du théorème de densité de Chebotarev), on obtient le résultat cherché.
   
   \smallskip
   
   \n {\it Remarque.} On peut remplacer {\it toute place} par {\it  toute place sauf un nombre fini}, ou même par {\it toutes les places d'un ensemble de densité analytique $> 1/2$ }; la démonstration est la même. Une telle amélioration n'est pas possible en caractéristique $\ne 2. $
                  
    \medskip
    
    \n 6.2. {\bf  Elimination des extensions radicielles}.            
        
        \smallskip 
                  
     \n {\bf Proposition 6.2.1} - {\it Soient $L$ et $L'$ deux $G$-algèbres galoisiennes sur $k$. Supposons que
les $G$-formes  $q_L$ et $q_{L'}$ deviennent isomorphes sur une extension radicielle de $k$. Alors $q_L$ et $q_{L'}$ sont isomorphes.}
   
     % \n   [En abrégé: les extensions radicielles n'ont pas d'importance.]     

 \smallskip 
\n {\it Démonstration}. On peut supposer que l'extension radicielle en question est de degré fini sur $k$, donc contenue dans $k^{1/q}$, où  $q$ est une puissance convenable de  2.  Posons
$ L^{1/q} = L \otimes_k k^{1/q}$, et définissons de même   $L'^{1/q}$. Par hypothèse, il existe un isomorphisme de $G$-formes $ \theta: L^{1/q}  \to L'^{1/q}$. 

 \smallskip 
 
 \n {\bf Lemme 6.2.2} - {\it Il existe un isomorphisme
de $G$-formes $\theta_0 : L \to L'$ et un seul tel que l'on ait}

\smallskip

(6.2.3) \ $ \theta_0(x^q) =   \theta(x)^q $  \ {\it pour tout}  $x \in  L^{1/q}$.               
                  
           \smallskip 
          \n {\it Démonstration du lemme 6.2.2.} Remarquons d'abord que l'application $x \mapsto x^q$ est une bijection de $L^{1/q}$ sur $L$ (ce qui justifie la notation); cela se vérifie en se ramenant au cas d'un corps, qui est bien connu (cf. e.g. [A V,  \S6.7 et \S15.4]). Ce fait, appliqué à $L$ et à $L'$, montre qu'il existe une
          bijection unique $\theta_0 : L \to L'$ ayant la propriété (6.2.3). Il est clair que $\theta_0$ est $k$-linéaire, et commute à l'action de $G$; le fait qu'elle transforme $q_L$ en $q_L'$ résulte de la formule
          $\Tr(x^q) = \Tr(x)^q$, qui elle-même se démontre en se ramenant au cas où $L$ est un corps
          (ou bien en écrivant  $\Tr(x)$ comme $\sum_{g\in G} gx$.)
          
          \smallskip
          
          Il est clair que le lemme 6.2.2 entraîne la proposition 6.2.1.
          
          \smallskip
          
          \n {\it Remarque.} La proposition 6.2.1 est spéciale aux $G$-formes trace, autrement dit aux $G$-formes associées à des $G$-algèbres galoisiennes; elle ne s'étend pas à des $G$-formes quelconques, comme le montre déjà le cas $G = \{1\}$, où de telles formes correspondent aux
          éléments de $k^\times/k^{\times2}$. Nous reviendrons là-dessus au \S7.4.

       \medskip
       
       \n 6.3. {\bf Démonstration du théorème 6.1.5.}
       
       \smallskip
                   
                  Soient $L$ et $L'$ deux $G$-algèbres galoisiennes sur $k$. Nous devons montrer l'équivalence de :
                  
                  \medskip
                  
                  (6.3.1) {\it Les $G$-formes  $q_L$ et $q_{L'}$ sont isomorphes.}
                  
\smallskip

\n et

\medskip
 
 (6.3.2) {\it Les $G/G_0$-algèbres galoisiennes $L^{G_0}$ et $L'^{G_0}$ sont isomorphes.}

                  \smallskip
                  
                  \n 6.3.3. {\it Réduction au cas où le corps $k$ est parfait.} Soit $k_i = k^{2^{-\infty}}$ la clôture parfaite du corps $k$. Si la condition  (6.3.1) (resp. la condition  (6.3.2)) est satisfaite sur $k$, elle l'est sur $k_i$. Inversement, si elle est satisfaite sur $k_i$, elle l'est sur $k$: c'est clair pour (6.3.2) puisque $\Gamma_{k_i} = \Gamma_k$, et pour (6.3.1) cela résulte de la proposition 6.2.1.
                  
                    On déduit de là que l'on peut supposer que $k = k_i$ , autrement dit que $k$ est parfait.
                    Nous ferons cette hypothèse dans la suite de la démonstration.
                    
                    \smallskip
                    
                    \n 6.3.4. {\it Traduction en termes de cohomologie galoisienne.} 
                    
                                        \smallskip
                    
                    Comme on l'a rappelé au \S6.1.6, la $G$-algèbre $L$ est déterminée à isomorphisme près
                    par un homomorphisme continu       $\varphi_L : \Gamma_k \to G$, défini à conjugaison près, autrement dit par un élément de  $H^1(k,G)$. Soit $f_L$ cet élément. L'injection $G \to U_G$ transforme 
                    $f_L$ en un élément $u_L$ de $H^1(G,U_G)$. On définit de même $f_{L'} \in H^1(k,G)$ et
                     $u_{L'} \in H^1(k,U_G)$.  On sait - voir par exemple [BFS 94, proposition 1.5.1]\footnote{Dans [BFS 94],
                     la caractéristique de $k$ est supposée $\ne 2$; toutefois cette hypothèse ne joue aucun rôle
                     dans la démonstration, à condition d'interpréter $q_L$ comme une $G$-forme bilinéaire symétrique, et non comme une $G$-forme quadratique.} - que la propriété (6.3.1)  est équivalente à:
                     
                     \smallskip
                     (6.3.5) $u_L = u_{L'}$ \ {\it dans}  $H^1(k,U_G)$.
                               
                  \smallskip
                  
                       \n 6.3.6. {\it Fin de la démonstration.}

                 \smallskip 
                 \n  Soient $u_L^0$ et $u_{L'}^0$ les images de $u_L$ et de $u_{L'}$ dans $H^1(k,U_G/U_G^0)$.
                  D'après le corollaire 4.9.2 (qui est applicable puisque $k$ est parfait), (6.3.5) équivaut à:

                  \smallskip
                  (6.3.7)  $u_L^0 =u_{L'}^0$ \ {\it dans} $H^1(k,U_G/U_G^0)$.
                  
                  \smallskip
                  
                \n  D'après (5.4.2), on a $U_G/U_G^0 = G/G_0 \times E_G$, \ d'où:
                  
                  \smallskip
                  
                   (6.3.8) \ $H^1(k,U_G/U_G^0) = H^1(k,G/G_0) \times H^1(k,E_G)$.
                   
                   \smallskip
             \n     L'élément $u_L^0$ de $H^1(k,U_G/U_G^0)$ a deux composantes: 
   l'une dans $H^1(k,G/G_0)$ et l'autre dans $H^1(k,E_G)$. 
    La seconde est triviale car $G \to U_G \to~E_G$ applique $G$  en l'élément neutre de $E_G$. La première est à valeurs dans $H^1(k,G/G_0) = \Hom_{\rm cont}(
                   \Gamma_k,G/G_0)$; c'est l'homomorphisme noté $\varphi^0_L$ dans le théorème 6.1.7,
                  autrement dit le composé  $\Gamma_k \stackrel{\varphi_L}{\rightarrow} G  \to G/G_0$. On déduit de là que (6.3.7) équivaut à $\varphi^0_L = \varphi^0_{L'}$,
                 ce qui est équivalent à (6.3.2), et termine la démonstration.

                   \bigskip
                      
 \centerline    {\bf \S7 - Compléments}
     
     \smallskip

  \smallskip

  \n 7.1. {\bf Notations.}
  
  \smallskip
  Ce sont les mêmes qu'aux \S\S5,6: $G$ est un groupe fini, $k$ est un corps de caractéristique 2, et 
    $A = k[G]$. On note $\t : A \to k$  l'homomorphisme d'augmentation, i.e. l'unique forme linéaire 
  sur $A$ telle que  $\t(g) = 1$ pour tout  $g \in G$. On note $\sigma_G$ l'élément $\sum_{g\in G} g$ de $A$.
  
  \smallskip
  
   \n Si $x = \sum_{g \in G} x_g g$  est un élément de $A$, on pose:
    
    \smallskip
    
    (7.1.1)  \ $\delta_g(x) = x_g$;  
    
    \smallskip
    
  \n  lorsque $g=1$, on écrit  $\delta$ à la place de  $\delta_1$.
    
    \smallskip
    
  \n    On a :
      \smallskip
      
      (7.1.2) \ $\delta(xy) = \delta(yx)$  \ et  \ $\delta(x^*) = \delta(x)$ \ si $x,y \in A$.
      
      \smallskip
      
      \n et
      
      \smallskip
      
      (7.1.3) \ $\t = \sum_{g\in G} \delta_g$.
      
      \medskip
  
  On s'intéressera à des {\it $G$-formes} $(V,q)$ au sens suivant:
  
  \smallskip
(7.1.4) {\it $V$ est un $k$-espace vectoriel muni $:$

\ \ $\bullet$ d'une action de $G$ qui en fait un $A$-module libre de rang $1,$

\ \ $\bullet$ d'une forme bilinéaire symétrique $q$ non dégénérée et invariante par} $G$.
  
  \medskip
  
   On a  $\dim V = |G|$.
   
   \medskip
   
   \n{\it Remarque}. La condition que $q$ soit invariante par $G$ est équivalente à :
   
   \smallskip
   
   (7.1.5)  \ $q(ax,y) = q(x,a^*y)$ \ quels que soient  \ $x,y \in V$ et $a\in A$.
   
   \medskip
  
 \n {\it Exemple: la forme unité}. C'est le cas $(V,q) = (A,q_1)$, où $q_1(x,y) = \delta(xy^*)$ avec les notations ci-dessus. Si $g$ et $g'$ sont deux éléments de $G$, on a  $q_1(g,g') = 1$ si $g=g'$ et $q_1(g,g')=0$ sinon; autrement dit, les éléments de $G$ forment une BNA de $(A,q_1)$.
   
  \medskip

  \n 7.2. {\bf Deux questions.}
  
  \smallskip
       
           Nous avons vu au \S6 que, si $L$ est une $G$-algèbre galoisienne, la $G$-forme  $q_L$
ne dépend que de la $G/G_0$-algèbre galoisienne  $L^{G_0}$, ou, ce qui revient au même, du
composé $\Gamma_k \to G \to G/G_0$ (ou encore de l'homomorphisme $\alpha_L : X_G \to k/\wp(k)$ du \S6.1.7).
           
        Il est naturel de se poser la question suivante:

      \medskip
      
      (7.2.1) {\it Lorsqu'on connait $ \Gamma_k  \to G/G_0, $ peut-on décrire explicitement la $G$-forme
   trace correspondante ?} 
      
       \medskip

       Autre question, tout aussi naturelle:
       
        \medskip

  (7.2.2) {\it Comment caractériser les $G$-formes} (au sens de (7.1.4)) {\it qui  sont des formes trace ?}
  
     \medskip

La suite de ce \S \  donnera des réponses partielles à ces questions; voir notamment les théorèmes 7.3.2, 7.10.4 et 7.10.10.

\bigskip

\n 7.3 - {\bf Exemple de réponse à la question $(7.2.1)$.}

\smallskip
  Lorsque l'image de $\Gamma_k$   dans $ G/G_0$ est triviale, la $G$-forme trace est la forme unité, nous l'avons vu. Supposons maintenant que cette image soit non triviale, mais aussi petite que possible. Cela revient à faire l'hypothèse suivante:

\medskip
 
 (7.3.1) {\it  L'image de $\Gamma_k  \to G/G_0$ est un sous-groupe d'ordre 2 de $G/G_0$}.

  \smallskip
   Notons $\{1,\gamma\}$ l'image de $\Gamma_k $ dans $G/G_0$. Le noyau de $\Gamma_k  \to \{1,\gamma\}$ est alors un sous-groupe ouvert d'indice 2 de $\Gamma_k$. Il correspond à une extension quadratique de  $k$, que l'on écrit à la Artin-Schreier comme  $k(t)$  avec  $t^2 + t = z$, \  $z \in k$. Choisissons un élément  $s \in G$, d'ordre
  une puissance de 2, dont l'image dans $G/G_0$ soit $\gamma$ (un tel élément existe car, si $S$  est un 2-Sylow de  $G$, l'homomorphisme  $S \to G \to G/G_0$ est surjectif). 
  
  \medskip
  
  \n {\bf Théorème 7.3.2} -  {\it  Avec les hypothèses et notations ci-dessus, il existe une base $\{v\}$ du $A$-module $L$  telle que $:$
  
  \smallskip
   $ (7.3.3)$ \ \ $q_L(v,v) = 1, \ q_L(v,sv) = z $ \  et \ $q_L(v,gv) = 0$ \  pour tout} $g \ne 1,s,s^{-1}$.
    
    \smallskip
       \n [La forme $q_L$ est donc presque la forme unité, la différence étant dictée
    par le caractère quadratique de $\Gamma_k$ fourni par $\Gamma_k  \to G/G_0$.]

    \medskip
    
  \n  La démonstration sera donnée au \S7.10.6.
    
    \medskip
 \n 7.4 - {\bf Les propriétés particulières des $G$-formes trace.}

 \medskip
 
  Soit $(V,q)$ une $G$-forme. On va donner trois  {\it conditions nécessaires} pour que  $(V,q)$  soit isomorphe à une $G$-forme trace\footnote{En fait, comme on le verra, la troisième condition entraîne la seconde.}.
 Dans quelques cas (pas très nombreux), nous montrerons plus loin que ces conditions sont suffisantes.

 \medskip
 
\n  $\bullet$ {\it Première condition}:

\smallskip

(7.4.1) \ {\it  On a  $q(x,sx) = 0$ \ pour  tout  $s \in G$ d'ordre 2 et tout  $x \in V$.}
  
   \medskip
    C'est vrai si $(V,q) = (L,q_L)$, car le produit $y = x.sx$ appartient à la sous-algèbre $L^s$ de  $L$ 
    fixée par  $s$, et l'on a  $\Tr_{L/k}(y)= \sum_{g\in G}  gy = 0$ car les termes relatifs à  $g$  et à  $gs$  se détruisent deux à deux.
    
    (Variante: la trace  $\Tr_{L/k}$ se factorise en  $\Tr_{L^s/k} \circ  \Tr_{L/L^s}$ et l'on a
    $\Tr_{L/L^s}(y) = 2y =0$.)

    C'est cette condition qui est la plus importante (et qui est spéciale à la caractéristique  2); noter que c'est une condition ``géométrique" : elle est invariante par extension du corps de base.
    
 \medskip

\n  $\bullet$ {\it Deuxième condition}:

\smallskip

 (7.4.2) \ {\it Il existe $e\in V$ tel que  \  $\sigma_Gx = q(x,e)e$ \ pour tout $x\in V$.}
 
 \n [Rappelons que $\sigma_G = \sum_{g\in G} g$.]

\medskip

\n {\it Exemple.} Lorsque $(V,q)$ est la forme unité $(A,q_1)$, on a $e= \s_G$.

 \medskip

 Dans le cas $(V,q) = (L,q_L)$, on prend  $e=1$. Noter que  $e$, s'il existe,  est non nul, car $\sigma_Gx$ n'est pas toujours  0  puisque  $V$  est $k[G]$-libre de rang 1; de plus, $e$ est unique, et fixé par $G$; nous l'appellerons  {\it l'élément canonique} de $(V,q)$.

 \medskip
 
 \n {\it Remarques.}  
 
 \n  (7.4.2.1) Puisque $V$ est $k[G]$-libre de rang 1, le sous-espace $V^G$ de $V$ fixé par $G$ est de dimension 1; si $\{v\}$ est une base de $V^G$, il existe une unique forme linéaire  $\l(x)$ sur  $V$ telle que  $\sigma_G x = \l(x)v$ pour tout $x\in V$. Cette forme linéaire est $G$-invariante; elle s'écrit donc $\l(x) = \lambda q(v,x)$ avec $ \lambda \in  k^\times$; la condition (7.4.2) équivaut à dire que {\it $\lambda$ est un carré}, auquel cas on a $e=\lambda^{1/2}v$. Noter que cette condition est satisfaite après extension du corps de base à $k^{1/2}$.

\medskip

 \n (7.4.2.2)  Voici une autre caractérisation de l'élément  $e$ (en supposant (7.4.1)):
  
  \medskip
  
  ($7.4.2'$) \ {\it On a  $ q(e,x)^2 = q(x,x)$ \ pour tout  $x \in V$.}
  
  \medskip
  
  \n  En effet, on déduit de (7.4.2):
    
    \medskip
    
    $q(e,x)^2 = q(q(e,x)e,x) = q(\sigma_Gx,x)  =  \sum q(gx,x)$.
    
    \medskip

Dans  $\sum q(gx,x)$, il y a trois types de termes: celui avec $g=1$, qui donne  $q(x,x)$;
 ceux où $g$ est d'ordre 2, qui donnent 0 d'après (7.4.1); ceux avec  $g$  d'ordre $>2$,
qui se détruisent deux à deux, car  $q(gx,x) = q(g^{-1}x,x)$.
On a donc $\sum q(gx,x) = q(x,x)$, ce qui démontre ($7.4.2'$).

\medskip

\n {\it Une conséquence de} (7.4.2):

\smallskip

\n {\bf Proposition 7.4.3} - {\it Soit $(V,q)$ une $G$-forme satisfaisant à} (7.4.2). {\it Supposons qu'il existe une
extension de $k$ sur laquelle $(V,q)$ devienne isomorphe à la forme unité. Il existe alors une telle extension qui est séparable et finie sur $k$.}

\n [Autrement dit, $(V,q)$ peut s'obtenir par {\it descente galoisienne} à partir de la forme unité.]

\smallskip
\n {\it Démonstration}. Soit $e$ l'élément canonique de $(V,q)$. Soit $P$ le schéma des isomorphismes du triplet  $(A,q_1,\sigma_G)$  sur le triplet $(V,q,e)$. C'est un torseur sous le groupe des 
automorphismes de $(A,q_1,\sigma_G)$. Or le groupe des automorphismes de $(A,q_1)$ est $U_G^{\sch}$, et celui de $(A,q_1,\sigma_G)$ est $U_G$ : cela se voit en remarquant que, si $u$ est un point de $U_G^{\sch}$,
on a $u\sigma_G = \t(u)\sigma_G$, donc  $u$ fixe $\sigma_G$ si et seulement si $\t(u)=1$. Ainsi,
$P$ est un $U_G$-torseur, donc est lisse, donc possède des points sur une extension finie séparable du corps de base (cf. par exemple [BLR 90, \S2.2, cor.13]).

\n[Voici une description concrète de l'ensemble $P(k')$ des points de $P$ à valeurs dans une $k$-algèbre commutative $k'$: un point de $P(k')$ est un point $v$ de $k' \otimes_k V$ tel que $q(gv,g'v) = 1$ si
$g=g'$, $q(gv,g'v) = 0$ si $g\neq g'$, et $\sigma_Gv = 1 \otimes e$.]

\bigskip

\n  $\bullet$ {\it Troisième condition}:

 \smallskip
 
 (7.4.4) \ {\it Il existe une application additive $F : V \to V$ jouissant des propriétés suivantes}:
 
\n  (7.4.4.1) (semi-linéarité) \ $F(\lambda x) = \lambda^2F(x)$ \ {\it pour tout $\lambda \in k$ et tout $x\in V$.}
  
\n  (7.4.4.2) (compatibilité avec $G$)  \ $F(gx) = gF(x)$ \ {\it pour tout $g \in G$ et tout $x\in V$.}
  
\n  (7.4.4.3) (compatibilité avec $q$) \ $q(Fx,Fy) = q(x,y)^2$ \ {\it pour tous} $x,y \in V$.

 \medskip
 
Lorsque  $(V,q) = (L,q_L)$, on prend pour $F$ l'application $x \mapsto x^2$; les propriétés (7.4.4.1) et (7.4.4.2)
sont immédiates, et (7.4.4.3) résulte de ce que  Tr$_{L/k}(x^2) = $ Tr$_{L/k}(x)^2$ pour tout $x \in L$.

\medskip
 
\n Noter que $F$ est injective (et même bijective si  $k$  est parfait). 

\smallskip
Il est commode de reformuler (7.4.4) comme:

\medskip
 
$(7.4.4')$ \ {\it Soit $\pi : k \to k$ l'application $\lambda \mapsto \lambda^2$ et soit  $(V,q)_\pi$ la $G$-forme déduite de $(V,q)$ par le changement de base $\pi$. On a $(V,q) \simeq_G (V,q)_\pi$.}

\medskip
 
   Une autre façon de formuler $(7.4.4')$ consiste à se placer sur le corps  $k^{1/2}$; on dispose alors de deux
   $G$-formes: celle obtenue à partir de $(V,q)$ par le changement de base correspondant à 
   l'inclusion $\iota:  k \to k^{1/2}$, et celle obtenue par l'isomorphisme $\s : k \to k^{1/2}$ donné par $\lambda \mapsto \lambda^{1/2}$. L'énoncé devient alors:
   
   \medskip
   
   $(7.4.4'')$ \ {\it Les deux $G$-formes ainsi définies sur $k^{1/2}$ sont isomorphes.}
   
   \medskip

   \n {\bf Proposition 7.4.5} - {\it On a $(7.4.4) \Rightarrow (7.4.2).$}
   
   \smallskip
   \n {\it Démonstration}. Supposons que $(V,q)$ satisfasse à (7.4.4). D'après $(7.4.4'')$ les deux $G$-formes
   sur $k^{1/2}$ obtenus par les changements de base  $\iota:  k \to k^{1/2}, \ \lambda \mapsto \lambda$ et
   $\sigma : k \to k^{1/2}, \  \lambda \mapsto \lambda^{1/2}$ sont isomorphes. D'après (7.4.2.1), la première
   satisfait à (7.4.2). Il en est donc de même de la seconde; comme $\sigma : k \to k^{1/2}$ est un isomorphisme ,
    il en va de même pour $(V,q)$, par transport de structure.

   \medskip
   \n {\bf Proposition 7.4.6} - {\it Soient $(V,q)$ et $(V',q')$ deux $G$-formes satisfaisant à $(7.4.4)$. Si ces 
   $G$-formes deviennent isomorphes sur une extension radicielle de $k$, elles sont $k$-isomorphes.}
   
    [Dans le cas des formes trace, on retrouve la proposition 6.2.1.]
   
   \medskip
   
   \n {\it Démonstration}. En raisonnant par récurrence, on se ramène au cas où l'extension radicielle considérée est  $k^{1/2}$. On applique alors $(7.4.4'')$ comme dans la démontration précédente:  les deux $G$-formes deviennent isomorphes sur $k^{1/2}$ lorsqu'on fait le changement de base $\s$, et comme $\s$ est un isomorphisme, le même énoncé vaut sur $k$. [On peut aussi procéder de façon plus explicite, comme dans la démonstration de la proposition 6.2.1.]
      \bigskip

\n   7.5 - {\bf Exemples et contre-exemples relatifs à la question (7.2.2).}
   
    \smallskip
    
   Commençons par une question de nature ``géométrique'':   
   
   \smallskip
    (7.5.1) - {\it Si $k$ est algébriquement clos, est-il vrai que toute $G$-forme satisfaisant à $(7.4.1)$  est isomorphe à la forme unité ?}   
    \smallskip
    
    Je ne connais pas de contre-exemple.
    
    \medskip
    
    \n {\it Remarque.} Si la réponse est ``oui'', alors toute $G$-forme (sur un corps $k$  quelconque) satisfaisant à (7.4.1) et (7.4.2) s'obtient par torsion galoisienne à partir de la $G$-forme unité, donc correspond à un élément de $H^1(k,U_G)$; cela résulte de la proposition 7.4.3.
        
     \medskip
   
   \n {\bf Théorème 7.5.2} - {\it La réponse à $(7.5.1)$ est `` oui ''  lorsque le groupe $G$ a l'une des trois propriétés suivantes} $:$
   
  (i) {\it  \ il est commutatif} ;
   
   (ii){\it \  son ordre est une puissance de} 2;
   
   (iii) {\it  son ordre est impair}.  
   
    \medskip 
    
 \n  Pour la démonstration, voir \S 7.7.
 
 \medskip
 
   Lorsqu'on ne suppose pas que $k$ est algébriquement clos, il n'est pas difficile de donner des exemples de $G$-formes
     satisfaisant à (7.4.1) et (7.4.2) qui ne sont pas des formes trace, même si $k$ est parfait. Par exemple:
                    
      \medskip
      
      \n {\bf Proposition 7.5.3} - {\it Si $E_G \neq \{1\}$ $($i.e. si $U_G \neq G.U^0_G$, cf. \S5.4$)$,il existe une extension finie $k$ de $\F_2$ et une $G$-forme $(V,q)$ sur $k$ telle que} :
      
       (a) {\it $(V,q)$ satisfait à (7.4.1) et (7.4.2)};
       
       (b) {\it $(V,q)$ n'est pas isomorphe à une $G$-forme trace.}
       
       \smallskip
       
       \n {\it Démonstration.} Choisissons une extension finie $k$ de $\F_2$ assez grande pour que $U_G(k)$ contienne un sous-groupe cyclique $H$ non contenu dans $G.U_G(\kbar)$; c'est possible puisque $E_G \neq \{1\}$. Quitte à agrandir $k$, on peut aussi supposer que $\Gamma_k$ opère trivialement sur $E_G(\kbar)$.
       Soit $\varphi: \Gamma_k \to H$ un homomorphisme continu et surjectif; on peut voir $\phi$ comme un
       1-cocycle de $\Gamma_k$ à valeurs dans $U_G(\kbar)$; notons $[\varphi]$ sa classe dans $H^1(k,U_G)$.
       L'image de $[\varphi]$ dans $H^1(k,E_G) = \Hom_{\rm cont}(\Gamma_k,E_G(\kbar))$ est non triviale. Il en résulte que $[\varphi]$ n'est pas contenue dans l'image de $H^1(k,G) \to H^1(k,U_G)$; si $(V,q)$ est la $G$-forme obtenue en tordant la forme unité par $[\varphi]$, cela montre que $(V,q)$ n'est pas isomorphe à une
       $G$-forme trace.
       
       \smallskip
       
       \n {\it Remarques.}
       
        i) Si l'on prend $G = \widehat{{ S}}_3$ (cf. \S5.5.6), la construction ci-dessus est possible sur $k=\F_2$,
        et la $G$-forme ainsi obtenue satisfait non seulement à (7.4.1) et (7.4.2), mais aussi à (7.4.4).
        
        ii) Lorsque $G$ est un 2-groupe, la proposition 7.5.3 admet la réciproque suivante:
        
        \smallskip
        
        \n {\bf Proposition 7.5.4} - {\it Si $G$ est un $2$-groupe, et si $E_G = \{1\}$,  toute $G$-forme satisfaisant à (7.4.1) et (7.4.2) est isomorphe à une $G$-forme trace.}
        
        \smallskip
        
        \n {\it Démonstration.}  D'après le théorème 7.5.2 et la remarque qui le précède, les $G$-formes satisfaisant à (7.4.1) et (7.4.2) correspondent aux éléments de $H^1(k,U_G)$. Celles qui sont isomorphes à une $G$-forme trace correspondent aux éléments de l'image de l'application naturelle $ i :H^1(k,G) \to H^1(k,U_G)$. Nous devons donc montrer que $i$ est surjective. On a un diagramme commutatif:
   
   \bigskip

\hspace{13mm} \ \ $H^1(k,G)   \hspace{6mm} \stackrel{i} \rightarrow  \quad H^1(k,U_G)$

\vspace{2mm}
(7.5.5) \hspace{8 mm}  \ $ \pi \downarrow  \ \hspace{25 mm} \downarrow \pi'$

\vspace{2mm}

\hspace{13mm} \ $H^1(k,G/G_0) \  \stackrel{j} \rightarrow \ \ H^1(k,U_G/U^0_G).$

\medskip

De plus:

\smallskip

\n (7.5.6) {\it L'application \ $\pi: H^1(k,G) \to H^1(k,G/G_0)$ \ est surjective.} Cela résulte de ce que
$cd_2(\Gamma_k) \leqslant 1$, cf. [Se 64, II, \S2.2 et I, \S3.1].

\smallskip

\n (7.5.7) {\it L'application \ $ j:H^1(k,G/G_0) \to H^1(k,U_G/U_G^0)$ \ est bijective.} Cela résulte de ce que
$G/G_0 \to U_G/U_G^0$ est un isomorphisme, puisque $E_G =  \{1\}$.

\smallskip

\n (7.5.8) {\it L'application \ $\pi': H^1(k,U_G) \to H^1(k,U_G/U^0_G)$ est injective.} En effet, comme $G$
est un 2-groupe, les groupes $U_G$ et $U^0_G$ sont unipotents. De plus, $U^0_G$ est {\it  déployé} (``split'',
i.e. extension successive de groupes isomorphes à $\G_a$, cf. [Sp 98, \S12.3.5]), car il provient par extension des scalaires d'un groupe défini sur un corps parfait, à savoir $\F_2$. Cela entraîne que ses tordus (au sens de la cohomologie galoisienne) sont déployés, car la propriété d'être déployé se teste sur la clôture séparable du corps de base, cf. [Sp 98, th.14.3.8 (iii)]; leur $H^1$ est donc trivial; l'injectivité de $\pi'$ en résulte d'après la suite exacte de cohomologie non abélienne, cf. [Se 64, I, cor.2 à la prop.39].

\smallskip

  Ces trois propriétés, jointes au fait que $\pi' \circ i = j \circ \pi$ entraînent que $\pi'$ est bijectif et que $i$ est surjectif. D'où la proposition.
  
  \medskip
  
  \n (7.5.9) {\it Question}. {\it Existe-t-il un 2-groupe fini $G$ tel que} $E_G \neq  \{1\}$ ? 
     
         \bigskip

\n   7.6 - {\bf Traductions hermitiennes.}
   
    \medskip
    
      \n   7.6.1. {\it La $G$-forme associée à un élément hermitien.}

   \smallskip
    Soit $h$ un élément hermitien de $A$. La forme bilinéaire $q_h$ sur l'espace vectoriel $A$ définie par
    
     \smallskip
    (7.6.2)  $q_h(a,b) = \delta(ahb^*), \ \  a,b \in A,$

     \smallskip
    
  \n  est $G$-invariante et symétrique [rappelons, cf. \S7.1, que, si $x \in A$, $\delta(x)$ désigne le coefficient de 1 dans $x$]; elle est non dégénérée si et seulement si $h$ est inversible; dans ce cas, $(A,q_h)$
    est une $G$-forme.
    
     \smallskip
    
              Si $h'$ est un autre élément hermitien  inversible de $A$, les $G$-formes  $(A,q_h)$ et $(A,q_{h'})$ sont isomorphes
      si et seulement si il existe $a\in A^\times$ tel que  $h' = aha^*$; on dit alors que 
      $h$ et $h'$ sont {\it équivalents}, ce que nous écrirons  $h \sim h'$. 
      
      \smallskip
      
   \n   7.6.3. {\it Passage des $G$-formes aux éléments hermitiens.} Soit $(V,q)$ une $G$-forme. Choisissons 
      une base  $\{v\}$  du $k[G]$-module $V$; on associe à ces données l'élément hermitien
      $h=h_v$ défini par:
      
      \medskip
      
      (7.6.4)  $h = \sum_{g \in G} q(v,gv) g$, \ i.e. $\delta_g(h) = q(v,gv) $ pour tout $g\in G$.
      
      \medskip 
      
        L'application  : $A \to V$  donnée par  $a \mapsto av$ est un $G$-isomorphisme de $A$ sur $V$
        qui transforme $q_h$ en $q$; cela résulte de la formule (7.6.2) appliquée à $a=1, b=g$. Ainsi, 
        {\it toute $G$-forme est isomorphe à la forme $(A,q_h)$ associée à un élément hermitien inversible de $A$, bien déterminé à équivalence près.}
        
        \medskip

       Soit $h$ un élément hermitien inversible de $A$. Nous allons voir comment se traduisent pour la $G$-forme $(A,q_h)$  les propriétés (7.4.1),
        (7.4.2) et (7.4.4) .    
        
          \medskip
   
\n  7.6.5. {\it Traduction de} (7.4.1).

\smallskip

   \n {\bf Proposition 7.6.6} - {\it  Soit  $h$ un élément hermitien de $A$. Pour que $(A,q_h)$ satisfasse à $(7.4.1)$, il faut et il suffit que l'on ait}:   
   
   \smallskip
   (7.6.7) \  $\delta_s(h) = 0$  \ {\it pour tout $s\in G$ d'ordre} 2.
   
   \medskip
   
   \n {\it Démonstration}. Soit $s$ un élément de $G$ d'ordre 2. On a $\delta_s(h)= q_h(s.1,1)$; cela montre que (7.4.1) entraîne (7.6.7). Inversement, supposons que l'on ait $\delta_{s'}(h) = 0$ pour tout conjugué  $s'$ de $s$, et montrons que l'on a $q_h(sx,x) = 0$ pour tout  $x \in A$. La fonction  $f : x \mapsto q_h(x,sx)$ est une forme quadratique sur l'espace vectoriel $A$. Sa forme bilinéaire associée est
     $(x,y) \mapsto q_h(sx,y)+q_h(x,sy)$, qui est 0 puisque $s^2=1$. 
     Pour prouver que $f=0$ il suffit donc de prouver que $f$ s'annule sur la base $G$ de $A$. Or, si  $g$ est un élément de $G$, on a
   
   \smallskip  
   
   $q_h(g,sg) =  \delta(ghg^{-1}s)  = \delta(hg^{-1}sg) =  \delta(hs') = \delta_{s'}(h)$, \ avec $s'= g^{-1}sg$.
   
   \smallskip
   
 \n  Comme  $ \delta_{s'}(h) = 0$ par hypothèse, on en déduit que  $f(g)=0$ pour tout $g\in G$, d'où $f = 0$. 
  
  \medskip

  Un élément hermitien $h$ satisfaisant à la condition (7.6.7) sera dit {\it spécial}.  Si l'on utilise la partition  $G \sm G_2 = \Sigma\ \sqcup \ \Sigma^{-1}$ introduite dans la démonstration du théorème 5.1.2, cela revient à dire que $h$ s'écrit sous la forme   
  
  \medskip

  (7.6.8) $h = \lambda.1 + \sum_{\sigma \in \Sigma} h_{\sigma} (\sigma + \sigma^{-1})$, \ avec $\lambda, h_\sigma \in k$,
  
   \medskip

 \n  ce qui est équivalent à :
 
 \medskip
 
  (7.6.9) {\it Il existe  $\lambda \in k$ et $a \in A$ tels que}  $h = \lambda.1 + a + a^*$.
  
  \medskip
  
  \n Cela entraîne:
  
  \medskip

  (7.6.10)  $ \lambda = \t(h) = \delta(h)$,
   
   \medskip
   
   \n puisque $\delta(a) = \delta(a^*)$ et $\t(a) = \t(a^*). $ Nous dirons que $h$ est {\it normalisé} s'il est spécial et 
   si $\lambda = 1$, autrement dit si $h$ est de la forme $1+a+a^*$, avec $a\in A$.
   
   \medskip

\n 7.6.11.  {\it Traduction de } (7.4.2).
  
  \medskip

    \n {\bf Proposition 7.6.12} - {\it Supposons que $h$ satisfasse à la condition $(7.6.7)$. On a $:$
    
    \smallskip
    $(7.6.12)  \ \ q_h(x,x) = \t(h)\t(x)^2$  pour tout  $x \in A$.
    
    \smallskip
  \n Pour que
    $(A,q_h)$ satisfasse à $(7.4.2)$ il faut et il suffit que $\t(h)$ soit un carré dans $k$, et l'on a alors $:$
    
    \smallskip
    
    $(7.6.13) \ \  e =\t(h)^{-1/2}\s_G$,  \ où}  \ $\s_G = \sum_{g\in G} g$.
    
  \smallskip
 [Noter que $\t(h)$ est $\neq 0$ puisque $h$ est inversible et que $\t: A \to k$ est un homomorphisme
 d'algèbres.]

    \smallskip

\n {\it Démonstration.}   Les deux formes quadratiques

    \smallskip

\quad $x \mapsto q_h(x,x)$ \ \ et  \ \ $x \mapsto \t(h)\t(x)^2$

    \smallskip

\n sont additives: c'est clair pour la seconde, et, pour la première, cela résulte de ce que  $q_h(x,y)$ est symétrique. Pour montrer qu'elles sont égales, il suffit donc de le vérifier pour les éléments
$g$ de $G$. Or on a $q_h(g,g)= \delta(ghg^{-1}) = \delta(h)$ d'après (7.1.2) et (7.6.2),  et $\t(h)\t(g)^2 = \t(h) = \delta(h) $ d'après (7.6.10). Cela démontre (7.6.12). 

  Supposons maintenant que $k$ soit parfait, et définissons  $e\in A$ par la formule $e =\t(h)^{-1/2}\s_G$. On a:
  
  \smallskip
  
(7.6.14) \  $ q_h(e,x) = \delta(ehx^*) = \t(h)^{-1/2}\delta(\sigma_Ghx^*)$.

 \smallskip
\n Or on a:

 \smallskip
 
 (7.6.15) \ $\delta(\sigma_Ga) = \t(a)\sigma_G$  \ pour tout  $a \in A$.
 
  \smallskip
  
  \n On peut donc récrire (7.6.14) sous la forme:
  
  \smallskip

  (7.6.16) \   $ q_h(e,x) = \t(h)^{-1/2}\t(h)\t(x) = \t(h)^{1/2}\t(x)$,
  
  \smallskip
  
  \n d'où : 
  
  \smallskip
  
  (7.6.17) \ $q_h(e,x)e = \t(x)\sigma_G = \sigma_Gx$.

  \smallskip
  
  En comparant avec (7.4.2), on voit que $e$ n'est autre que ``l'élément canonique''  de $(A,q_h)$, ce qui démontre (7.6.13) lorsque $k$ est parfait.
  
    Dans le cas général, on applique ce qui précède à une extension parfaite $k'$ de $k$. On en déduit que
   l'élément canonique $e$ de $ k' \otimes_k A$ appartient à $A$ si et seulement si  $\t(h)^{1/2}$ appartient à $k$,
   i.e. si $\t(h)$ est un carré dans $k$, et l'on a alors (7.6.13), ce qui achève la démonstration de la proposition
   7.6.12.
   
   \medskip
   
   \n {\it Remarque.}  Si $h$ est normalisé (i.e. de la forme $1 + a + a^*$), la proposition précédente montre que la condition (7.4.2) est satisfaite, et que l'on a $e = \sigma_G$. Inversement, si (7.4.1) et (7.4.2) sont satisfaites, on peut écrire $h$ sous la forme $\mu h'$, avec $\mu = \t(h)^{1/2}$, et $h'$ normalisé;  noter que $h' \simeq h$. Ainsi, {\it toute $G$-forme satisfaisant à 
   $(7.4.1)$ et $(7.4.2)$ est isomorphe à une forme $(A,q_h)$, avec  $h$  normalisé}. En outre, {\it deux telles
   formes $(A,q_{h_1})$ et $(A,q_{h_2})$ sont isomorphes si et seulement si il existe $a\in A^\times$ tel que 
   $\t(a)= 1$ et} $h_2 = ah_1a^*$.
   
   \medskip

\n 7.6.18. {\it Traduction de } (7.4.4).

\medskip 
  Notons $F : A \to A$ l'endomorphisme de Frobenius de $A$:
  
  \medskip
  
  (7.6.19) \ \  $F(\sum x_gg) = \sum x_g^2 g$.
  
  \medskip
  
  \n {\bf Proposition 7.6.20} - {\it Pour que $(A,q_h)$ satisfasse à $(7.4.4)$ il faut et il suffit que $F(h) \sim h$.}
  
  \smallskip
  
  En effet, la condition $F(h) \sim h$ n'est qu'une traduction de  $(7.4.4').$
  
  \bigskip
  
  \n 7.7. {\bf Démonstration du théorème 7.5.2.}
  
  \medskip
  
    Dans ce \S, $k$ est supposé algébriquement clos, et nous en profiterons pour identifier un schéma
    lisse à l'ensemble de ses points.
    
    \smallskip
    
    Il s'agit de donner une réponse positive à la question (7.5.1) dans chacun des trois cas suivants:
    
     (i) $G$ est commutatif;
     
     (ii) $G$ est un 2-groupe;
     
     (iii) $|G|$ est impair.
     \smallskip
  
\n  Vu les traductions hermitiennes du \S7.6, cela revient à prouver que tout  hermitien $h$, qui est spécial et inversible, est tel que $h \sim 1$, autrement dit est de la forme $aa^*$, avec $a \in A^\times$.
  
  \medskip
  Notons $H_s$ l'ensemble des hermitiens spéciaux, autrement dit de la forme $\lambda.1 + a + a^*$, avec $\lambda \in k$ et  $a \in A$; c'est un sous-espace vectoriel de $A$, et d'après (7.6.8), on a:
  
  \smallskip
  
  (7.7.1)  $\dim H_s = 1 + |\Sigma| = 1 + \frac{1}{2} (|G| - |G_2|)$, \ \ avec les notations du \S7.6.5.
  
  \medskip
  
 Les éléments inversibles de $H_s$ forment un ouvert dense de $H_s$, que nous noterons $H_s^\times$.
 Faisons opérer le groupe algébrique $A^\times$ sur $A$ par  $a\bullet b = aba^*$. C'est une action linéaire;
 elle laisse stable  $H_s$  et $H_s^\times$:  cela se voit, soit par calcul direct, soit en invoquant la proposition 7.6.6. Soit $H'_s$ l'orbite de 1 par cette action; c'est l'ensemble des hermitiens de la forme $a\bullet1 = aa^*$.
 La question (7.5.1) a une réponse positive si et seulement si l'on a $H'_s = H_s^\times$, autrement dit si
 $A^\times$ agit transitivement sur $H_s^\times$.
 
 \medskip
   
   \n {\bf Proposition 7.7.3} - {\it L'orbite $H'_s$ est ouverte} (et donc dense) {\it dans} $H_s^\times$.
              
              \smallskip
              
              \n {\it Démonstration}. Le fixateur du point 1 est le groupe $U_G^{\sch}$. On a donc
              
              \smallskip
              
              $\dim H'_s = \dim A^\times - \dim U_G = |G| - \frac{1}{2} (|G| + |G_2|) +1 =  1 + \frac{1}{2} (|G| - |G_2|)$,
              
              \smallskip
              
              \n d'après la proposition 5.1.3. Vu (7.7.1), cela montre que $\dim H'_s = \dim H_s^\times$, d'où la proposition.

                     \medskip
                     
                        \n {\bf Corollaire 7.7.4} - {\it Pour que (7.5.1) ait une réponse positive, il faut et il suffit que
                        $H'_s$ soit fermé dans} $H_s^\times$.

                     \smallskip
                     
                     C'est clair.
                     
                     \medskip
                     
                      \n {\bf Corollaire 7.7.5} - {\it Tout élément de $H_s^\times$ qui appartient au centre de $A$
                      appartient à $H'_s$.}
                     
                     \smallskip
                     
                       \n {\it Démonstration}. Soit $h$ un tel élément. Son fixateur dans $A^\times$ est l'ensemble
                       des $a$ tels que $a \bullet h = h$, i.e. $aha^*=h$, ou encore $aa^*=1$ puisque $h$ et $a$ commutent;  c'est donc le groupe $U_G$. L'orbite de $h$ a donc la même dimension que celle de 1; cela entraîne que c'est une orbite ouverte, ce qui n'est possible que si elle est égale à  $H'_s$.
                       
                       \smallskip
                         \n {\bf Corollaire 7.7.6} - {\it La réponse à (7.5.1) est `` oui '' dans le cas} (i) (autrement dit, quand $G$ est commutatif).

                   \smallskip
                  
                  Cela résulte du corollaire précédent, qui donne même un résultat légèrement plus fort: il suffit que les hermitiens spéciaux de $A$ soient contenus dans le centre de $A$.
                  
                  \medskip
                  
                  \n 7.7.7. {\it Démonstration du théorème 7.5.2 dans le cas} (ii).
                  
                  \smallskip
                  
                  Ici, $G$ est un 2-groupe. L'algèbre $A$ est alors une algèbre locale (non commutative en général)
                  dont l'idéal bilatère maximal ${\frak m}$ est le noyau de $\t: A \to k$. Le noyau $A_1$ de $A^\times \to k^\times $ est un groupe unipotent connexe: cela se voit, par exemple, en filtrant $A$ par les puissances de ${\frak m}$ [cela résulte aussi de ce que la variété sous-jacente est isomorphe à un espace affine]. Le groupe $A_1$
                  opère sur la variété $H_1$ des éléments hermitiens normalisés au sens du \S7.6.5. Le fixateur
                  de 1 dans cette action est $U_G$. L'orbite de 1 est donc de dimension $\dim A_1 - \dim U_G$,
                  ce qui est aussi la dimension de $H_1$. Cette orbite est donc ouverte dans $H_1$. Mais les orbites d'un groupe unipotent agissant sur une variété affine (ou même seulement quasi-affine) sont fermées,
                  d'après un théorème de Rosenlicht ([Ro 61, th.2], voir aussi [Bo 91, prop.4.10]). On en conclut que $A_1$ opère transitivement sur $H_1$, autrement dit que tout hermitien normalisé $h$ est de la forme $aa^*$ avec $a  \in A_1$. Mais tout hermitien spécial inversible est un multiple d'un hermitien normalisé, et lui est donc équivalent puisque tout élément de $k$ est un carré; il est donc de la forme $aa^*$, avec $a\in A^\times$,
                  ce qui démontre le théorème 7.5.2 dans le cas (ii).
                  
                  \medskip

                  \n 7.7.8. {\it Démonstration du théorème 7.5.2 dans le cas} (iii).
                  
                  \smallskip

                  Ici, $G$ est d'ordre impair. L'algèbre $A$ est semi-simple. D'après le \S5.5.1, on peut  la décomposer en $A = A_0 \times \prod A_i$ avec:

                  \smallskip

                    $A_0 = k$,

                  \smallskip

                   $A_i = \M_{n_i} \times \M_{n_i}^{\rm opp}$, où l'involution est $(x,y)\mapsto (y,x)$.

                   \smallskip
                   
                   Or, dans chacune de ces algèbres à involution, tout élément hermitien $h$ s'écrit sous la forme  $h=aa^*$ :
                   c'est clair pour $A_0$, et, pour $A_i$, on a $h = (x,x)$ avec $x \in  \M_{n_i}(k)$ et l'on prend $a=(x,1)$. Le même résultat est donc vrai pour~$A$.

                                   \medskip
                                   
                                   \n {\bf 7.8. L'invariant d'une $G$-forme associé à un caractère essentiel.}
                                   
                                   \smallskip
                                   
                         Soit $\e : G \to \{1,c\}$ un caractère essentiel de $G$, cf. \S5.2, distinct du caractère unité. Nous allons voir comment on peut utiliser
                         $\e$ pour définir un {\it invariant} d'une $G$-forme, à valeurs dans $k/\wp(k)$.
                         
                         \smallskip
                          
                          \n 7.8.1. {\it Définition de $h_\e$.} Comme au \S5.2, notons $G_1$ (resp. $G_c$) l'ensemble des $g \in G$ tels que $\e(g)=1$ (resp. $\e(g)=c$), et choisissons une partie $S$ de $G_c$ telle que $G_c = S \sqcup S^{-1}$.
                          Si $h$ un élément hermitien de $A$, nous définirons $h_\e \in k$ par la formule:
                          
                          \smallskip
                          
             (7.8.2)             $h_\e = \sum_{s \in S} \delta_s(h)$.
             
             \smallskip
             
           Du fait que $h$ est hermitien, cette somme ne dépend pas du choix de  $S$.
           
           \medskip
           
          \n 7.8.3. {\it Une formule de transformation}.
          
          \smallskip

       Si $a = \sum_{g\in G} a_g g$ est un élément de $A$, nous poserons:
       
       \smallskip
         
         (7.8.4) \ $a_+ = \sum_{g\in G_1} a_gg$  \ et    \  $a_- = \sum_{g\in G_c}a_gg$;
         
         \smallskip
         
       \n  on a  :
         
         \smallskip
         (7.8.5) \ $a = a_+ + a_-$.
         
         \medskip
         
         \n {\bf Proposition 7.8.6} - {\it Si $h$ est un élément hermitien de $A$, on a}:
         
         \smallskip
         
          (7.8.7)  \ $(aha^*)_\e = \t(a)^2h_\e + \t(h)\t(a_+)\t(a_-)$ \ \ {\it pour tout} $a \in A$.
          
          \smallskip
          
          \n {\it Démonstration.} Nous allons employer une méthode qui permet d'interpréter $h_\e$
          comme une ``demi-somme'', ce qui évite d'avoir à utiliser l'ensemble auxiliaire $S$. Pour cela, choisissons
          un anneau commutatif $\tilde{k}$ muni d'une surjection $\tilde{k} \to k $, et tel que
          
            \smallskip

          (7.8.8) \ $4.\tilde{k} = 0$;  
          
          \smallskip

          (7.8.9) le noyau de $\tilde{k} \to k$ est $2.\tilde{k}$;

            \smallskip

          (7.8.9) l'application $x \mapsto 2x$ définit par passage au quotient une injection $\iota:k
          \to 2.\tilde{k}$.
                                        
                               \smallskip
                               Un tel anneau existe: il suffit par exemple de prendre pour  $\tilde{k}$ l'anneau des vecteurs de Witt de longueur  2  sur  $k$. On peut même demander que l'injection de (7.8.9) soit
                               une bijection, autrement dit que  $\tilde{k}$ soit un $\Z/4\Z$-module libre, cf. [AC IX, App., prop.4].

                                    \smallskip
                                    
                                    Posons $ \widetilde{A} = \tilde{k}[G]$, et choisissons des éléments $\tilde{h}$ et $\tilde{a}$ de $ \widetilde{A}$ dont les images dans $A$ sont $h$ et $a$, avec  $\tilde{h}$ hermitien.  Utilisons pour  $ \widetilde{A}$ les mêmes notations $\t, \tilde{a}_+, \tilde{a}_-$ que pour $A$. On a:
                                    
                                    \smallskip
                                    
         (7.8.10)                           $\iota(h_\e) = \t(\tilde{h}_-)$  \  et  \  $\iota((\tilde{a}h\tilde{a}^*)_\e) = \t((\tilde{a}\tilde{h}\tilde{a}^*)_-)$,
                                    
                                    \smallskip
\n ce qui montre que $h_\e$ et $(aha^*)_\e$ sont des ``demi-sommes'' dans $\tilde{k}$. De plus:
                                   
    \smallskip
    
        $  (\tilde{a}\tilde{h}\tilde{a}^*)_- \ =    \  \tilde{a}_-\tilde{h}_-\tilde{a}^*_- \ + \ \tilde{a}_+\tilde{h}_-\tilde{a}^*_+  \ + \ \tilde{a}_+\tilde{h}_+\tilde{a}^*_{-} \ + \ \tilde{a}_-\tilde{h}_+\tilde{a}^*_ +,$ 
    
     \n  d'où:
      \smallskip
      
      (7.8.11)  $\t((\tilde{a}\tilde{h}\tilde{a}^*)_-) =  \  \t(\tilde{h}_-)(\t(\tilde{a}_+)^2 + \t(\tilde{a}_-)^2) \ + \  2\t(\tilde{h}_+)\t(\tilde{a}_+)\t(\tilde{a}_-)$.
      
      \smallskip
      
    \n  Comme $\tilde{h}_-$ est hermitien, $\t(\tilde{h}_-)$ est divisible par 2; de plus, on a :
      
      \smallskip

 \quad      $\t(\tilde{a}_+)^2 + \t(\tilde{a}_-)^2 \equiv \t(\tilde{a})^2  \pmod 2$.

      \smallskip 
      
       \n On peut donc récrire (7.8.11) sous la forme   
      
      \smallskip

      (7.8.12) \  $\t((\tilde{a}\tilde{h}\tilde{a}^*)_-) =  \  \t(\tilde{h}_-)\t(\tilde{a})^2  + 2\t(\tilde{h}_+)\t(\tilde{a}_+)\t(\tilde{a}_-)$.

      \smallskip
      
   \n  On en déduit, d'après (7.8.10):
     
     \smallskip
     
     (7.8.13) \ $\iota((aha^*)_\e) =  \  \iota(h_\e\t(a)^2 + \t(h_+)\t(a_+)\t(a_-))  $,
     
     \smallskip
     
     \n d'où, d'après (7.8.9):
     
     \smallskip
     
     (7.8.14)  \    $(aha^*)_\e =  \  \   h_\e\t(a)^2 \ + \      \t(h_+)\t(a_+)\t(a_-)$,

     \smallskip
     
     \n et comme $\t(h_+) = \t(h)$, cela donne  (7.8.7).    
     
     \smallskip
     
     \n 7.8.15. {\it Définition de l'invariant}.
     
     \smallskip
     
       Supposons maintenant que $h$ soit inversible. On a alors $\t(h) \neq 0$, ce qui permet de définir l'invariant normalisé :
       
       \smallskip
       
    (7.8.16) \ $   h_\e^{\rm norm} = h_\e/\t(h)$.
       
           \smallskip
       
       \n {\bf Proposition 7.8.16} - {\it Soient $h$ et $h'$ deux éléments hermitiens inversibles de~$A$. Si
       $h \sim h'$, on a}:
       \smallskip

       (7.8.17) \  $h^{\rm norm}_\e \equiv  h'^{\rm norm}_\e \ \  {\rm mod} \ \wp(k).$
                                      
                                   \smallskip

       \n {\it Démonstration.}   Par hypothèse, il existe $a\in A^\times$ tel que $h'=aha^*$; on a alors
       $\t(h') = \t(h)\t(a)^2$, et la proposition 7.8.6 montre que:
       
       \smallskip
       
        (7.8.18)  \     $h'^{\rm norm}_\e  = h^{\rm norm}_\e  + \t(a_+)\t(a_-)/\t(a)^2.$    
        
        \smallskip
        
        \n Comme $\t(a)= \t(a_+)+\t(a_-)$, le terme  $\t(a_+)\t(a_-)/\t(a)^2$ peut s'écrire $\lambda + \lambda^2$,
        avec $\lambda = \t(a_+)/\t(a)$. D'où la proposition.
        
        \medskip
        
          Vu le dictionnaire entre classes d'hermitiens et $G$-formes, on voit que {\it l'on a  attaché à toute $G$-forme  $(V,q)$ un élément $(V,q)_\e$ de} $k/\wp(k)$, à savoir la classe mod  $ \wp(k)$ de $ h^{\rm norm}_\e $, où $h$ est l'un quelconque des hermitiens associés à $(V,q)$.      
          
          \medskip
          
          \small{\it Variante.} Nous venons de définir l'invariant $(V,q)_\e$ dans $k/\wp(k)$ en utilisant le dictionnaire hermitien. On peut aussi procéder de façon plus directe:
          
            Le sous-espace $V_1$ de $V$ fixé par $G_1= \Ker \e$ est de dimension 2. Si $x \in V_1$, il existe $y \in V$ tel que $x = \sum_{g\in G_1} gy$, et un simple calcul montre que la somme $\sum_{s\in S} q(y,sy)$ ne dépend, ni du choix de $S$, ni du choix de $y$; si on la note $q_\e(x)$, on peut prouver que $q_\e$ est une forme quadratique non dégénérée sur $V_1$, et que son invariant d'Arf dans $k/\wp(k)$ est égal à  $(V,q)_\e$:
            c'est la définition directe mentionnée plus haut.}

          \medskip
          
          \n {\it Remarque.} Lorsque $\e$ est le caractère unité, on convient que $h_\e=0$. Avec cette convention,
          on a:
          
          \smallskip
          
          \n {\bf Proposition 7.8.19} - {\it Soient $\e_1$ et $\e_2$ deux caractères essentiels. Alors $:$
          
          \smallskip
            
             $(7.8.20)$ \ \  $h_{\e_1\e_2} = h_{\e_1} +h_{\e_2}$ \ \ et} \ \ $h^{\rm norm}_{\e_1\e_2} = h^{\rm norm}_{\e_1} +h^{\rm norm}_{\e_2}$.
     
              \smallskip
         
         \n {\it Démonstration.} Soit $\e_3=\e_1\e_2$. Le cas où l'un des $\e_i$ est égal à 1 est immédiat. Supposons donc $\e_i \neq 1$ pour tout $i$. Notons $H_1$ l'ensemble des $g\in G$ tels que
    $\e_1(g)=1$ et $\e_2(g)=\e_3(g)\neq 1$, et définissons de même $H_2$ et $H_3$ (cf. démonstration du théorème 5.2.7). Décomposons chaque $H_i$ en $H_i = S_i \sqcup S_i^{-1}$. L'ensemble des $g\in G$ tels que $\e_1(g)\neq 1$ est $S_2 \sqcup S_2^{-1} \sqcup S_3 \sqcup S_3^{-1}$. Si l'on pose $x_i = \sum_{g\in S_i} \delta_g(h)$, on a
    
    \smallskip 
    
    (7.8.21) \  $h_{\e_1} = x_2+x_3$,  \ et de même \  $h_{\e_2} = x_3+x_1$  \ et \ $h_{\e_3} = x_1+x_2$.
    
    \smallskip
    
    En ajoutant ces trois égalités, on obtient $ h_{\e_1} + h_{\e_2} + h_{\e_3} = 0$, ce qui équivaut à la première formule de
    (7.8.20); la seconde en résulte en divisant par $\t(h)$.
    
    \smallskip
    
        \n {\bf Corollaire 7.8.22} - {\it L'invariant $(V,q) \mapsto (V,q)_\e$ est un homomorphisme de $X_G$ dans $k/\wp(k)$.}
        
  \n      [Rappelons que $X_G$ désigne le groupe des caractères essentiels de $G$, cf. \S5.3.]

                                   \medskip
                                   
                                   \n 7.8.23. {\it Le cas des $G$-formes trace.}
                                   
           Supposons maintenant que $(V,q)$ soit une $G$-forme trace $(L,q_L)$. Dans ce cas, il y a une façon  plus directe d'associer à $(V,q)$ et à un caractère essentiel $\e \neq 1$ un élément de $k/\wp(k)$: en effet, la sous-algèbre $L^\e$ de $L$ formée des éléments fixés par $G_1 = \Ker \e$ est une $k$-algèbre étale de degré 2, donc est isomorphe
           à  $k[x]/(x^2+x+z)$, où $z$ est bien déterminé mod $\wp(k)$; l'image de $z$ dans $k/\wp(k)$ est un invariant de $(L,q_L)$. En fait:
           
           \smallskip
           
           \n {\bf Proposition 7.8.24} - {\it Les deux invariants définis ci-dessus coïncident.}
           
           \smallskip
           
           \n {\it Démonstration}. Soit $\{v\}$ une base du $A$-module $L$, et soit $h$ l'élément hermitien correspondant. Quitte à multiplier $v$ par un scalaire, on peut supposer que $\Tr_{L/k}(v)=1$, auquel cas $h$
           est un hermitien normalisé: on a $\t(h)=1$. Posons:
           
            \smallskip
            
           (7.8.25)\ $y = \sum_{g\in G_1} gv$ \ et \ $z= \sum_{g\in G_c} gv$.

            \smallskip
            
        \n    On a:
        
               \smallskip

            (7.8.26) \  $y,z \in L^\e$,  \  $gy=z$  et $gz=y$  si $g\in G_c$,  \ et  \ $y+z=\Tr_{L/k}(v)=1$.
            
                   \smallskip

\n Le produit $yz$ est un élément de $k$. De façon plus précise:

       \smallskip

     \n {\bf Lemme 7.8.27} - {\it On a $yz = h_\e$.}
     
            \smallskip
            
            \n {\it Démonstration du lemme 7.8.27}. On a:
            
                   \smallskip
                   
                   (7.8.28)  \  $yz = \sum_{a\in G_1,b\in G_c} av.bv.$
                   
                          \smallskip
     
     Soit $\cal{M}$ l'ensemble des parties à deux éléments de $G$ rencontrant à la fois $G_1$ et $G_c$.
     Si $M=\{a,b\}$ est un élément de  $\cal{M}$, posons $v_M = av.bv$, de sorte que (7.8.28) peut se récrire sous la forme
     
       \smallskip
       
       (7.8.29) \ $yz = \sum_{M\in \cal{M}} v_M.$

       \smallskip
     
     Le groupe $G$ opère par translations à gauche sur l'ensemble $\cal{M}$. Cette action est libre, et chaque orbite contient un point et un seul de la forme $\{1,s\}$ avec $s\in S$. Cela permet
     de récrire (7.8.29) comme:
     
      \smallskip

\n (7.8.30) \ $yz = \sum_{g\in G, s\in S}  gv.gsv = \sum_{s\in S}\Tr_{L/k}( v.sv) = \sum_{s\in S} q_L(v,sv) = h_\e$,

\smallskip

\n ce qui démontre (7.8.28).

\smallskip

\n {\it Fin de la démonstration de la proposition 7.8.24}. Les formules (7.8.26) et (7.8.28) donnent la structure
de l'algèbre quadratique $L^\e$: elles montrent que cette algèbre est isomorphe à $k[x]/(x^2+x+h_\e)$; son invariant d'Artin-Schreier dans $k/\wp(k)$ est donc la classe de $h_\e$, c'est-à-dire $(V,q)_\e$.

\medskip

         \n {\bf 7.9. Application: le cas des $2$-groupes.}
         
         \smallskip
         
          \n {\bf Proposition 7.9.1} - {\it Supposons que $G$ soit un $2$-groupe tel que $E_G= \{1\}$. Soient
          $h$ et $h'$ deux hermitiens normalisés de $A$. Alors}:
          
          \smallskip
          
 \n  \quad   \ $h \sim h' \  \Longleftrightarrow \  h_\e  \equiv  h'_\e \ \  {\rm mod} \ \wp(k) $ \ {\it pour tout caractère essentiel $\e$ de} $G$.

\smallskip

\n {\it Démonstration.} D'après la proposition 7.5.4, les $G$-formes $(A,q_h)$ et $(A,q_{h'})$ sont isomorphes à des $G$-formes trace $(L,q_L)$ et $(L',q_{L'})$. Soient $\varphi_L, \varphi_{L'}: \Gamma_k \to G$ les homomorvarphismes (uniques à conjugaison près) associés à $L, L'$. D'après la proposition 6.6.17, on a
$(L,q_L) \simeq (L',q_{L'})$ si et seulement si les composés de $\varphi_L$ et $\varphi_{L'}$ avec $G \to G/G_0$ coïncident; cela revient à demander que, pour tout caractère essentiel $\e$ de $G$, on ait $\e \circ \varphi_L = \e \circ \varphi_{L'}$; d'après la proposition 7.8.24, cela revient à demander que $ h_\e  \equiv  h'_\e \ \  {\rm mod} \ \wp(k) $.

     \medskip
     
     \n 7.9.2. {\it Exemple $:$ le groupe quaternionien}. Supposons que $G$ soit le groupe quaternionien d'ordre 8.   Avec les
     notations du \S5.5.5, on peut écrire $h$ et $h'$ sous la forme:
     
     \smallskip

      $h = 1 + a_1(u+u^{-1}) + a_2(v+v^{-1}) + a_3(w+w^{-1})$ \ \ avec $a_i \in k$,

           \smallskip

           $h' = 1 + a_1'(u+u^{-1}) + a_2'(v+v^{-1}) + a_3'(w+w^{-1})$ \ \ avec $a'_i \in k$.
           
           \smallskip
           
   \n   Il y a      a trois caractères essentiels $\neq 1$; les $h_\e$ correspondants sont $a_1+a_2, \ a_2+a_3$ et $a_3+a_1$;
           ceux de $h'$ sont $a_1'+a'_2, \ a_2'+a_3'$ et $a_3'+a_1'$. Comme $E_G =  \{1\}$, cf. \S5.5.5, la proposition 7.9.1 s'applique; elle montre que:

           \smallskip
           
  (7.9.3)         \ $h \sim h' \  \Longleftrightarrow \ a_1+a_1' \equiv a_2+a_2' \equiv a_3+a_3' \ \ {\rm mod} \  \wp(k).$    
           
           \medskip
           
           \n 7.9.4. {\it Exemple $: \ G$ cyclique d'ordre $n=2^m, m \geqslant 2$}. Ici encore, on sait que  $E_G =  \{1\}$, cf. \S5.5.2. Soit  $s$ un générateur de $G$. Ecrivons $h$ et $h'$ sous la forme:
           
           \smallskip 
           
           $h = 1 + \sum_{0 < i < n/2} a_i(s^i + s^{-i})$ \ \ avec $a_i \in k$,

           \smallskip 
           
           $h' = 1 + \sum_{0 < i < n/2} a_i'(s^i + s^{-i})$ \ \ avec $a_i' \in k$.
           
           \smallskip
           
         \n    Soit $\e$ l'unique caractère essentiel non trivial de $G$. On a:
             
             \smallskip
             
 (7.9.5)       $h_\e = \sum_{i \ {\rm impair}} a_i$ \ et \     $h'_\e = \sum_{i \ {\rm impair}} a_i'$.

  \smallskip
           
    \n       La proposition 7.9.1 se traduit par:
           
           \smallskip
   (7.9.6)  \    $h \sim h' \  \Longleftrightarrow \    \sum_{i \ {\rm impair}} a_i \equiv \sum_{i \ {\rm impair}} a_i'     \ \ {\rm mod} \  \wp(k).$    
    
    \medskip
    
        \n {\bf 7.10. Application: description des $G$-formes trace en termes d'hermitiens.}
         
         \smallskip
       
         Etant donnés une $G$-algèbre galoisienne $L$ et un hermitien normalisé  $h = \sum  \lambda_g g$ de $A$, on aimerait savoir à quelle condition $h$ correspond à la $G$-forme trace  $q_L$, autrement dit dans quel cas on a:
         
           \smallskip
         (7.10.1)  $(L,q_L) \simeq_G (A,q_h)$,
         
           \smallskip
         
         \n ou, de façon équivalente:
         
           \smallskip
         
          (7.10.2) {\it Il existe une base $\{v\}$ du $A$-module $L$ telle que $q_L(v,g.v) = \lambda_g$ pour tout} $g\in G$.
          
            \smallskip
    Comme on l'a vu plus haut, il y a une condition nécessaire, imposée par les caractères essentiels. Rappelons-la:
        
     Soit $\e$ un caractère essentiel. Par la formule (7.8.2) on associe à $h$  un élément $h_\e$ de $k$,
    d'où un élément, noté $h(\e)$, de $k/\wp(k)$. D'autre part, on associe à $L$ l'élément de $k/\wp(k)$, noté $L(\e)$,
    qui correspond par Artin-Schreier au composé $\Gamma_k \stackrel{\varphi_L} \to G \stackrel{\e} \to \Z/2\Z$. La condition nécessaire en question est celle de la proposition 7.8.24, autrement dit:
    
    \smallskip
    
    (7.10.3) \ $h(\e) = L(\e)$ {\it dans $k/\wp(k)$ \ pour tout caractère essentiel $\e$ de} $G$.
    
    \smallskip
    
\n    Voici un cas où cette condition est suffisante:
    
    \smallskip
    
    \n {\bf Théorème 7.10.4} - {\it Supposons (7.10.3) satisfaite, ainsi que}:
    
    \smallskip
    
    (7.10.5) {\it Il existe un $2$-sous-groupe $H$ de $G$ tel que $E_H =  \{1\}$ et que $\lambda_g = 0$ pour tout $g \not \in H$.}
    
    \smallskip
    
    \n {\it Alors $h$ et $q_L$ se correspondent, autrement dit on a (7.10.1) et (7.10.2).}
    
    \smallskip
    
    \n {\it Démonstration}. Choisissons un sous-groupe $H$ de $G$ satisfaisant à (7.10.5). D'après la proposition 7.5.4
    il existe une $H$-algèbre galoisienne $M$ telle que $(M,q_M)$ soit $H$-isomorphe à $(k[H],h)$. Soit $L'$ la $G$-algèbre galoisienne induite de $M$. Si $\e$ est un caractère essentiel de $G$, on a $L'(\e) = h(\e) = L(\e)$. D'après le théorème 6.1.7 cela entraîne que les $G$-formes $q_L$ et $q_L'$ sont isomorphes; d'où (7.10.1).
    
    \medskip
    
    \n 7.10.6. {\it Application $:$ démonstration du théorème 7.3.2}. Reprenons les notations $(L,\gamma,z,...)$ de ce théorème, et soit $h = 1+ z(\gamma + \gamma^{-1})$; c'est un hermitien normalisé. Il nous faut montrer que $h$ et $q_L$ se correspondent. Cela résulte du théorème 7.10.4 appliqué au groupe cyclique $H$ engendré par $\gamma$: la condition (7.10.3) est vraie par construction, et la condition $E_H =  \{1\}$ est vraie d'après le \S5.5.2.
    
    \medskip
    
      \n 7.10.7. {\it Autre application.} Soit $S$ un $2$-sous-groupe de Sylow de $G$.  Faisons l'hypothèse:
      
      \smallskip
      
       (7.10.8)  $E_S =  \{1\}$ . 
       
       \smallskip
       
         Nous allons en déduire une famille d'éléments hermitiens de $A$ qui décrivent, de façon essentiellement unique, toutes les $G$-formes trace.
           
             Pour cela, choisissons une base $\xi_1,...,\xi_n$ du $\Z/2\Z$-espace vectoriel $G/G_0$, et soit $\e_1,...,\e_n$
             la base duale du groupe des caractères essentiels de  $G$. Pour chaque  $i$, choisissons un représentant  $\gamma_i$ de $\xi_i$ dans $S$. Considérons les hermitiens de la forme:
             
             \smallskip
             
             (7.10.9) \ \  $h(z_1,...,z_n) = 1 + \sum z_i(\gamma_i+ \gamma_i^{-1})$, \ avec \ $z_i \in k$.
             
    \smallskip
      Ce sont des hermitiens normalisés. Nous allons voir qu'ils suffisent:
      
      \smallskip
      
      \n    \n {\bf Théorème 7.10.10} - {\it Faisons l'hypothèse (7.10.8). Alors toute $G$-forme trace $q_L$ correspond à un hermitien  $h(z_1,...,z_n)$ tel que l'image de $z_i$ dans $k/\wp(k)$ soit égale à $L(\e_i)$.   }
      
      \smallskip
      
      \n {\it Démonstration.} Cela résulte du théorème 7.10.4, appliqué en prenant $H=S$.
      
      \medskip
      
      On obtient ainsi (hélas, sous l'hypothèse (7.10.8)) une description explicite des $G$-formes trace. D'où l'intérêt du calcul de $E_S$ lorsque $S$ est un $2$-groupe.

    \bigskip

\bigskip

\centerline{\bf Références}

\bigskip

{\small

 [A V]  N. Bourbaki, {\it Algèbre, Chapitre} V, {\it Corps commutatifs}, Masson, Paris, 1981 et Springer-Verlag, 2006; traduction anglaise, Springer-Verlag, 1998.

 [A VIII]   ------- , {\it Algèbre, Chapitre} VIII, {\it Anneaux et modules semi-simples}, nouvelle édition révisée, Springer-Verlag, 2011.
  
 [AC IX]  ------- , {\it Algèbre Commutative, Chapitre} IX, {\it Anneaux locaux noethériens complets}, Masson,
 Paris, 1983 et Springer-Verlag, 2006.
 
 [Arf 41] C. Arf, {\it Untersuchungen über quadratische Formen in Körpern der Charakteristik $2$}, J. Crelle {\bf 183} (1941), 148-167.

[Ba 89]  E. Bayer-Fluckiger, {\it Self-dual normal bases} I, Indag. Math. {\bf 51} (1989), 379-383.

[Bar 13] M. Barakat, {\it Computations of unitary groups in characteristic 2}, 
\url{http://www.mathematik.uni-kl.de/~barakat/__for_JP_Serre/UnitaryGroup.pdf}

[BFS 94] E. Bayer-Fluckiger \& J-P. Serre, {\it Torsions quadratiques et bases normales autoduales}, Amer. J. Math. {\bf116} (1994), 1-64.

[BL 90]  E. Bayer-Fluckiger \& H.W. Lenstra, Jr, {\it Forms in odd degree extensions and self-dual normal bases}, Amer. J. Math. {\bf112} (1990), 359-373.

[BLR 90] S. Bosch, W. Lütkebohmert \& M. Raynaud, {\it Néron Models}, Ergebn. Math. (3 Folge) {\bf21}, Springer-Verlag, 1990.

[Bo 91] A. Borel, {\it Linear Algebraic Groups}, second enlarged edition, Springer-Verlag, 1991.

[BoS 64] A. Borel \& J-P. Serre, {\it Théorèmes de finitude en cohomologie galoisienne}, Comm. Math. Helv.
{\bf 39} (1964), 111-164 (= A.Borel, Oe.64).

[BPS 13]  E. Bayer-Fluckiger, R. Parimala \& J-P. Serre, {\it Hasse principle for
$G$-trace forms}, Izvestiya RAS/Ser. Math. {\bf 77} (2013), 5-28. 

[BR 00]  V. Bovdi \& A.L. Rosa, {\it On the order of the unitary subgroup of a modular group algebra}, Comm. Algebra {\bf 28} (2000), 1897–1905. 

[BT 65] A. Borel \& J. Tits, {\it Groupes réductifs}, Publ. Math. IHES {\bf 27} (1965), 55-150 (= A. Borel, Oe.66 et  J. Tits, Oe.61).

[Bu 11] W. Burnside, {\it Theory of Groups of Finite Order}, second edition, Cambridge Univ. Press, 1911 et Dover Publ. 1955.

[CGP 10] B. Conrad, O. Gabber \& G. Prasad, {\it Pseudo-reductive Groups}, New math. monographs {\bf 17}, Cambridge Univ. Press, 2010.

[Ch 58] C. Chevalley, {\it Classification des groupes de Lie algébriques},
Séminaire ENS 1956-1958, Secrétariat math., IHP, 1958; édition révisée
par P. Cartier, {\it Classification des Groupes Algébriques Semi-simples},
Springer-Verlag, 2005.

[CR 62] C.W. Curtis \& I. Reiner, {\it Representation Theory of Finite Groups and Associative Algebras}, Pure and Applied Math. XI, Intersc. Publ. 1962.

     [EVT II] N. Bourbaki, {\it Espaces Vectoriels Topologiques, Chapitre} II{\it, Ensembles convexes et espaces localement convexes},
     Masson, Paris, 1981 et Springer-Verlag, 2007; traduction anglaise, Springer-Verlag, 1987.
     
     [Fe 82] W. Feit,  {\it The Representation Theory of Finite Groups}, North-Holland Math. Library {\bf 25}, 1982.
     
     [Fu 84] W. Fulton, {\it Intersection Theory}, Springer-Verlag, 1984.
     
     [GM 13] P. Gille \& L. Moret-Bailly, {\it Actions algébriques de groupes arithmétiques}, in  {\it Torsors, Étale Homotopy and Applications to Rational Points}, éditeur A. Skorobogatov,   L.M.S. Lect. Notes {\bf 405} (2013), 231-249.
     
     [Hu 67] B. Huppert. {\it Endliche Gruppen} I, Springer-Verlag, 1967.

     [Ja 03] J.C. Jantzen, {\it Representations of Algebraic Groups}, seconde édition, Math. Surveys {\bf 107}, AMS, 2003.
     
       [KMRT 98] M. Knus, A. Merkurjev, M. Rost \& J--P. Tignol, {\it The Book of Involutions}, AMS Colloquium Publications {\bf 44}, 1998. 
      
      [La 56] S. Lang, {\it Algebraic groups over finite fields}, Amer. J. Math. {\bf 78} (1956), 555-563.

     [MT 12] G. Malle \& D. Testerman, {\it Linear Algebraic Groups and Finite Groups of Lie Type}, Cambridge Studies in Advanced Mathematics {\bf 133}, Oxford, 2011.
     
     [Ro 61]  M. Rosenlicht, {\it On quotient varieties and the affine embedding of certain homogeneous spaces},
     Amer. J. Math.  {\bf 83}  (1961), 211-223.

[Se 62] J-P. Serre, {\it Cohomologie galoisienne des groupes algébriques linéaires}, Colloque sur la théorie des groupes algébriques, Bruxelles (1962), 53-68 (= Oe.53).

[Se 64] -------, {\it Cohomologie Galoisienne}, Springer Lect. Notes {\bf 5} (1964); cinquième édition révisée et complétée, Springer-Verlag, 1994; traduction anglaise, {\it Galois Cohomology}, Springer-Verlag, 1997.

[Se 68] -------, {\it Représentations linéaires des groupes finis}, Hermann, Paris, 1968; cinquième édition
corrigée et augmentée, Hermann, Paris, 1998; traduction anglaise, {\it Linear Representations of Finite Groups}, Springer-Verlag, 1977.

[Se 05] -------, {\it BL-bases and unitary groups in characteristic 2}, Oberwolfach Reports {\bf2} (2005), 37-40.

  [SGA 3] M. Demazure \& A. Grothendieck, {\it Schémas en Groupes}, 3 vol., Springer Lect. Notes {\bf 151, 152, 153} (1970); édition révisée par P. Gille et P. Polo, {\it Documents Mathématiques},  SMF, 2012, $201?$, 2012.
  
  [Sp 98] T.A. Springer, {\it Linear Algebraic Groups}, second edition, Birkhäuser Boston, 1998.

[St 67] R. Steinberg, {\it Lectures on Chevalley Groups}, Notes polycopiées,  Yale University, 1967.

[St 75]  -------, {\it Torsion in Reductive Groups}, Adv. in Math. {\bf15} (1975), 63-92 (= C.P. 415-444).

[Va 05] A. Vasiu, {\it Normal, unipotent subgroup schemes of reductive groups}, C. R. Acad. Sci. Paris {\bf 341} (2005), 79-84.

[Wa 79] W.C. Waterhouse, {\it Introduction to affine group schemes}, Springer-Verlag, 1979.

[We 61] A. Weil, {\it Adeles and Algebraic Groups}, I.A.S. Princeton, 1961; Progress in Math. {\bf 23}, Birkhäuser Boston, 1982.}
     
    \bigskip

     \bigskip
     
    \n Collège de France
    
    \n 3, rue d'Ulm
    
    \n 75005 PARIS

     \end{document}